\documentclass[11pt]{article}
\usepackage[english]{babel}
\usepackage[left=1in,right=1in,top=1in,bottom=1.5in]{geometry} 
\usepackage[skip=5pt, indent=10pt]{parskip}
\usepackage{amssymb, amsmath, amsthm}
\usepackage{hyperref}
\usepackage[dvipsnames]{xcolor}
\usepackage{colortbl}
\usepackage{nicefrac}
\usepackage[margin=1cm]{caption}
\usepackage{natbib}
\usepackage{bbm}
\usepackage[scr]{rsfso}
\usepackage{comment}
\usepackage{algorithm}
\usepackage{algpseudocode}
\usepackage{times}
\usepackage{graphicx}
\usepackage{stackrel}
\usepackage{cases}

\definecolor{darkred}{HTML}{880000}
\definecolor{darkblue}{HTML}{000088}

\hypersetup{
  colorlinks=true,
  linkcolor=RoyalBlue,
  citecolor=RoyalBlue,
  urlcolor=RoyalBlue
}

\newcommand{\myendproof}{$\hfill\square$}

\renewcommand{\leq}{\leqslant}

\renewcommand{\geq}{\geqslant}

\renewcommand{\tilde}{\widetilde}

\newcommand{\eps}{\varepsilon}

\newcommand{\dd}{{\mathrm{d}}}
\newcommand{\Var}{{\rm Var}}

\newcommand{\sfp}{\mathsf{p}}

\newcommand{\NN}{\mathsf{NN}}

\newcommand{\E}{\mathbb E}
\newcommand{\N}{\mathbb N}
\newcommand{\p}{\mathbb P}
\newcommand{\R}{\mathbb R}
\newcommand{\Z}{\mathbb Z}
\newcommand{\1}{\mathbbm 1}

\newcommand{\bk}{{\mathbf k}}

\newcommand{\cB}{{\mathcal{B}}}

\newcommand{\cF}{{\mathcal{F}}}
\newcommand{\cG}{{\mathcal{G}}}
\newcommand{\cH}{{\mathcal{H}}}
\newcommand{\cK}{{\mathcal{K}}}

\newcommand{\cL}{{\mathcal{L}}}

\newcommand{\cN}{{\mathcal{N}}}
\newcommand{\cO}{{\mathcal{O}}}

\newcommand{\cV}{{\mathcal{V}}}

\newcommand{\floor}[1]{{\lfloor #1 \rfloor}}

\def\argmin{\operatornamewithlimits{argmin}}

\newcommand{\gelu}{\mathrm{GELU}}

\newcommand{\sdata}{\sigma_{\mathrm{data}}}
\newcommand{\smin}{\sigma_{\mathrm{min}}}

\newcommand{\integral}[1]{\int\limits_{#1}}

\newtheorem{Th}{Theorem}[section]
\newtheorem{Lem}[Th]{Lemma}
\newtheorem{Def}[Th]{Definition}
\newtheorem{Prop}[Th]{Proposition}

\newtheorem{Rem}[Th]{Remark}

\newtheorem{As}[Th]{Assumption}

\title{Escaping the Curse of Dimensionality in One-Step Flow-Based Generative Models}

\author{
Konstantin Yakovlev\thanks{HSE University, Russian Federation. \texttt{kdyakovlev@hse.ru}}
}

\date{}

\begin{document}

\maketitle

\begin{abstract}
    Flow and diffusion-based generative models achieve state-of-the-art quality, but require multiple forward passes at inference time, making them computationally costly.
    Consistency models overcome this by learning a flow map directly, enabling few-step generation and dramatically reducing the computation needed at test time.
    However, formal statistical guarantees for such one-step generators are scarce.
    Existing results either suffer from the curse of dimensionality or rely on impractical loss functions.
    We investigate whether one-step flow-based models can escape the curse of dimensionality under the assumption that the data distribution is generated by a low-dimensional surface perturbed with Gaussian noise.
    In this setup, we derive approximation and generalization error bounds for the flow map estimate that circumvent the curse of dimensionality.
    Our findings provide the first principled explanation for practical success of one-step flow-based generative models.
\end{abstract}

\section{Introduction}

Generative models based on flow \citep{lipman2023flow,liu2023flow,albergo2025stochastic} and diffusion \citep{ho2020denoising,song2021a,karras2022elucidating} have demonstrated exceptional performance across various domains, including images \citep{esser2024scaling}, audio \citep{lu2022conditional}, and video \citep{ho2022video}.
However, these models require a large number of evaluations of the learned model, presenting a significant computational challenge.
A promising avenue is to estimate the flow map rather than the instantaneous velocity, which is accomplished within the framework of Consistency Models \citep{song2023consistency,song2024improved}.
This approach enables few-step generation, resulting in a reduction in sampling time by several orders of magnitude compared to traditional models.
The recently developed methods for learning flow maps \citep{kim2024consistency,boffi2025flow,boffi2025how,geng2025mean,lu2025simplifying,frans2025one,zhang2026alphaflow} effectively bridge the gap between few-step models and multi-step diffusion and flow-based models.

Despite the empirical success of these methods, theoretical guarantees regarding the statistical efficiency of distribution estimation using few-step generative models are scarce.
\cite{dou2024theory} made the first effort: under mild
assumptions they proved convergence rates for a one-step generator, but
the bounds suffer from the curse of dimensionality and rely on an impractical
Wasserstein-loss consistency criterion.
\cite{li2025consistency} obtained quantitative guarantees for distribution estimation of a one-step consistency model.
However, their iterative training framework
deviates from the framework of \cite{song2023consistency} and requires training $\cO(D^{5 / 2} / \eps)$ models from scratch, where $D$ is the ambient dimension and $\eps$ is the target accuracy.
This is highly impractical.
\cite{ding2024characteristic} studied statistical efficiency of flow-based one-step generative models and derived a convergence rate that suffers from the curse of dimensionality.
Under the assumption that the data distribution is supported on a low-dimensional linear subspace, convolved with Gaussian noise, the bound was substantially improved.
However, this assumption is not satisfied in practice.

Overall, the aforementioned works do not explain the success of flow map learning in practice.
Note that in many practical scenarios the data distribution has an underlying low-dimensional structure \citep{bengio2013representation,pope2021the,brown2023verifying}.
It has been recognized that
generative models which exploit such intrinsic structure can escape the
curse of dimensionality \citep{stephanovitch2024wasserstein,
tang2024adaptivity, yakovlev2025generalization,
chakraborty2025statistical}.
For these models, the exponent of $n$ in the convergence rate is determined by the intrinsic dimension.
This naturally raises the following research question:

\begin{quote}
    \centering
    \textit{
    Can one-step flow-based models escape the curse of dimensionality?
    }
\end{quote}

To answer this question, we consider the self-distillation framework introduced by \cite{boffi2025flow,boffi2025how}.
We focus on this framework because it admits a well-defined population objective and has proven effective for learning flow maps in practice.

\noindent
\textbf{Contributions.}\quad
In this paper, we affirmatively answer the above question.
Specifically, we investigate the statistical efficiency of self-distillation \citep{boffi2025how} under the noisy setting as in \cite{yakovlev2025generalization}.
Our main contributions are:

\begin{enumerate}
    \item We establish that a flow map can be approximated without the curse of dimensionality, an insight that extends beyond the scope of this work.
    \item We demonstrate that the flow map estimator obtained via self-distillation  achieves a convergence upper bound of $\tilde\cO(n^{-\frac{2\beta}{5\beta + 2d}})$, thereby  avoiding the curse of dimensionality.
      Here, $d$ denotes the intrinsic dimension, which can be substantially smaller than the ambient dimension $D$.
    \item Finally, we show that, for a wide
    class of stochastic interpolant coefficients, the learned flow map recovers the target distribution without the curse of dimensionality in terms of the 2-Wasserstein distance.
\end{enumerate}

\noindent
\textbf{Notation.}\quad
Throughout the paper, we use $f \lesssim g$ and $g \gtrsim f$ to denote $f = \cO(g)$.
If both $f \lesssim g$ and $g \lesssim f$ hold, then we write $f \asymp g$.
Furthermore, the notation $\tilde\cO(\cdot)$ is defined similarly to $\cO(\cdot)$ except that the logarithmic factors are hidden.
For any $R > 0$, the Euclidean ball in $\R^r$ centered at the origin is denoted by $\cB(0, R)$.
We denote the set of non-negative integers by $\mathbb{Z}_+$.
For a multi-index $\bk = (k_1, \dots, k_r) \in \Z_+^r$ and a vector $v = (v_1, \dots, v_r) \in \R^r$, we define $v^\bk = v_1^{k_1}v_2^{k_2}\dots v_r^{k_r}$ and $|\bk| = k_1 + k_2 + \dots + k_r$.
We denote the 2-Wasserstein distance between probability measures $\mu, \nu$ by $W_2(\mu, \nu)$ and the distribution of a random vector $X$ by $\mathrm{Law}(X)$.
Lastly, we utilize the operators $a \vee b = \max(a, b)$, $a \wedge b := \min(a, b)$, and $\floor{\beta}$ for the largest integer strictly less than $\beta > 0$.

\subsection{Other Related Works}

\noindent
\textbf{Statistical theory of diffusion and flow-based models.}\quad
The estimation of the velocity field and the score function (namely, the gradient of the log-density) is a central topic in recent generative modeling literature \citep{oko2023diffusion, wibisono2024optimal,zhang2024minimax, fukumizu2025flow}.
While the manifold hypothesis is commonly employed to avoid the curse of dimensionality \citep{tang2024adaptivity,azangulov2024convergence}, such approaches rely on restrictive geometric or density assumptions.
\cite{yakovlev2025generalization}
recently proposed a more general assumption that better suits real-world data.
However, existing theory focuses mainly on instantaneous velocity learning.
Our work extends this to the flow map learning, which is a more challenging problem.

\noindent
\textbf{Convergence of consistency models.}\quad
A separate line of research focuses on the iteration complexity of consistency models
\citep{lyu2024sampling, yang2025improved,chen2025convergence}, typically assuming that the flow map has been accurately estimated.
However, these works do not address the flow map estimation problem itself, leaving a critical gap in the theoretical understanding of statistical efficiency of such models.

\section{Preliminaries and notations}
\label{sec:prelim_not}

\noindent
\textbf{Norms.}\quad
For a vector $v$ and a matrix $A$, we denote their Euclidean and operator norms by $\|v\|$ and $\|A\|$, respectively.
The notations $\|\cdot \|_0$ and $\| \cdot \|_\infty$ refer to the number of non-zero entries and the maximum absolute entry value of a vector or matrix.
For a vector-valued function $f$ defined on a set $\Omega \subseteq \R^r$ and a non-negative weighting function $p : \Omega \to \R$, we denote the $L^\infty$ and weighted $L^2$ norms as
\begin{align*}
    \|f\|_{L^\infty(\Omega)} = \sup_{x \in \Omega}\|f(x)\|
    \quad \text{and} \quad
    \|f\|_{L^2(\sfp)} = \left\{\integral{\Omega}\|f(x)\|^2 \sfp(x) \, \dd x\right\}^{1/2} .
\end{align*}
Finally, for a random variable $\xi$, we denote its $\psi_1$-norm by $\|\xi\|_{\psi_1} = \inf\{t > 0  :  \E\exp(|\xi| / t) \leq 2\}$.

\noindent
\textbf{Stochastic Interpolants.}\quad
Let $X_1 \sim \sfp^*$, where $\sfp^*$ is the probability distribution on $\R^D$ we aim to sample from.
The \emph{stochastic interpolant} framework \citep{albergo2025stochastic} operates as follows.
The stochastic interpolant is a stochastic process $(X_t)_{0 \leq t \leq 1}$ defined by independently sampling $X_0 \sim \cN(0, I_D)$, $X_1 \sim \sfp^*$ and setting
\begin{align}
\label{eq:X_t_def}
    X_t = \alpha_t X_0 + \beta_t X_1.
\end{align}
Here, $\alpha_t$ and $\beta_t$ are continuous on $[0,1]$ and differentiable on $[0,1)$, satisfying the boundary conditions $\alpha_0 = \beta_1 = 1$ and $\alpha_1 = \beta_0 = 0$.
Let $\sfp_t^*$ be the distribution of $X_t$.
Under the regularity and integrability conditions of \cite[Theorem 6]{albergo2025stochastic}, the velocity field $b^*(y, t) = \E[\dot{X}_t \, | \, X_t = y]$ generates a deterministic flow between $\cN(0, I_D)$ and $\sfp^*$.
Specifically, for
\begin{align*}
    \dot{Z}_t = b^*(Z_t, t),
    \quad 0 \leq t \leq 1,
    \quad Z_0 \sim \cN(0, I_D)
\end{align*}
it holds that $Z_t$ and $X_t$ have the same distributions for every $t \in [0, 1]$.
The challenge is that the velocity field $b^*$ is unknown and should be estimated from data using flow matching objective.
Specifically, given i.i.d. samples $Y_1, Y_2, \dots, Y_n$ from $\sfp^*$, the velocity field estimate $\widehat{b}(y, t)$ is defined as the empirical risk minimizer
\begin{align}
\label{eq:fm_loss_def}
    \widehat{b}\in \argmin_{b \in \cB} \frac{1}{n}\sum_{i = 1}^n \ell_{FM}(b, Y_i),
    \quad \ell_{FM}(b, y) = \integral{0}^{1 - t_0} \E[\|b(X_t, t) - \dot{X}_t\|^2 \, | X_1 = y] \, \dd t ,
\end{align}
where $\cB$ is the hypothesis class and $t_0 \in [0, 1)$ is the early stopping time.
We use this objective only when $\int_0^{1-t_0}\E\|\dot X_t\|^2\,\dd t<\infty$. Early stopping is needed for some schedules, including the diffusion interpolant, but is not needed for the linear interpolant in our positive-noise setting. Its role in related estimation problems is discussed in \cite{fukumizu2025flow,oko2023diffusion,yakovlev2025generalization}.
The flow matching objective has the following property.
\begin{Prop}
\label{prop:fm_loss_dif}
    For every $b \in \cB$ such that the maximum of $\E\|b(X_t, t)\|^2$ over $t \in [0, 1 - t_0]$ is finite, the following holds:
    \begin{align}
    \E[\ell_{FM}(b, X_1) - \ell_{FM}(b^*, X_1)] = \integral{0}^{1 - t_0} \E \|b(X_t, t) - b^*(X_t, t)\|^2 \, \dd t .
\end{align}
\end{Prop}
The proof of Proposition \ref{prop:fm_loss_dif} is deferred to Appendix \ref{sec:prop_fm_loss_dif_proof}.

\medskip
\noindent
\textbf{Self-Distillation.}\quad
The self-distillation framework of \cite{boffi2025how} accelerates sampling by directly learning a flow map that transports a simple base distribution to the target distribution in a single step.
For an arbitrary $y \in \R^D$ and $0 \leq s \leq t \leq 1 - t_0$, consider the flow generated by the true velocity field $\dot Z_\tau = b^*(Z_\tau, \tau)$ with $s \leq \tau \leq t$ and $Z_s = y$.
For every $0 \leq s < t \leq 1 - t_0$ and $y \in \R^D$, we define the actual flow map
\begin{align}
\label{eq:v_star_def}
    \Psi^*_{s, t}(y) = y + \integral{s}^t b^*(Z_\tau, \tau) \, \dd \tau
    = y + (t - s) v_{s, t}^*(y),
    \quad v^*_{s, t}(y) = \frac{1}{t - s}\integral{s}^t b^*(Z_\tau, \tau) \, \dd \tau .
\end{align}
Following the notation in \cite{geng2025mean}, we refer to $v^*_{s, t}(y)$ as the \emph{average velocity}.
On the diagonal, define $\Psi^*_{t,t}(y)=y$ and $v^*_{t, t}(y) = b^*(y, t)$ by continuous extension, for all $y \in \R^D$ and $t \in [0, 1 - t_0]$.
Moreover, \cite{boffi2025how} demonstrate that $\Psi^*_{s, t}$ solves the Lagrangian equation
\begin{align}
\label{eq:psi_star_lagrange}
    \partial_t \Psi^*_{s, t}(y) = b^*(\Psi^*_{s, t}(y), t),
    \quad \text{for } 0\leq s<t\leq1-t_0 \text{ and } y\in\R^D .
\end{align}
Inspired by this property, we parametrize a candidate map as $\Psi_{s, t}(y) = y + (t - s)v_{s, t}(y)$ and define the self-distillation objective
\begin{align}
\label{eq:self_dist_loss_def}
    \ell_{SD}(v, y)
    = \ell_{FM}(v, y)
    + \integral{0}^{1 - t_0} \E \|\partial_t\Psi_{0, t}(X_0) - v_{t, t}(\Psi_{0, t}(X_0)) \|^2 \, \dd t .
\end{align}
Here, the first term is the flow matching loss defined in \eqref{eq:fm_loss_def}, where we identify the velocity candidate $b(y, t)$ with $v_{t, t}(y)$.
The second term, which we denote by $R(v)$, penalizes deviations from the Lagrangian consistency condition \eqref{eq:psi_star_lagrange}.
This penalty controls the discrepancy between the time derivative of the candidate trajectory and the instantaneous velocity.
Consequently, the loss function in \eqref{eq:self_dist_loss_def} is expressed as $\ell_{SD}(v, y) = \ell_{FM}(v, y) + R(v)$.
Moreover, since $v^*$ satisfies the Lagrangian equation \eqref{eq:psi_star_lagrange}, the consistency penalty vanishes, i.e., $R(v^*) = 0$.
Given i.i.d. data points $Y_1, Y_2, \dots, Y_n \sim \sfp^*$, the learned flow map $\widehat{\Psi}_{0, t}(y) = y + t\widehat{v}_{0, t}(y)$ is defined through empirical risk minimization
\begin{align}
\label{eq:avg_vel_erm_def}
    \widehat{v} \in \argmin_{v \in \cV} \left\{ \frac{1}{n}\sum_{i = 1}^n \ell_{SD}(v, Y_i) \right\} .
\end{align}
Here, $\cV$ denotes the set of average velocity estimators.
Finally, we emphasize that in the original formulation of \cite{boffi2025how} the consistency term was integrated over both time variables $s$ and $t$ to enable few-step sampling.
In this work, we show that integrating only over $t$ is sufficient to establish one-step generation guarantees.

\medskip
\noindent
\textbf{Neural Networks.}\quad
In this paper we use the Gaussian error linear unit (GELU) activation \citep{hendrycks2016gaussian}, defined as
$\gelu(x) = x \cdot \Phi(x)$, where $\Phi(x)$ denotes the cumulative distribution function of the standard normal distribution.
We choose GELU because it is infinitely differentiable and has bounded derivatives.
This allows us to invoke the framework of \cite{yakovlev2026simultaneous} to construct function approximations with an explicitly controlled Lipschitz constant, which is crucial for our estimation-error analysis.
This choice is also supported by the widespread use of GELU in modern neural network architectures \citep{brown2020language, peebles2023scalable, biderman23pythia}.

For a shift vector $b = (b_1, \dots, b_r) \in \R^r$ we use the shifted GELU, which maps a vector $x \in \R^r$ to $\R^r$ as $\gelu_b(x) = \big( \gelu(x_1 - b_1), \dots, \gelu(x_r - b_r) \big)$.
For $L \in \N$ and an architecture vector $W = (W_0, W_1, \dots, W_L) \in \N^{L + 1}$, a feedforward neural network of depth $L$ and architecture $W$ is a function $f : \R^{W_0} \to \R^{W_L}$ that is defined as the composition of linear maps and GELU nonlinearities:
\begin{align}
    \label{eq:feed_forward_nn_def}
    f(x) = -b_L + A_L \circ \gelu_{b_{L - 1}} \circ A_{L - 1} \circ \gelu_{b_{L - 2}} \circ \ldots \circ A_2 \circ \gelu_{b_1} \circ A_1 \circ x ,
\end{align}
where $A_j \in \R^{W_{j} \times W_{j - 1}}$ is a weight matrix and $b_j \in \R^{W_j}$ is a bias vector for every $j \in \{1, \dots, L\}$.
The maximum number of neurons in each layer is defined as $\|W\|_\infty$, which stands for the width of the network.
We further restrict the class of neural networks $f$ of the form \eqref{eq:feed_forward_nn_def} to those with at most $S \in \N$ non-zero parameters, with all weights and biases bounded by $B > 0$:
\begin{align}
    \label{eq:nn_class_def}
    \NN(L, W, S, B) = \Bigg\{
    f
    % \text{ as in \eqref{eq:feed_forward_nn_def}}
    :
    \sum_{j = 1}^L \big( \|A_j\|_0 + \|b_j\|_0 \big) \leq S, \,
    \max_{1 \leq j \leq L} \left\{ \|A_j\|_\infty \vee \|b_j\|_\infty \right\} \leq B \Bigg\}.
\end{align}

\noindent
\textbf{Smoothness classes.}\quad
Let $\Omega \subseteq \R^r$.
For a multi-index $\bk \in \Z_+^r$, the partial derivative of a function $f: \Omega \to \R$ is denoted by
\begin{align*}
    \partial^\bk f(x) = \frac{\partial^{|\bk|} f(x)}{\partial x_1^{k_1} \partial x_2^{k_2} \dots \partial x_r^{k_r}},
    \quad x \in \Omega .
\end{align*}
For $\beta > 0$, the \emph{Hölder class} $\cH^\beta(\Omega, \R, H)$ is the set of functions $f$ such that
\begin{align*}
    \max_{\substack{\bk \in \Z_+^r \\ |\bk| \leq \floor{\beta}}} \|\partial^\bk f\|_{L^\infty(\Omega)} \leq H
    \quad \text{and} \quad
    \max_{\substack{\bk \in \Z_+^r \\ |\bk| = \floor{\beta}}}
    \sup_{\substack{x, y \in \Omega \\ x \not= y}}
    \frac{|\partial^\bk f(x) - \partial^\bk f(y)|}{1 \wedge \|x - y\|_\infty^{\beta - \floor{\beta}}} \leq H .
\end{align*}
A vector-valued function $f: \Omega \to \R^m$ belongs to $\cH^\beta(\Omega, \R^m, H)$ if every component function is in $\cH^\beta(\Omega, \R, H)$.

Additionally, for $k \in \Z_+$, we denote by $W^{k, \infty}(\Omega)$ the Sobolev space of vector-valued functions whose derivatives $\partial^\bk f$ are in $L^\infty(\Omega)$ for all $\bk \in \Z_+^r$ with $|\mathbf{k}| \leq k$.
For $j \in \{0, \dots, k\}$, the associated Sobolev seminorm is defined as
\begin{align*}
    |f|_{W^{j, \infty}(\Omega)} = \max\left\{\|\partial^\bk f\|_{L^\infty(\Omega)} \, : \, \bk \in \Z_+^r, \; |\bk| = j \right\},
    \quad f \in W^{k, \infty}(\Omega) .
\end{align*}

\section{Results}
\label{sec:results}

In this section, we present the main findings of our paper.
Before proceeding with upper bounds on the flow map approximation and estimation accuracy, we first detail the assumptions on the data distribution under which the bounds hold.

\subsection{Assumptions}

The key to sharp error estimates is the following structural assumption.
As in \cite{yakovlev2025generalization,yakovlev2025implicit}, we model the target distribution $\sfp^*$ as a nonparametric Gaussian mixture, where the distribution of means is supported on a smooth low‑dimensional surface.

\begin{As}[relaxed manifold assumption, \citep{yakovlev2025generalization}]
\label{asn:relax_man}
    Let $g^* \in \cH^\beta([0, 1]^d, \R^D, H)$, for some $H > 0$ and $d, D \in \N$, with $\|g^*\|_{L^\infty([0, 1]^d)} \leq 1$.
    We consider the model where the observed samples $Y_1, \dots, Y_n$ are independent copies of a random element $X_1 \in \R^D$, generated as follows:
    \begin{align*}
        X_1 = g^*(U) + \sdata Z,
    \end{align*}
    where $U \sim \mathrm{Un}([0, 1]^d)$ and $Z \sim \cN(0, I_D)$ are independent and $0 < \smin \leq \sdata < 1$.
\end{As}

Under Assumption \ref{asn:relax_man}, the target distribution is supported near a low-dimensional surface $\mathrm{Im}(g^*)$.
Note that the effective dimension $d$ can be substantially smaller than the ambient dimension $D$.
Crucially, this assumption corresponds to scenarios where the data distribution exhibits low-dimensional structure and well-separated components, as supported by empirical studies \citep{khayatkhoei2018disconnected,pope2021the,brown2023verifying}.
In addition, this assumption relaxes the standard manifold hypothesis employed in the analysis of generative models \citep{dahal2022manifold,tang2024adaptivity,azangulov2024convergence} by eliminating the requirement of a density lower bound, a condition that is non-typical for real-world data (see the discussion in \cite{zhang2024minimax}).
Compared with \cite{yakovlev2025generalization}, we further impose that the Gaussian
noise level is strictly positive, $\sdata \geq \smin > 0$.
Positive $\smin$ allows for better spatial Lipschitz constant of the velocity field $b^*(y, t)$, which facilitates our analysis and allows us to employ a broad class of interpolant coefficients.
We acknowledge that the assumption excludes scenarios involving noiseless data, and we leave this aspect for future research.
Finally, note that Assumption \ref{asn:relax_man} aligns with Assumption 2 in \cite{ding2024characteristic}, where the target distribution is a Gaussian convolution.

We now impose regularity conditions on the interpolant coefficients $\alpha_t$ and $\beta_t$.

\begin{As}
\label{asn:flow_matching_coefs}
    The functions $\alpha_t, \beta_t : [0, 1] \to \R$ satisfy
    \begin{itemize}
        \item[(i)] $\alpha_0 = \beta_1 = 1$ and $\alpha_1 = \beta_0 = 0$.
        \item[(ii)] $\frac{1}{2}\sdata^2 \leq \alpha_t^2 + \beta_t^2 \sdata^2 \leq 1$ and $0 \leq \beta_t \leq 1$ for all $0 \leq t \leq 1$.
        \item[(iii)] $\alpha\in C^2([0,1))\cap C([0,1])$ and $\alpha^2, \beta \in C^3([0, 1])$.
        \item[(iv)] the function $\alpha_t$ is monotonically decreasing on the interval $[0, 1]$.
    \end{itemize}
\end{As}

Note that Assumption \ref{asn:flow_matching_coefs} covers a broad class of commonly used interpolants \citep{albergo2025stochastic}: the linear interpolant ($\alpha_t = 1 - t$ and $\beta_t = t$), the diffusion interpolant ($\alpha_t = \sqrt{1 - t^2}$ and $\beta_t = t$), and the trigonometric interpolant ($\alpha_t = \cos(\pi t / 2)$ and $\beta_t = \sin(\pi t / 2)$).
Our condition (iii) is stronger than the temporal regularity condition in \cite{ding2024characteristic}.
The additional derivative is used in the second-order discretization analysis.

\subsection{Flow map approximation}

In this section we introduce a class of neural networks that will be used to approximate the flow map.
We then show that, within this class, the flow map can be approximated with a rate that does not suffer from the curse of dimensionality.

Before presenting the approximation result, we first derive an explicit representation of the true velocity field, which is stated in the following proposition.

\begin{Prop}
\label{prop:true_vel_repr}
Suppose that Assumptions \ref{asn:relax_man} and \ref{asn:flow_matching_coefs} hold.
For $t \in [0, 1]$, set $\tilde\alpha_t^2 = \alpha_t^2 + \beta_t^2 \sdata^2$.
Then, the true velocity field admits the representation
\begin{align*}
    b^*(y, t) = c_1(t, \sdata) y + c_2(t, \sdata) f^*(y, t),
    \quad (y, t) \in \R^D \times [0, 1] ,
\end{align*}
where
\begin{align*}
  c_1(t, \sigma) = \frac{\dot \alpha_t \alpha_t + \dot \beta_t \beta_t \sigma^2}{\alpha_t^{2} + \beta_t^2 \sigma^2},
  \quad
  c_2(t, \sigma) = \frac{\alpha_t(\dot \beta_t \alpha_t - \beta_t \dot \alpha_t)}{\alpha_t^2 + \beta_t^2\sigma^2},
  \quad \sigma \in (0, 1], \; t \in [0, 1] .
\end{align*}
Furthermore,
\begin{align}
\label{eq:f_star_def}
    f^*(y, t) = \frac{\integral{[0, 1]^d}g^*(u)\exp\left(-\frac{\|y - \beta_t g^*(u)\|^2}{2\tilde\alpha_t^2}\right) \dd u }{ \integral{[0, 1]^d}\exp\left(-\frac{\|y - \beta_t g^*(u)\|^2}{2\tilde\alpha_t^2}\right) \dd u },
    \quad (y, t) \in \R^D \times [0, 1] .
\end{align}

\end{Prop}

The proof of Proposition \ref{prop:true_vel_repr} is given in Appendix \ref{sec:prop_true_vel_repr_proof}.
The proposition reveals that the true velocity field is a sum of a linear and a bounded term.
This representation, already exploited in the analysis of diffusion models \citep{yakovlev2025generalization}, is the key insight that allows us to design the class of average velocities.
Specifically, we approximate the bounded component of the velocity field using a cubic
B-spline quasi-interpolant.
The formal definition of this quasi-interpolant is detailed below.

\begin{Def}[B-spline quasi-interpolant]
\label{def:b_spline_interp}
    For a fixed $t_0 \in [0, 1)$, $N \in \N$, and a set of points $\{X_k\}_{k = 0}^N \subseteq \R^D$, let the coefficients $\{\alpha_k\}_{k = -1}^{N + 1}$ be defined as $\alpha_k = X_k$ for $k \in \{0, 1, \dots, N\}$, with boundary values $\alpha_{-1} = 2X_0 - X_1$ and $\alpha_{N + 1} = 2X_N - X_{N - 1}$.
    The cubic B-spline quasi-interpolant $I : [0, 1 - t_0] \to \R^D$ is defined by
    \begin{align*}
        I(t; X_0, X_1, \dots, X_N) = \sum_{k = -1}^{N + 1} \alpha_k B_3(h^{-1}t - (k - 2)),
        \quad t \in [0, 1 - t_0] ,
    \end{align*}
    where $B_3$ is the cubic B-spline basis function supported on $[0,4]$ and $h=(1-t_0)/N$, with grid points $\tau_k=kh$.
\end{Def}

Cubic B-splines are selected specifically for their temporal regularity, which is essential for an accurate approximation of the consistency penalty $R(v)$.
We can now define the class of average velocities.

\begin{Def}[average velocities class]
\label{def:av_vel_class}
Let $N \in \N$, $t_0 \in [0, 1)$, and $L_f > 0$ be given.
Suppose that $\{f_k \, : \, 0 \leq k \leq N\} \subseteq \NN(L, W, S, B)$ are $L_f$-Lipschitz functions satisfying
\begin{align*}
    (i) &\quad \|f_k\|_{L^\infty(\R^D)}  \lesssim \sqrt{D},
    \qquad
    (ii) \quad \|\nabla f_k(y)\| \lesssim \smin^{-2},
    \quad \text{for all } \|y\| \leq R .
\end{align*}
Set $h = (1 - t_0) / N$ and $\tau_k = hk$ for $k \in \{0, 1, \dots, N\}$.
For any $y\in\R^D$ and $\sigma\in[\smin,1]$, put
\[
    b_k(z)=c_1(\tau_k,\sigma)z+c_2(\tau_k,\sigma)f_k(z),
    \quad 0\leq k\leq N.
\]
Define the Heun iterates by
\begin{align*}
    Y_{k+1}=Z_k+h b_k(Z_k),
    \quad
    Z_{k+1}=Z_k+\frac h2\left(b_k(Z_k)+b_{k+1}(Y_{k+1})\right),\quad Z_0=y.
\end{align*}
Here, the functions $c_1$ and $c_2$ are defined in Proposition \ref{prop:true_vel_repr}.
For all $(y, t) \in \R^D \times [0, 1 - t_0]$ we define the velocity field
\[
    v_{t, t}(y) = c_1(t, \sigma) y + c_2(t, \sigma) f(y, t),
    \quad f(y, t) = I(t; f_0(y), f_1(y), \dots, f_N(y)) .
\]
Similarly, for all $t \in (0, 1 - t_0]$ and $y \in \R^D$, we define the average velocity as
\[
    v_{0, t}(y) = t^{-1} \left( I(t; Z_0, Z_1, \dots, Z_N) - y\right).
\]
Here, $I(t;\cdot)$ denotes the cubic B-spline quasi-interpolant from Definition \ref{def:b_spline_interp}.
We denote the resulting class of average velocity candidates by $\cV(L, W, S, B, N, R, L_f)$.

\end{Def}

Having formally defined the class of average velocities, we now discuss the motivation behind our architectural choices.

The representation of the flow map via the average velocity in \eqref{eq:v_star_def} motivates the use of residual neural networks \citep{he2016resnet} for $v_{0, t}$ and
$v_{t, t}$ (Definition \ref{def:av_vel_class}), rather than the standard feedforward architectures.
Residual architectures have proven effective for approximating trajectories of ODEs \citep{baggenstos2023discrete} and possess favorable approximation complexity bounds \citep{muller2019resnets,oono2019resnets,liu2024resnets}.
Moreover, the iterative construction in Definition \ref{def:av_vel_class} explicitly implements Heun's method, a second-order Runge-Kutta discretization scheme for ODEs.

Although higher-order Runge-Kutta methods could potentially refine the approximation, their analysis is considerably more involved.
Specifically, analyzing higher-order schemes requires determining not only the order of accuracy but also the precise dependence on the ambient dimension $D$, the stopping time $t_0$, and the noise parameter $\sigma_{\min}$.
We find that a second-order method offers a good balance between numerical accuracy and theoretical tractability for our purposes.
Similarly, we use cubic B-splines because they provide the temporal regularity needed to control the consistency penalty $R(v)$ while keeping the analysis tractable.
Although B-splines are not, to our knowledge, explicitly used in existing architectures, we introduce them solely to guarantee the temporal smoothness of the flow-map estimator, which is essential for our estimation error analysis.
Moreover, cubic B-splines can themselves be approximated by neural networks, so the approximation class in Definition \ref{def:av_vel_class} can be realized by a ResNet architecture.
We keep the B-spline formulation because it makes the approximation arguments substantially more transparent.

It is important to note that the neural networks in Definition \ref{def:av_vel_class} are required to be Lipschitz continuous.
This property is essential for the subsequent estimation analysis.
In practice, this constraint can be enforced using spectral normalization \citep{miyato2018spectral} or Lipschitz neural-network architectures \citep{anil2019sorting, qi2023lipsformer}.
This requirement is consistent with other recent theoretical studies, including the work of \cite{ding2024characteristic}, which also emphasizes the importance of Lipschitz continuity.

We now present one of our main results, which establishes that the class of average velocities, as defined in Definition \ref{def:av_vel_class}, is sufficiently rich to approximate the target flow map without the curse of dimensionality.

\begin{Th}[approximation of the flow map]
\label{th:avg_vel_approx}
    Assume Assumptions \ref{asn:relax_man} and \ref{asn:flow_matching_coefs} hold.
    Fix $t_0\in[0,1/2]$.
    Then, for any $\eps\in(0,1/2]$ satisfying $\eps^\beta\leq c\smin^4$ and
    \begin{align*}
        D^{3/2} \exp\{\cO(\smin^{-2})\} (\eps^{\beta}\vee\eps)\left(\log(1 / \eps) + \log(D \smin^{-2})\right)
        \lesssim 1 ,
    \end{align*}
    there exists an average velocity candidate $\tilde{v}$ as in Definition \ref{def:av_vel_class} with the parameters
    \begin{align*}
        N \asymp \eps^{-\beta / 2}\smin^{-2},
        \quad
        L_f \asymp (\smin^{-2} D \log(1 / \eps))^{\cO(1)},
        \quad
        R \asymp \exp\{\cO(\smin^{-2})\} D \log(\eps^{-1}D \smin^{-2})
    \end{align*}
    such that
    \begin{align*}
        \E[\ell_{SD}(\tilde v, X_1) - \ell_{SD}(v^*, X_1)]
    \lesssim (D \log(1 / \eps))^{\cO(1)} \exp\{\cO(\smin^{-4})\} \eps^\beta .
    \end{align*}
    Moreover, the architecture of $\tilde v$ is as follows:
    \begin{align*}
        &L
        \lesssim \smin^{-2} \log(D \smin^{-2} \log(1 / \eps)),
        \quad
        \|W\|_\infty \vee S
        \lesssim \eps^{-d}\exp\{\cO(\smin^{-2})\} (D \log(1 / \eps))^{\cO(1)},
        \\&
        \log B
        \lesssim \exp\{\cO(\smin^{-2})\} (D \log(1 / \eps))^{\cO(1)}
    \end{align*}
    The hidden constants behind $\cO(\cdot)$ may depend only on $d$, $\beta$, and $H$.
    
\end{Th}

The proof of Theorem \ref{th:avg_vel_approx} is deferred to Appendix \ref{sec:th_avg_vel_approx_proof}.
The proof uses the relation between the score function $s^*(y, t) = \nabla \log\sfp_t(y)$ and the velocity field $b^*(y, t)$.
Specifically, $b^*(y, t)$ can be expressed as a linear combination of $y$ and $s^*(y, t)$ (see, for instance, \cite[Proposition 3.1]{ding2024characteristic} and \cite[Section 10]{lipman2024flow}).
Building on this relation, we utilize the simultaneous approximation of the score function and its derivatives \citep{yakovlev2026simultaneous}.
This ensures that the resulting average velocity approximation is spatially Lipschitz, a property that is essential for the estimation analysis below.

Theorem \ref{th:avg_vel_approx} establishes that the average velocity can be approximated to accuracy $\tilde\cO(\eps^\beta)$ by the residual neural network architecture defined in Definition \ref{def:av_vel_class}, using only
$\tilde\cO(\eps^{-d - \beta/2})$ non-zero parameters.
Notably, the dependence on the ambient dimension $D$ is polynomial, ensuring that the bound remains informative even in high-dimensional regimes where $D \gtrsim \log n$.
Although the bounds scale exponentially with $\sigma_{\min}^{-2}$, relaxing this dependence to a polynomial one is an important direction for future work.

Under similar assumptions on the data distribution, \cite{yakovlev2025generalization, yakovlev2026simultaneous} derived approximation bounds for the score function $s^*(y, t) = \nabla \log \sfp_t(y)$ that avoid the curse of dimensionality.
Specifically, achieving an accuracy of $\tilde\cO(\eps^{2\beta})$ for the score function requires $\tilde\cO(\eps^{-d})$ non-zero parameters, which is slightly lower than the complexity presented in Theorem \ref{th:avg_vel_approx}.
We attribute this gap to the fact that approximating the score function corresponds to estimating the instantaneous velocity, whereas estimating the average velocity is inherently more challenging, as it involves the analysis of discretization error in the underlying ODE.

Compared with \cite{ding2024characteristic}, we observe that our approach avoids the curse of dimensionality under more realistic assumptions.
More specifically, by setting $\beta \geq 1$ and assuming $g^*(u) = Au$ for some $A \in \R^{D \times d}$ with orthonormal columns, we recover Assumption 3 of \cite{ding2024characteristic}, which describes the target distribution as the convolution of a distribution supported on a low-dimensional linear subspace with Gaussian noise.
In their Lemma H.7, they achieved $\cO(\eps^2)$ approximation accuracy for their semi-group-penalized risk using $\tilde\cO(\eps^{-d - 2})$ parameters.
In contrast, Theorem \ref{th:avg_vel_approx} provides similar approximation guarantees under more realistic assumptions by treating the low-dimensional structure of the data distribution more generally.

Finally, we expect that similar approximation guarantees can be established for other one-step flow-based objectives such as Mean Flow \citep{geng2025mean}, shortcut loss \citep{frans2025one}, and Alpha Flow loss \citep{zhang2026alphaflow}.
Crucially, the framework in Theorem \ref{th:avg_vel_approx} is not specific to self-distillation and is therefore of independent interest beyond the present setting.

\subsection{Flow map estimation}

In this section, we establish the sample complexity of learning the flow map via the self-distillation objective given by \eqref{eq:self_dist_loss_def}.
The following theorem constitutes one of the central results of this paper.
It demonstrates that the flow map can be estimated without the curse of dimensionality.

\begin{Th}
\label{th:estim_main}
Assume Assumptions \ref{asn:relax_man} and \ref{asn:flow_matching_coefs} hold.
Fix $t_0\in[0,1/2]$.
Assume further that the sample size is sufficiently large in the sense that it satisfies
\begin{align*}
    D^{3/2} e^{\cO(\smin^{-2})} \log(nD\smin^{-2}) \left(\frac{\lambda(t_0)}{n}\right)^{\frac{2(\beta\wedge1)}{5\beta + 2d}} \lesssim 1 ,
    \quad
    \left(\frac{\lambda(t_0)}{n}\right)^{\frac{2\beta}{5\beta+2d}} \lesssim \smin^4,
    \quad n \gtrsim \lambda(t_0),
\end{align*}
where $\lambda(t_0) = |\alpha_t|_{W^{1, \infty}([0, 1 - t_0])} \vee e$.
Then, with probability at least $(1 - \delta)$, the empirical risk minimizer \eqref{eq:avg_vel_erm_def} over the class $\cV(L, W, S, B, N, R, L_f)$ given by Definition \ref{def:av_vel_class} with the architecture
\begin{align*}
    &L \lesssim \smin^{-2} \log(D \smin^{-2} \log n),
    \quad \|W\|_\infty \vee S \lesssim \left( \frac{n}{\lambda(t_0)} \right)^{\frac{2d}{5\beta + 2d}} \exp\{\cO(\smin^{-2})\} (D \log n)^{\cO(1)},
    \\
    &\log B \lesssim \exp\{\cO(\smin^{-2})\} (D \log n)^{\cO(1)}
\end{align*}
and the parameters
\begin{align*}
  L_f \asymp (\smin^{-2}D \log n)^{\cO(1)},
  \quad N \asymp (n / \lambda(t_0))^{\beta / (5\beta + 2d)} \sigma_{\min}^{-2},
  \quad 
  R \asymp \exp\{\cO(\smin^{-2})\}D \log(nD\smin^{-2}) 
\end{align*}
satisfies the inequality
\begin{align*}
    &\E_{X_1}[\ell_{SD}(\widehat{v}, X_1) - \ell_{SD}(v^*, X_1)]
    \lesssim (D \, e^{\sigma_{\min}^{-4}} \log n)^{\cO(1)}
    \lambda(t_0)
    \left(\frac{\lambda(t_0)}{n} \right)^{\frac{2\beta}{5\beta + 2d}} \log(\lambda(t_0)) \log(4 / \delta) .
\end{align*}
Here, the hidden constant may depend on $d$, $\beta$, and $H$ only.
\end{Th}

The proof of Theorem \ref{th:estim_main} is provided in Appendix \ref{sec:th_estim_main_proof}.
The proof uses an $\epsilon$-net argument in conjunction with the Bernstein-type inequality for unbounded random variables \citep{Lecue12oracle}.
Furthermore, we establish that the loss landscape is well-behaved in the sense that it satisfies the Bernstein condition, which enables sharp convergence rates.

Theorem \ref{th:estim_main} establishes an excess risk rate of $\tilde\cO(n^{-2\beta / (5\beta + 2d)})$, with a polynomial dependence on the ambient dimension $D$ and $\lambda(t_0)$.
Under similar data distribution assumptions, this behaviour is consistent with results from the score estimation literature \citep{chen2023score,yakovlev2025generalization}.
To the best of our knowledge, Theorem \ref{th:estim_main} is the first result to provide convergence rates for flow map estimation that do not suffer from the curse of dimensionality under reasonable assumptions on the low-dimensional structure of the data distribution.
The rate exhibits exponential dependence on $\smin^{-4}$, consistent with Theorem \ref{th:avg_vel_approx}.
Whether this bound can be improved remains an important open problem for future work.

Finally, we observe that the rates in the score estimation literature \citep{yakovlev2025generalization,yakovlev2025implicit}, under similar assumptions on the data distribution, scale as $\tilde\cO(n^{-2\beta / (2\beta + d)})$, which decays faster than the rate in Theorem \ref{th:estim_main}.
As discussed following Theorem \ref{th:avg_vel_approx}, this discrepancy arises because flow map estimation is a more challenging task, requiring a more involved approximation error analysis.

\subsection{Distribution estimation}

In this section, we analyze the estimation of the target distribution with respect to the 2-Wasserstein distance.
The following theorem demonstrates that the one-step generator, learned via the self-distillation objective \eqref{eq:self_dist_loss_def}, recovers the target distribution
without suffering from the curse of dimensionality.

\begin{Th}[distribution estimation]
\label{th:distr_estim}
    Let $\cV(L, W, S, B, N, R, L_f)$ be the class defined in Theorem \ref{th:estim_main}, and let $\widehat{\Psi}_{0, t}(y) = y + t \cdot \widehat{v}_{0, t}(y)$ be a flow map, where $\widehat{v}$ is the empirical risk minimizer from \eqref{eq:avg_vel_erm_def}.
    Define $\widehat \sfp_{1 - t_0} = \mathrm{Law}(\widehat \Psi_{0, 1 - t_0}(X_0))$.
    Then, with probability at least $(1 - \delta)$, it holds that
    \begin{align*}
        W_2\left(\widehat \sfp_{1 - t_0}, \sfp^* \right)
        &\lesssim D^{\cO(1)} \exp\{\cO(\smin^{-4})\}\sqrt{\lambda(t_0)}
        \left(\frac{\lambda(t_0)}{n} \right)^{\frac{\beta}{5\beta + 2d}} \mathscr{L}(t_0, n, \delta)
        \\&\quad
        + \sqrt{D}(\alpha_{1 - t_0} \vee (1 - \beta_{1 - t_0})),
    \end{align*}
    where $\mathscr{L}(t_0, n, \delta) = (\log n)^{\cO(1)}\sqrt{\log(\lambda(t_0)) \log(4 / \delta)}$ .
    The hidden constant behind $\cO(\cdot)$ depends only on $d$, $\beta$, and $H$.
\end{Th}

The proof of Theorem \ref{th:distr_estim} is moved to Appendix \ref{sec:th_distr_estim_proof}.
The primary technical challenge is that the standard analysis based on the Alekseev-Gr\"obner lemma \citep{benton2024error,
fukumizu2025flow} requires the velocity field estimate to be globally Lipschitz, with the $W_2$ distribution error scaling exponentially with this Lipschitz constant.
A direct application of these results would not yield a consistent estimate, as the global Lipschitz constant guaranteed by Theorem \ref{th:estim_main} grows at least logarithmically with the sample size $n$.
However, the local Lipschitz constant provided in Theorem \ref{th:estim_main} is well-behaved.
By extending the standard analysis, we demonstrate that local Lipschitzness is sufficient to obtain a sharp convergence rate.

Theorem \ref{th:distr_estim} has several important implications for different choices of the interpolant coefficients $\alpha_t$ and $\beta_t$, as detailed in the following remark.

\begin{Rem}
\label{rem:distr_estim}
Under the assumptions of Theorem \ref{th:distr_estim}, consider the following two cases:
\begin{enumerate}
    \item \textbf{Linear Interpolant:} For the linear interpolant ($\alpha_t = 1- t$ and $\beta_t = t$), we have $\lambda(t_0) \asymp 1$ for all $t_0 \in [0, 1)$.
    In this case, for $t_0 = 0$ and all sufficiently large $n$, the 2-Wasserstein distance scales as $W_2(\widehat \sfp_1, \sfp^*) = \tilde\cO(n^{-\beta / (5\beta + 2d)})$.

    \item \textbf{Diffusion Interpolant:} For the diffusion interpolant ($\alpha_t = \sqrt{1 - t^2}$ and $\beta_t = t$), we have $\lambda(t_0) \lesssim t_0^{-1/2}$ and $\alpha_{1 - t_0} \vee (1 - \beta_{1 - t_0}) \lesssim \sqrt{t_0}$.
    Therefore, setting $t_0 = n^{-4\beta / (17\beta + 6d)}$ for all sufficiently large $n$, we obtain a convergence rate of $W_2(\widehat \sfp_{1 - t_0}, \sfp^*) = \tilde{\cO}(n^{-2\beta / (17\beta + 6d)})$.
\end{enumerate}
In both cases the dependence on the ambient dimension $D$ is polynomial, which is crucial in high-dimensional settings.
    
\end{Rem}

Building on Theorem \ref{th:distr_estim} and Remark \ref{rem:distr_estim}, we compare our findings with the results presented in \cite{ding2024characteristic}.
In Corollary 3.16 of their work (convergence of linear interpolant under low-dimensional hypothesis), the authors obtain a rate for the squared $W_2$ distance of $\tilde{\cO}(n^{-4 / (3d + 9)})$.
Consequently, if $\beta > 4d / (3d - 1)$, our approach achieves a faster rate of convergence.
Note that such a choice of $\beta$ is admissible in the linear case, where $g^*$ is infinitely smooth.
Furthermore, our guarantees hold for a significantly broader class of data distributions, extending beyond linear subspaces.
To the best of our knowledge, Theorem \ref{th:distr_estim} is the first such result in the literature on one-step flow-based and diffusion-based generative models.
Whether these distribution estimation bounds are optimal remains an important open question for future research.

Finally, we discuss the exponential dependence on $\sigma_{\min}^{-4}$ in Theorems \ref{th:estim_main} and \ref{th:distr_estim}.
To our knowledge, this dependence can be mitigated under substantially stronger assumptions on the data distribution.
For example, if the target distribution $\sfp^*$ is weakly log-concave, the result of \cite[Theorem 3]{stephanovitch2026lipschitz} suggests that the factor $\exp\{\cO(\smin^{-4})\}$ can be replaced by a constant.
However, that result does not apply under Assumption \ref{asn:relax_man}, and extending it to our setting remains an open problem.
More broadly, obtaining simultaneously polynomial dependence on the ambient dimension $D$, optimal statistical rates in the sample size $n$, and only mild dependence on $\sigma_{\min}$ appears to be a fundamental theoretical challenge.
Our analysis establishes favourable dependence on $D$ and $n$.
We view this as a meaningful step toward understanding the statistical properties of one-step flow-based models.

\section{Conclusion}
\label{sec:conclusion}

In this work, we developed a theoretical framework for the estimation of flow maps in one-step generative models.
We derived approximation and estimation bounds that avoid the curse of dimensionality, assuming that the target distribution is the convolution of a distribution supported on a low-dimensional surface with Gaussian noise.
These findings provide a first step toward explaining the empirical success of one-step flow-based generative models.
Future work will focus on further refining the derived bounds, particularly by exploring the use of higher-order B-splines and more sophisticated Runge-Kutta integration schemes.

\bibliographystyle{abbrvnat}
\bibliography{references}

%%%%%%%%%%%%%%%%%%%%%%%%%%%%%%%%%%%%%%%%%%%%%%%%%%%%%%%%%%%%

\appendix

\section{Proofs of supporting propositions}

This section collects the proofs of Propositions \ref{prop:fm_loss_dif} and \ref{prop:true_vel_repr}.

\subsection{Proof of Proposition \ref{prop:fm_loss_dif}}
\label{sec:prop_fm_loss_dif_proof}

From the definition of the flow matching loss defined in \eqref{eq:fm_loss_def} we obtain that
\begin{align*}
    &\E[\ell_{FM}(b, X_1) - \ell_{FM}(b^*, X_1)]
    = \integral{0}^{1 - t_0} \E \|b(X_t, t) - b^*(X_t, t)\|^2 \, \dd t
    \\&
    + 2\integral{0}^{1 - t_0} \E \langle b(X_t, t) - b^*(X_t, t) , b^*(X_t, t) - \dot X_t \rangle \, \dd t .
\end{align*}
Taking into account that $b^*(y, t) = \E[\dot X_t \, | \, X_t = y]$ and using the tower property of the conditional expectation, we conclude that
\begin{align*}
    \E \langle b(X_t, t) - b^*(X_t, t) , b^*(X_t, t) - \dot X_t \rangle
    = \E \left[ \langle b(X_t, t) - b^*(X_t, t), b^*(X_t, t) - \E[\dot X_t \, | \, X_t] \rangle \right]
    = 0 .
\end{align*}
Consequently, the claim readily follows.

\myendproof

\subsection{Proof of Proposition \ref{prop:true_vel_repr}}
\label{sec:prop_true_vel_repr_proof}

Recall that $b^*(y, t) = \E[\dot X_t \, | \, X_t = y]$ for $t<1$.
Therefore, it holds that
\begin{align*}
    b^*(y, t) = \E[\E[\dot \alpha_t X_0 + \dot \beta_t(g^*(u) + \sdata Z) \, | \, X_t = y, U = u] \, | \, X_t = y] .
\end{align*}
Next, observe that for a fixed $U = u$, it holds that
\begin{align*}
    (X_0, Z, X_t)^\top \, | \, (U = u)
    \sim \cN\left(
    \begin{pmatrix}
        0_D \\ 0_D \\ \beta_t g^*(u)
    \end{pmatrix},
    \begin{pmatrix}
        I_D & 0_{D \times D} & \alpha_t I_D \\
        0_{D \times D} & I_D & \beta_t \sdata I_D \\
        \alpha_t I_D & \beta_t \sdata I_D & \tilde\alpha_t^2 I_D
    \end{pmatrix}
    \right) .
\end{align*}
Therefore, using the conditional law of a Gaussian vector, we deduce that
\begin{align*}
    \E[X_0 \, | \, X_t = y, \, U = u] = \frac{\alpha_t (y - \beta_t g^*(u))}{\tilde\alpha_t^2},
    \quad
    \E[Z \, | \, X_t = y, U = u] = \frac{\beta_t \sdata (y - \beta_t g^*(u))}{\tilde\alpha_t^2} .
\end{align*}
Hence, we arrive at
\begin{align}
\label{eq:b_star_aux_cond}
    b^*(y, t) = \E\left[ \frac{(\alpha_t\dot\alpha_t + \beta_t \dot\beta_t \sdata^2)y}{\tilde\alpha_t^2}  + \left(\dot \beta_t - \frac{\beta_t(\alpha_t \dot \alpha_t + \beta_t\dot\beta_t\sdata^2)}{\tilde\alpha_t^2}\right) g^*(U) \, \bigg| \, X_t = y \right] .
\end{align}
We next notice that $X_t | (U = u) \sim \cN(\beta_t g^*(u), \tilde\alpha_t^2 I_D)$ and, hence, $\E[g^*(U) \, | \, X_t = y] = f^*(y, t)$, where $f^*(y, t)$ is defined in \eqref{eq:f_star_def}.
Substituting this identity into \eqref{eq:b_star_aux_cond} yields the desired result.

\myendproof

\section{Proof of Theorem \ref{th:avg_vel_approx}}
\label{sec:th_avg_vel_approx_proof}

The proof constructs $v\in\cV(L,W,S,B,N,R,L_f)$ by controlling both terms of
$\ell_{SD}(v,\cdot)-\ell_{SD}(v^*,\cdot) =\ell_{FM}(v,\cdot)-\ell_{FM}(v^*,\cdot) + R(v)$.
The flow matching term requires a uniform in time approximation of the velocity field $b^*$, while $R(v)$ also requires temporal smoothness and an accurate trajectory.
An approximation for a fixed time alone is therefore insufficient.
We first approximate $b^*(\cdot,\tau_k)$ at each grid time by GELU networks, then interpolate them with cubic B-splines to obtain a uniform in time bound and the temporal regularity needed for $R(v)$.
Next, Heun's method and a second B-spline quasi-interpolation yield the average velocity and the flow map, with controlled $L^2$ and derivative errors.
Combining these estimates yields the desired approximation error bound.

For clarity, we split the proof of Theorem \ref{th:avg_vel_approx} into several steps.

\noindent
\textbf{Approximation of the velocity field.}\quad
We begin by introducing a uniform grid $\tau_k = kh$ with $k \in \{0, \dots, N\}$, where the step size  $h = (1 - t_0) / N$.
Here, $N \in \N$ will be specified later in the proof.
Our first objective is to construct an approximation of the velocity field at each grid point $\tau_k \in [0, 1 - t_0]$, a result established in the following lemma.

\begin{Lem}[approximation of the velocity field]
\label{lem:vel_approx}
    Suppose Assumptions \ref{asn:relax_man} and \ref{asn:flow_matching_coefs} hold, and fix $t \in [0, 1 - t_0]$.
    Assume further that $\eps \in (0,1/2]$ satisfies the condition
    \begin{align}
    \label{eq:vel_appr_eps_constr}
        D^{3/2} \exp\{\cO(\smin^{-2})\} (\eps^{\beta}\vee\eps)\left(\log(1 / \eps) + \log(D \smin^{-2})\right)
        \lesssim 1 .
    \end{align}
    Then, there exists a GELU neural network $\bar f \in \NN(L, W, S, B)$ such that the corresponding velocity field approximation, defined by
    \begin{align*}
        \bar b(y) = c_1(t, \sdata) y + c_2(t, \sdata) \bar f(y),
        \quad y \in \R^D,
    \end{align*}
    satisfies the following properties:
    \begin{align*}
        (i) &\quad \max_{0 \leq s \leq 1 - t_0} \|\bar f-f^*(\cdot,t)\|^2_{L^2(\sfp_s^*)}
        \vee \|\bar b-b^*(\cdot,t)\|^2_{L^2(\sfp_s^*)}
        \lesssim \left(D e^{\smin^{-2}}\log(1/\eps)\right)^{\cO(1)}\eps^{2\beta},
        \\
        (ii) &\; \sup_{y \in \cK_t}\|\nabla \bar f(y)\| \lesssim \sdata^{-2},
        \; |\bar f|_{L^\infty(\R^D)} \lesssim \sqrt{D} ,
        \; \max_{1 \leq l \leq D}|\bar f_l|_{W^{1, \infty}(\R^D)} \lesssim \left(\smin^{-2} D \log(1 / \eps)\right)^{\cO(1)} .
    \end{align*}
    Here, $\cK_t$ is defined as
    \begin{align*}
            \cK_t = \left\{y \in \R^D \, : \, \min_{u \in [0, 1]^d}\|y - \beta_t g^*(u)\| \leq R \right\} ,
    \end{align*}
    with $R \asymp \exp\{\cO(\smin^{-2})\} D \left(\log(1 / \eps) + \log(D \smin^{-2})\right)$.
    Moreover, the architecture of $\bar f$ is given by
    \begin{align*}
        &L
        \lesssim \smin^{-2} \log(D \smin^{-2} \log(1 / \eps)),
        \quad
        \|W\|_\infty \vee S
        \lesssim \eps^{-d}\exp\{\cO(\smin^{-2})\} (D \log(1 / \eps))^{\cO(1)},
        \\&
        \log B
        \lesssim \exp\{\cO(\smin^{-2})\} (D \log(1 / \eps))^{\cO(1)}
    \end{align*}

\end{Lem}

The proof of Lemma \ref{lem:vel_approx} is provided in Appendix \ref{sec:lem_vel_approx_proof} and builds on \cite{yakovlev2026simultaneous}.
We choose a substantially larger radius $R$ than in \cite{yakovlev2026simultaneous} because our analysis requires the velocity field approximation to be Lipschitz on a ball whose radius grows as $\exp\{\cO(\smin^{-2})\}$.
Now, we analyze the derivatives of the constructed velocity field approximation.
Let
\begin{align*}
    \bar b(y, \tau_k) = c_1(\tau_k, \sdata) y + c_2(\tau_k, \sdata) \bar f(y, \tau_k)
\end{align*}
be the velocity field approximation provided by Lemma \ref{lem:vel_approx} for each $k \in \{0, 1, \dots, N\}$ and $\eps \in (0, 1)$ satisfying \eqref{eq:vel_appr_eps_constr}.
Apply Lemma \ref{lem:vel_approx} with a tube radius at least twice the class radius below.
Introducing $\tilde f_k(y)=\bar f(y,\tau_k)$ and invoking claim (ii), we can set
\begin{align}
\label{eq:R_avg_vel_def}
    R \asymp \exp\{\cO(\smin^{-2})\} D \left(\log(1 / \eps) + \log(D \smin^{-2})\right)
\end{align}
in Definition \ref{def:av_vel_class}.
Therefore, for every $k \in \{0, \dots, N\}$, the function $\tilde f_k$ is Lipschitz on $\cB(0, R)$.
The Lipschitz constant scales as $\sdata^{-2}$.
Furthermore, $\tilde f_k$ is globally $L_f$-Lipschitz, where $L_f$ is given by
\begin{align*}
    L_f \asymp (\smin^{-2} D \log(1 / \eps))^{\cO(1)} .
\end{align*}
Next, observe that, by Assumption \ref{asn:flow_matching_coefs}, the functions $c_1$ and $c_2$, defined in Proposition \ref{prop:true_vel_repr}, satisfy
\begin{align}
\label{eq:c_j_der_bounds}
    |c_j(\cdot, \sdata)|_{W^{k, \infty}([0, 1 - t_0])} \lesssim \sdata^{-2 - 2k},
    \quad j \in \{1, 2\}, \; k \in \{0, 1, 2\} .
\end{align}
Consequently, for all $y \in \R^D$, it follows that
\begin{align*}
    \|\bar b(y, \tau_k)\|
    &\leq |c_1(\cdot, \sdata)|_{W^{0, \infty}([0, 1 - t_0])} \|y\| + |c_2(\cdot, \sdata)|_{W^{0, \infty}([0, 1 - t_0])} \|\tilde f_k(y)\|
    \\&
    \lesssim \sdata^{-2} \|y\| + \sdata^{-2}\sqrt{D} .
\end{align*}
Furthermore, for all $y \in \cK_{\tau_k}$, we have that
\begin{align*}
    \|\nabla \bar b(y, \tau_k)\|
    \leq |c_1(\cdot, \sdata)|_{W^{0, \infty}([0, 1 - t_0])} + |c_2(\cdot, \sdata)|_{W^{0, \infty}([0, 1 - t_0])}\|\nabla \tilde f_k(y)\|
    \lesssim \sdata^{-4} .
\end{align*}
We now establish bounds on the derivatives of $f^*$, as defined in \eqref{eq:f_star_def}, which are provided in the following lemma.

\begin{Lem}
\label{lem:f_star_deriv}
    Suppose Assumptions \ref{asn:relax_man} and \ref{asn:flow_matching_coefs} hold, and let $f^*$ be defined as in \eqref{eq:f_star_def}.
    Then, for all $y \in \R^D$ and $t \in [0, 1 - t_0]$, the following bounds hold:
    \begin{align*}
        (i) &\quad \|f^*(y, t)\| \leq 1,
        \quad
        \|\nabla f^*(y, t)\| \leq 2\sdata^{-2},
        \\
        (ii) &\quad
        \|\partial_t f^*(y, t)\| \leq 2 \|y\| \cdot |\beta_t / \tilde\alpha_t^2|_{W^{1, \infty}([0, 1 - t_0])} + |\beta_t^2 / \tilde\alpha_t^2|_{W^{1, \infty}([0, 1 - t_0])} ,
        \\
        (iii) &\quad
        \|\nabla^2 f^*(y, t)\| \leq 32\sdata^{-4},
        \\
        (iv) &\quad
        \|\nabla \partial_t f^*(y, t)\|
        \leq 2|\beta_t/\tilde\alpha_t^2|_{W^{1,\infty}([0,1-t_0])}
        \\&\qquad+16\sdata^{-2}\left(\|y\|\,|\beta_t/\tilde\alpha_t^2|_{W^{1,\infty}([0,1-t_0])}
        +|\beta_t^2/\tilde\alpha_t^2|_{W^{1,\infty}([0,1-t_0])}\right),
        \\
        (v) &\quad
        \|\partial^2_{tt} f^*(y, t)\| \leq 2 |\beta_t / \tilde\alpha_t^2|_{W^{2, \infty}([0, 1 - t_0])} \|y\| + |\beta_t^2 / \tilde\alpha_t^2|_{W^{2, \infty}([0, 1 - t_0])}
        \\&\quad
        + 8 \left( \|y\| \cdot |\beta_t / \tilde\alpha_t^2|_{W^{1, \infty}([0, 1 - t_0])} + \frac{1}{2} |\beta_t^2 / \tilde\alpha_t^2|_{W^{1, \infty}([0, 1 - t_0])} \right)^2 .
    \end{align*}

\end{Lem}

The proof of Lemma \ref{lem:f_star_deriv} is postponed to Appendix \ref{sec:lem_f_star_deriv_proof}.
By Lemma \ref{lem:f_star_deriv}, we have that, for all $y \in \R^D$, the true velocity field $b^*(y, t)$ satisfies
\begin{align*}
    \| b^*(y, t)\|
    &\lesssim |c_1(\cdot, \sdata)|_{W^{0, \infty}([0, 1 - t_0])} \|y\| + |c_2(\cdot, \sdata)|_{W^{0, \infty}([0, 1 - t_0])} \|f^*(y, t)\|
    \\&
    \lesssim \sdata^{-2}(1 + \|y\|) .
\end{align*}
Furthermore, we have that
\begin{align}
\label{eq:true_vel_jac}
    \notag
    \|\nabla b^*(y, t)\|
    &\leq |c_1(\cdot, \sdata)|_{W^{0, \infty}([0, 1 - t_0])} + |c_2(\cdot, \sdata)|_{W^{0, \infty}([0, 1 - t_0])} \|\nabla f^*(y, t)\|
    \\&
    \lesssim \sdata^{-4}
\end{align}
and
\begin{align*}
    \|\nabla^2 b^*(y, t)\|
    \leq |c_2(\cdot, \sdata)|_{W^{0, \infty}([0, 1 - t_0])} \|\nabla^2 f^*(y, t)\|
    \lesssim \sdata^{-6} .
\end{align*}
Next, Assumption \ref{asn:flow_matching_coefs} ensures that
\begin{align}
    \label{eq:beta_t_alpha_t_der}
    |\beta_t^j / \tilde\alpha_t^2|_{W^{k, \infty}([0, 1 - t_0])}
    \lesssim \sdata^{-2k - 2},
    \quad \text{for all } j \in \{1, 2\},
    \quad k \in \{0, 1, 2\} .
\end{align}
Therefore, we have that
\begin{align}
\label{eq:true_vel_dt}
    \notag
    \|\partial_t b^*(y, t)\|
    &\leq |c_1(\cdot, \sdata)|_{W^{1, \infty}([0, 1 - t_0])} \|y\| + |c_2(\cdot, \sdata)|_{W^{1, \infty}([0, 1 - t_0])}
    \\&\quad
    + |c_2(\cdot, \sdata)|_{W^{0, \infty}([0, 1 - t_0])} \|\partial_t f^*(y, t)\|
    \\&
    \notag
    \lesssim \sdata^{-6}(1 + \|y\|) .
\end{align}
In addition, it holds that
\begin{align*}
    \|\nabla \partial_t b^*(y, t)\|
    &\leq |c_1(\cdot, \sdata)|_{W^{1, \infty}([0, 1 - t_0])} + |c_2(\cdot, \sdata)|_{W^{1, \infty}([0, 1 - t_0])}\|\nabla f^*(y, t)\|
    \\&\quad
    + |c_2(\cdot, \sdata)|_{W^{0, \infty}([0, 1 - t_0])}\|\partial_t \nabla f^*(y, t)\|
    \\&
    \lesssim \sdata^{-8}(1 + \|y\|) .
\end{align*}
Finally, we find that
\begin{align*}
    \|\partial^2_{tt} b^*(y, t)\|
    &\lesssim |c_1|_{W^{2, \infty}} \|y\| + |c_2|_{W^{2, \infty}} + |c_2|_{W^{1, \infty}} \|\partial_t f^*(y, t)\| + |c_2|_{W^{0, \infty}} \|\partial^2_{tt} f^*\|
    \\&
    \lesssim \sdata^{-10}(1 + \|y\|)^2 .
\end{align*}

\noindent
\textbf{Uniform in time approximation of the velocity field.}\quad
We now aim to extend the approximation of $b^*$ in such a way that it holds uniformly for all $t \in [0, 1 - t_0]$.
For this purpose, we employ cubic B-spline quasi-interpolants, whose approximation properties are detailed in the following lemma.

\begin{Lem}[approximation properties of B-splines]
\label{lem:b_splines_approx}
    Let $f \in C^2([0, 1 - t_0])$, where $t_0 \in [0, 1)$ is arbitrary.
    For some $K \in \N$ consider a uniform grid $\tau_k = kh$ for all $k \in \{0, 1, \dots, K\}$, where $h = (1 - t_0) / K$.
    Define the B-spline quasi-interpolant as
    \begin{align*}
        \tilde{f}(t) = \sum_{k = -1}^{K + 1} \alpha_k B_3(h^{-1} t - (k - 2)) ,
        \quad t \in \R ,
    \end{align*}
    where $\alpha_k = f(\tau_k)$ for every $0 \leq k \leq K$, and the boundary coefficients are $\alpha_{-1} = 2f(0) - f(\tau_1)$ and $\alpha_{K + 1} = 2 f(\tau_K) - f(\tau_{K - 1})$.
    Then it holds that
    \begin{align*}
        (i) &\quad |\tilde f - f|_{W^{0, \infty}([0, 1 - t_0])} \leq 4 |f|_{W^{2, \infty}([0, 1 - t_0])} h^2,
        \\
        (ii) &\quad |\tilde f - f|_{W^{1, \infty}([0, 1  - t_0])} \leq 8 |f|_{W^{2, \infty}([0, 1 - t_0])} h .
    \end{align*}

\end{Lem}

We move the proof of Lemma \ref{lem:b_splines_approx} to Appendix \ref{sec:lem_b_splines_approx_proof}.
Recall that the quasi-interpolant $\tilde f(t)$ in Lemma \ref{lem:b_splines_approx} coincides with $I(t;f(0),f(\tau_1),\dots,f(\tau_K))$, where $I(t;\cdot)$ is defined in Definition \ref{def:b_spline_interp}.
Introduce
\begin{align*}
    \tilde f(y, t) = \sum_{k = -1}^{N + 1} \alpha_k B_3(h^{-1}t - (k - 2)),
    \quad (t, y) \in [0, 1 - t_0] \times \R^D ,
\end{align*}
where $\alpha_k = \tilde f_k(y)$ for every $k \in \{0, 1, \dots, N\}$, and the boundary coefficients are $\alpha_{-1} = 2\tilde f_0(y) - \tilde f_1(y)$ and $\alpha_{N + 1} = 2 \tilde f_N(y) - \tilde f_{N - 1}(y)$.
In the following analysis, we utilize several key properties of B-splines, which are summarized in the following lemma.

\begin{Lem}[properties of uniform B-splines]
\label{lem:b_splines_basic_prop}
    Let $\{B_n \, : \, n \in \Z_+\}$ be the family of uniform B-splines, defined recursively by
    \begin{align*}
        B_n(x) = \frac{x}{n}B_{n - 1}(x) + \frac{n + 1 - x}{n}B_{n - 1}(x - 1) ,
        \quad n \in \N,
        \quad x \in \R ,
    \end{align*}
    where $B_0(x) = \1[0 < x \leq 1]$.
    Fix an arbitrary $K \in \N$ and $b > 0$, with $h = b / K$.
    Then, the following properties hold for all $x \in [0, b]$ and $y \in \R$:
    \begin{enumerate}
        \item $0 \leq B_n(y) \leq 1$;
        \item $\mathrm{supp}(B_n) = [0, n + 1]$;
        \item $\sum_{k = -1}^{K + 1} B_n(h^{-1} x - (k - 2)) = 1$ for $n\in\{2,3\}$;
        \item $\sum_{k = -1}^{K + 1}(kh - x) B_3(h^{-1}x - (k - 2)) = 0$;
        \item $B_n'(y)=B_{n-1}(y)-B_{n-1}(y-1)$ for $n\geq2$ at every $y$, and for $n=1$ almost everywhere.
    \end{enumerate}
\end{Lem}

The proof of Lemma \ref{lem:b_splines_basic_prop} is deferred to Appendix \ref{sec:lem_b_splines_basic_prop_proof}.
Consequently, defining
\begin{align*}
    \tilde b(y, t) = c_1(t, \sdata)y + c_2(t, \sdata) \tilde f(y, t),
\end{align*}
it follows from \eqref{eq:c_j_der_bounds} and Lemmata \ref{lem:vel_approx} and \ref{lem:b_splines_basic_prop} that, for all $y \in \R^D$,
\begin{align*}
    \| \tilde b(y, t) \|
    &\lesssim \sdata^{-2} \|y\| + \sdata^{-2} \max_{-1 \leq k \leq N + 1}\|\alpha_k\| \sum_{k = -1}^{N + 1}B_3(h^{-1}t - (k - 2))
    \\&
    \lesssim \sdata^{-2} \|y\| + \sdata^{-2} \sqrt{D} .
\end{align*}
Similarly, for all $y \in \R^D$ and $t \in [0, 1 - t_0]$, we have that
\begin{align}
\label{eq:tilde_f_lip_global}
    \|\nabla \tilde f(y, t)\| \lesssim \max_{0 \leq k \leq N} \|\nabla \tilde f_k(y)\| \sum_{k = -1}^{N + 1} B_3(h^{-1}t - (k - 2))
    \lesssim \left(\smin^{-2} D \log(1 / \eps)\right)^{\cO(1)}
\end{align}
and, hence,
\begin{align}
    \label{eq:tilde_b_lip_global}
    \|\nabla \tilde b(y, t)\|
    \lesssim \sdata^{-2} + \sdata^{-2} \cdot \left(\smin^{-2} D \log(1 / \eps)\right)^{\cO(1)}
    \lesssim \left(\smin^{-2} D \log(1 / \eps)\right)^{\cO(1)} .
\end{align}
We now establish the velocity field approximation error.
Specifically, for every $t \in [0, 1 - t_0]$ and every $y \in \R^D$, it holds that
\begin{align*}
    \|\tilde b(y, t) - b^*(y, t)\|
    &\lesssim |c_2(t, \sdata)|\sum_{k = -1}^{N + 1} \|\alpha_k - \alpha_k^*\| B_3(h^{-1}t - (k - 2))
    \\&\quad
    + \sdata^{-2}\left\| \sum_{k = -1}^{N + 1}\alpha_k^* B_3(h^{-1}t - (k - 2)) - f^*(y, t) \right\| ,
\end{align*}
where the coefficients $\alpha_k^*$ correspond to $f^*$, and are defined analogously to $\alpha_k$.
The direct denoiser bound in Lemma~\ref{lem:vel_approx}(i) gives, for every $0\leq k\leq N$ and every $t$,
\begin{align*}
 |c_2(t,\sdata)|\,\|\tilde f_k-f^*(\cdot,\tau_k)\|_{L^2(\sfp_t^*)}
 \lesssim (D\exp\{\smin^{-2}\}\log(1/\eps))^{\cO(1)}\eps^\beta.
\end{align*}
The triangle inequality gives the same bound, up to a fixed factor, for the two extrapolated boundary coefficients.
Combining this result with the B-spline properties detailed in Lemmata \ref{lem:b_splines_approx} and \ref{lem:b_splines_basic_prop}, it follows that
\begin{align*}
    \|\tilde b(\cdot, t) - b^*(\cdot, t)\|_{L^2(\sfp_t^*)}
    &\lesssim (D \exp\{\smin^{-2}\}\log(1 / \eps))^{\cO(1)} \eps^{\beta}
    \\&\quad
    + \sdata^{-2} \E^{1/2}|f^*(X_t, \cdot)|_{W^{2, \infty}([0, 1 - t_0])}^2 h^2 .
\end{align*}
Finally, invoking Lemma \ref{lem:f_star_deriv} and \eqref{eq:beta_t_alpha_t_der}, we arrive at the final bound
\begin{align}
    \label{eq:tilde_b_l2_approx}
    \|\tilde b(\cdot, t) - b^*(\cdot, t)\|_{L^2(\sfp_t^*)}
    \lesssim \left(D \exp\{\smin^{-2}\}\log(1 / \eps)\right)^{\cO(1)}(\eps^{\beta} + h^2) ,
    \quad t \in [0, 1 - t_0] .
\end{align}

\noindent
\textbf{ODE trajectory approximation.}\quad
We employ Heun's method to approximate the ODE trajectory on a uniform grid.
The following lemma quantifies the discretization error.

\begin{Lem}[Heun's Method error bound]
\label{lem:heun_error_bound}
    Let $N \in \N$ and $t_0 \in [0, 1)$ be arbitrary, and let $\tau_k = kh$ for all $0 \leq k \leq N$, where $h = (1 - t_0) / N$.
    Let $f$ and $f^*$ be mappings from $\R^D \times \R \to \R^D$ such that for all $y \in \R^D$ and $t \in [0, 1 - t_0]$ (and $\tau_k$ where applicable), the following hold:
    \begin{itemize}
        \item[(i)] $\|\nabla f(y, \tau_k)\| \vee \|\nabla f^*(y, t)\| \leq M_y$ for all $y \in \cB(0, R)$ for some $R > 0$ and for all $t \in [0, 1 - t_0]$.
        \item[(ii)] $\|f^*(y, t)\| \leq M(1 + \|y\|)$, $\|f(y, \tau_k)\| \leq M\|y\| + M_f$, $\|\partial_t f^*(y, t)\| \leq M_t(1 + \|y\|)$ .
        \item[(iii)] $\|\nabla^2 f^*(y, t)\| \leq M_{yy}$, $\|\nabla \partial_t f^*(y, t)\| \leq M_{ty}(1 + \|y\|)$, $\|\partial^2_{tt} f^*(y, t)\| \leq M_{tt}(1 + \|y\|)^2$ .
    \end{itemize}
    Let also $\{Z(t) \, : \, 0 \leq t \leq 1 - t_0\}$ be the solution to an ODE $\dot{Z}(t) = f^*(Z(t), t)$ with $Z(0) = y$.
    For $0 \leq k \leq N - 1$, consider the Heun discretization scheme:
    \begin{align*}
        Y_{k + 1} = Z_k + h f(Z_k, \tau_k),
        \quad Z_{k + 1} = Z_k + \frac{h}{2}(f(Z_k, \tau_k) + f(Y_{k + 1}, \tau_{k + 1})),
        \quad Z_0 = y .
    \end{align*}
    Assume further that $h$ is sufficiently small such that $h (M_y \vee M) \leq 1 $.
    Then, for all $y \in \R^D$ satisfying
    \begin{align*}
        (1 + \|y\|) e^M \vee \left( 2\exp\{3M /2 \} \left(\|y\| + \frac{3 M_f}{2}\right) + M_f h \right) \leq R ,
    \end{align*}
    it holds that
    \begin{align*}
        &\max_{0 \leq k \leq N} \|Z_k - Z(\tau_k)\|
        \\&
        \leq
        4e^{(3M_y / 2 + 2M ) (1 - t_0)}h^2(M_y(M_y M + M_t) + M_{yy}M^2 + 2M_{ty}M + M_{tt} )(1 + \|y\|)^2
        \\&\quad
        + 4 e^{3M_y(1 - t_0) / 2} h \sum_{k = 0}^N \|f(Z(\tau_k), \tau_k) - f^*(Z(\tau_k), \tau_k)\| .
    \end{align*}
    Furthermore, if we only assume the growth condition specified in $(ii)$:
    \[
        \|f^*(y, t)\| \leq M(1 + \|y\|),
        \; \|f(y, \tau_k)\| \leq M\|y\| + M_f,
        \quad (y, t) \in \R^D \times [0, 1 - t_0], \; k \in \{0, \dots, N\},
    \]
    and also $hM\leq1$, then, for all $y \in \R^D$ and $0 \leq t \leq 1 - t_0$, the following holds:
    \begin{align*}
        \max_{0 \leq k \leq N}\|Z_k\| \leq \exp\{3 M / 2\} \left( \|y\| + \frac{3 M_f}{2} \right),
        \quad \|Z(t)\| \leq (1 + \|y\|) e^{Mt} .
    \end{align*}

\end{Lem}

The proof of Lemma \ref{lem:heun_error_bound} is deferred to Appendix \ref{sec:lem_heun_error_bound_proof}.
We now apply Lemma \ref{lem:vel_approx} at each grid point $\tau_k = kh$ for $k \in \{0, \dots, N\}$, where $N$ will be specified later in the proof.
With the enlarged tube radius chosen above, $\cB(0,R)\subseteq\cK_{\tau_k}$ for all $0\leq k\leq N$.
Therefore, we can invoke Lemma \ref{lem:heun_error_bound} using the velocity field approximations $\{\bar b(\cdot, \tau_k) \, : \, 0 \leq k \leq N\}$ and the true velocity field $b^*(y, t)$.
The constants obtained from our construction are $M_y \asymp \sdata^{-4}$, $M \asymp \sdata^{-2}$, $M_f \asymp \sdata^{-2}\sqrt{D}$, $M_t \asymp \sdata^{-6}$, $M_{yy} \asymp \sdata^{-6}$, $M_{ty} \asymp \sdata^{-8}$, and $M_{tt} \asymp \sdata^{-10}$.
Choose the constant in the exponent defining $R$ larger than the constants in the growth bounds below.
With a sufficiently small constant in the initial radius condition, this ensures that all sequences under consideration remain inside the ball where the local Lipschitz bounds hold.
Therefore, we conclude that, for all $y \in \R^D$ satisfying
\begin{align}
    \label{eq:flow_approx_ball_def}
    \|y\| \lesssim \sqrt{D} \exp\{\cO(\sdata^{-2})\} \left( \log(1 / \eps) + \log(D\sdata^{-2}) \right) ,
\end{align}
it holds that
\begin{align}
\label{eq:Z_k_tk_diff}
    \notag
    &\max_{0 \leq k \leq N}\|Z_k - Z(\tau_k)\|
    \\&
    \lesssim \exp\{\cO(\sdata^{-4})\}(1+\|y\|)^2 h^2 + \exp\{\cO(\sdata^{-4})\} h \sum_{k = 0}^N \|\bar b(Z(\tau_k), \tau_k) - b^*(Z(\tau_k), \tau_k)\| .
\end{align}
Here, $h \lesssim \sdata^4$ and $\dot{Z}(t) = b^*(Z(t), t)$ with $Z_0 = y$.
Furthermore, for all $y \in \R^D$, it follows that
\begin{align}
\label{eq:max_Z_k_bound}
    \max_{0 \leq k \leq N} \|Z_k\| \lesssim \sqrt{D} \exp\{\cO(\sdata^{-2})\} (1 + \|y\|) .
\end{align}

\noindent
\textbf{Interpolation of the ODE trajectory.}\quad
In the previous step, we constructed an approximation of $Z(\tau_k)$ for all $\tau_k \in [0, 1 - t_0]$ with $k \in \{0, \dots, N\}$.
We now extend the approximation to the interval $[0, 1 - t_0]$.
To achieve this, we use B-splines.
Specifically, consider the approximation
\begin{align*}
    \tilde Z(t) = \sum_{k = -1}^{N + 1} \alpha_k(y) B_3(h^{-1} t - (k - 2)),
\end{align*}
where the coefficients are given by $\alpha_k(y) = Z_k$ for all $0 \leq k \leq N$ and $\alpha_{-1}(y) = 2 y - Z_1$, $\alpha_{N + 1}(y) = 2 Z_N - Z_{N - 1}$.
Similarly, we define $\alpha_k^* = Z(\tau_k)$ for every $0 \leq k \leq N$ and $\alpha_{-1}^*(y) = 2y - Z(\tau_1)$, $\alpha_{N + 1}^*(y) = 2 Z(\tau_N) - Z(\tau_{N - 1})$.
Hence, by Lemma \ref{lem:b_splines_approx}, for all $t \in [0, 1 - t_0]$, it holds that
\begin{align}
    \label{eq:tilde_Z_Z_approx_general}
    \notag
    \|\tilde Z(t) - Z(t)\|
    &\leq \sum_{k = -1}^{N + 1} \|\alpha_k(y) - \alpha_k^*(y)\| B_3(h^{-1}t - (k - 2))
    \\&\quad
    + \left\| \sum_{k = -1}^{N + 1}\alpha_k^*(y)B_3(h^{-1}t - (k - 2)) - Z(t) \right\|
    \\
    \notag
    &\lesssim \max_{0 \leq k \leq N} \|Z_k - Z(\tau_k)\| + |Z|_{W^{2, \infty}} h^2 ,
\end{align}
where the last inequality uses partition of unity property of B-splines outlined in Lemma \ref{lem:b_splines_basic_prop}.
Observe that, according to Lemma \ref{lem:heun_error_bound} and the bounds \eqref{eq:true_vel_jac} and \eqref{eq:true_vel_dt}, it holds that
\begin{align}
\label{eq:Z_w2_bound}
    |Z|_{W^{2, \infty}([0, 1 - t_0])}
    \lesssim \exp\{\cO(\sdata^{-2})\} (1 + \|y\|) ,
    \quad \text{for all } y \in \R^D .
\end{align}
Therefore, from \eqref{eq:Z_k_tk_diff} we conclude that, for all $y$ satisfying \eqref{eq:flow_approx_ball_def}, it follows that
\begin{align}
\label{eq:tilde_Z_Z_ball_approx}
    \notag
    |\tilde Z - Z|_{W^{0, \infty}([0, 1 - t_0])}
    &\lesssim \exp\{\cO(\sdata^{-4})\}(1+\|y\|)^2 h^2
    + \exp\{\cO(\sdata^{-2})\} (1 + \|y\|) h^2
    \\&\quad
    + \exp\{\cO(\sdata^{-4})\} h \sum_{k = 0}^N \|\bar b(Z(\tau_k), \tau_k) - b^*(Z(\tau_k), \tau_k)\| .
\end{align}
Furthermore, from Lemma \ref{lem:heun_error_bound} in conjunction with \eqref{eq:max_Z_k_bound}, \eqref{eq:tilde_Z_Z_approx_general}, and \eqref{eq:Z_w2_bound}, we find that
\begin{align}
\label{eq:tilde_Z_Z_lin_grow}
    |\tilde Z - Z|_{W^{0, \infty}([0, 1 - t_0])}
    \lesssim \sqrt{D} \exp\{\cO(\sdata^{-2})\} (1 + \|y\|),
    \quad \text{for all } y \in \R^D .
\end{align}
In addition, we have that $\tilde Z(t)$ also approximates the derivative of $Z(t)$.
Specifically, applying Lemmata \ref{lem:b_splines_approx} and  \ref{lem:b_splines_basic_prop}, we deduce that, for every $t \in [0, 1 - t_0]$,
\begin{align}
\label{eq:tilde_Z_Z_deriv_approx_general}
    \notag
    \|\partial_t \tilde Z(t) - \partial_t Z(t)\|
    &\leq h^{-1} \left\| \sum_{k = -1}^{N + 1}(\alpha_k(y) - \alpha_k^*(y)) (B_2(h^{-1}t - (k - 2)) - B_2(h^{-1}t - (k - 1))) \right\|
    \\&\quad
    + \left\| \partial_t\sum_{k = -1}^{N + 1} \alpha_k^*(y) B_3(h^{-1}t - (k - 2)) - \partial_t Z(t) \right\|
    \\&
    \notag
    \lesssim h^{-1} \max_{0 \leq k \leq N} \|Z_k - Z(\tau_k)\| + |Z|_{W^{2, \infty}([0, 1 - t_0])} h .
\end{align}
Therefore, from the derived bound, \eqref{eq:Z_k_tk_diff}, and \eqref{eq:Z_w2_bound}, we find that, for every $y$ satisfying \eqref{eq:flow_approx_ball_def}, the following holds:
\begin{align}
\label{eq:tilde_Z_Z_der_ball_approx}
    \notag
    |\tilde{Z} - Z|_{W^{1, \infty}([0, 1 - t_0])}
    &\lesssim \exp\{\cO(\sdata^{-4})\} \sum_{k = 0}^N \|\bar b(Z(\tau_k), \tau_k) - b^*(Z(\tau_k), \tau_k)\|
    \\&\quad
    + \exp\{\cO(\sdata^{-4})\}(1+\|y\|)^2 h
    + \exp\{\cO(\sdata^{-2})\} (1 + \|y\|) h .
\end{align}
Finally, from the bounds in \eqref{eq:max_Z_k_bound}, \eqref{eq:Z_w2_bound}, \eqref{eq:tilde_Z_Z_deriv_approx_general}, and Lemma \ref{lem:heun_error_bound}, we conclude that
\begin{align}
\label{eq:tilde_Z_Z_der_lin_grow}
    |\tilde Z - Z|_{W^{1, \infty}([0, 1 - t_0])}
    \lesssim h^{-1} \sqrt{D} \exp\{\cO(\sdata^{-2})\} (1 + \|y\|),
    \quad \text{for all } y \in \R^D .
\end{align}

\noindent
\textbf{Flow approximation in the $L^2$-norm.}\quad
Having derived bounds in Sobolev seminorms, we now derive the approximation guarantees in the $L^2$-norm.
Let $\tilde R$ denote the radius defined by the right-hand side of \eqref{eq:flow_approx_ball_def}.
We introduce the following approximations:
\begin{align*}
    \tilde v_{t, t}(y) &= \tilde b(y, t),
    \quad (t, y) \in [0, 1 - t_0] \times \R^D,
    \\
    \tilde v_{0, t}(y) &= \frac{\tilde Z(t) - y}{t},
    \quad (t, y) \in (0, 1 - t_0] \times \R^D .
\end{align*}
Recall that  $\tilde Z$ depends on the initial point $y$.
Let $\tilde\Psi_{0, t}(y) = y + t \, \tilde v_{0, t}(y) = \tilde Z(t)$.
From \eqref{eq:tilde_Z_Z_ball_approx} and \eqref{eq:tilde_Z_Z_lin_grow} we deduce that, for all $t \in [0, 1 - t_0]$,
\begin{align*}
    &\|\tilde\Psi_{0, t} - \Psi_{0, t}^*\|_{L^2(\sfp_0^*)}
    \\&
    \leq \left( \integral{\cB(0, \tilde R)} \|\tilde \Psi_{0, t}(y) - \Psi^*_{0, t}(y)\|^2 \sfp_0^*(y)  \dd y \right)^{1/2}
    + \left( \integral{\R^D \setminus \cB(0, \tilde{R}) } \|\tilde \Psi_{0, t}(y) - \Psi^*_{0, t}(y)\|^2 \sfp_0^*(y)  \dd y \right)^{1/2}
    \\&
    \lesssim \sqrt{D} \exp\{\cO(\sdata^{-4})\}(1+\tilde{R})^2 h^2 + \exp\{\cO(\sdata^{-4})\} h \sum_{k = 0}^N \| \bar b(\cdot, \tau_k) - b^*(\cdot, \tau_k)\|_{L^2(\sfp_{\tau_k}^*)}
    \\&\quad
    + \sqrt{D} \exp\{\cO(\sdata^{-2})\} \left( \integral{\R^D \setminus \cB(0, \tilde R)} (1 + \|y\|)^2 \, \sfp_0^*(y) \, \dd y \right)^{1 / 2} .
\end{align*}
By Lemma \ref{lem:vel_approx}, we have that
\begin{align*}
    \|\tilde\Psi_{0, t} - \Psi_{0, t}^*\|_{L^2(\sfp_0^*)}
    &\lesssim \sqrt{D} \exp\{\cO(\sdata^{-4})\} (1+\tilde R)^2 h^2 + \exp\{\cO(\smin^{-4})\} (D \log(1 / \eps))^{\cO(1)} \eps^\beta
    \\&\quad
    + \sqrt{D} \exp\{\cO(\sdata^{-2})\} \left( \integral{\R^D \setminus \cB(0, \tilde R)} (1 + \|y\|)^2 \, \sfp_0^*(y) \, \dd y \right)^{1 / 2} .
\end{align*}
Applying the concentration bound for Gaussian random vector (see Lemma \ref{lem:gaus_conc_norm}) together with the Hölder inequality and using \eqref{eq:flow_approx_ball_def}, we arrive at
\begin{align}
    \label{eq:tilde_psi_psi_star_aux}
    \notag
    \|\tilde\Psi_{0, t} - \Psi_{0, t}^*\|_{L^2(\sfp_0^*)}
    &\lesssim \sqrt{D} \exp\{\cO(\sdata^{-4})\} (1+\tilde R)^2 h^2 + \exp\{\cO(\smin^{-4})\} (D \log(1 / \eps))^{\cO(1)} \eps^\beta
    \\&\quad
    + D^{\cO(1)} \exp\{\cO(\sdata^{-2})\} \eps^\beta
    \\&
    \notag
    \lesssim (D \log(1 / \eps))^{\cO(1)} \exp\{\cO(\smin^{-4})\}(\eps^\beta + h^2) .
\end{align}
Similarly, we obtain the bound for the derivatives w.r.t. $t$.
Formally, from \eqref{eq:tilde_Z_Z_der_ball_approx} and \eqref{eq:tilde_Z_Z_der_lin_grow} we have that
\begin{align*}
    &\|\partial_t \tilde\Psi_{0, t} - \partial_t \Psi^*_{0, t}\|_{L^2(\sfp_0^*)}
    \lesssim \sqrt{D} \exp\{\cO(\sdata^{-4})\}(1+\tilde R)^2 h
    \\&
    + \exp\{\cO(\sdata^{-4})\} \sum_{k = 0}^N \| \bar b(\cdot, \tau_k) - b^*(\cdot, \tau_k)\|_{L^2(\sfp_{\tau_k}^*)}
    \\&
    + h^{-1}\sqrt{D}\exp\{\cO(\sdata^{-2})\} \left( \integral{\R^D \setminus \cB(0, \tilde R)} (1 + \|y\|)^2 \sfp_0^*(y) \, \dd y \right)^{1/2} .
\end{align*}
Now let $h \asymp \smin^{2} \eps^{\cO(1)}$, where the exact dependence on $\eps$ will be specified later.
By applying Lemma \ref{lem:vel_approx},  the Gaussian concentration bound from Lemma \ref{lem:gaus_conc_norm}, and the Hölder inequality, we conclude that
\begin{align*}
    \|\partial_t \tilde\Psi_{0, t} - \partial_t \Psi^*_{0, t}\|_{L^2(\sfp_0^*)}
    &\lesssim \sqrt{D} \exp\{\cO(\sdata^{-4})\}(1+\tilde R)^2 h
    + h^{-1}\exp\{\cO(\smin^{-4})\} (D \log(1 / \eps))^{\cO(1)} \eps^\beta
    \\&\quad
    + D^{\cO(1)} \exp\{\cO(\sdata^{-2})\} \eps^\beta .
\end{align*}
Thus, choose $N=\lceil C\smin^{-2}\eps^{-\beta/2}\rceil$. Since $1/2\leq1-t_0\leq1$, we have $h\asymp\smin^2\eps^{\beta/2}$; the smallness condition $\eps^\beta\leq c\smin^4$ ensures $h(M_y\vee M)\leq1$. This gives
\begin{align}
    \label{eq:psi_der_l2_approx}
    \|\partial_t \tilde\Psi_{0, t} - \partial_t \Psi^*_{0, t}\|_{L^2(\sfp_0^*)}
    \lesssim (D \log(1 / \eps))^{\cO(1)} \exp\{\cO(\smin^{-4})\} \eps^{\beta / 2} .
\end{align}
From \eqref{eq:tilde_psi_psi_star_aux} we obtain
\begin{align}
\label{eq:psi_l2_approx}
    \|\tilde\Psi_{0, t} - \Psi_{0, t}^*\|_{L^2(\sfp_0^*)}
    \lesssim (D \log(1 / \eps))^{\cO(1)} \exp\{\cO(\smin^{-4})\} \eps^\beta .
\end{align}
Finally, the choice of $h$ implies that $N \asymp \smin^{-2}\eps^{-\beta / 2}$ in Definition \ref{def:av_vel_class}.

\noindent
\textbf{Approximating the loss function.}\quad
From Proposition \ref{prop:fm_loss_dif} and the fact that $R(v^*) = 0$ we find that
\begin{align*}
    &\E[\ell_{SD}(\tilde v, X_1) - \ell_{SD}(v^*, X_1)]
    = \E[\ell_{FM}(\tilde v, X_1) - \ell_{FM}(v^*, X_1)] + R(\tilde v)
    \\&
    = \integral{0}^{1 - t_0} \E \|\tilde v_{t, t}(X_t) - b^*(X_t, t)\|^2 \, \dd t
    + \integral{0}^{1 - t_0} \E \|\partial_t \tilde\Psi_{0, t}(X_0) - \tilde v_{t, t}(\tilde \Psi_{0, t}(X_0))\|^2 \, \dd t .
\end{align*}
Taking into account that $\tilde v_{t, t}$ is globally Lipschitz, as established by \eqref{eq:tilde_b_lip_global}, we deduce that the last term on the right-hand side of the above expression is bounded by
\begin{align*}
    &\integral{0}^{1 - t_0} \E \|\partial_t \tilde\Psi_{0, t}(X_0) - \tilde v_{t, t}(\tilde \Psi_{0, t}(X_0))\|^2 \, \dd t
    \lesssim \integral{0}^{1 - t_0} \E \|\partial_t \tilde\Psi_{0, t}(X_0) - \partial_t \Psi^*_{0, t}(X_0)\|^2 \, \dd t
    \\&
    + \left(D \sdata^{-2} \log(1 / \eps)\right)^{\cO(1)} \integral{0}^{1 - t_0} \E \|\tilde \Psi_{0, t}(X_0) - \Psi^*_{0, t}(X_0)\|^2 \, \dd t
    + \integral{0}^{1 - t_0} \E \| \tilde v_{t, t}(X_t) - b^*(X_t, t) \|^2 \, \dd t .
\end{align*}
By applying the bounds \eqref{eq:tilde_b_l2_approx}, \eqref{eq:psi_der_l2_approx}, \eqref{eq:psi_l2_approx}, and considering the choice of $h$, we conclude that
\begin{align*}
    \E[\ell_{SD}(\tilde v, X_1) - \ell_{SD}(v^*, X_1)]
    \lesssim (D \log(1 / \eps))^{\cO(1)} \exp\{\cO(\smin^{-4})\} \eps^\beta .
\end{align*}
The proof is now complete.

\myendproof

\subsection{Proof of Lemma \ref{lem:vel_approx}}
\label{sec:lem_vel_approx_proof}

Notice that $X_t \sim \sfp_t^*$ satisfies Assumption 3.1 in
\cite{yakovlev2026simultaneous} with variance parameter
$\tilde\alpha_t^2=\alpha_t^2+\beta_t^2\sdata^2$.
Specifically, the distribution of $X_t$ coincides with that of
$\beta_t g^*(U)+\tilde\alpha_t Z$, where
$U\sim\mathrm{Un}([0,1]^d)$ and $Z\sim\cN(0,I_D)$ are independent.
Furthermore, $0\leq\beta_t\leq1$ and $\sdata^2/2\leq\tilde\alpha_t^2\leq1$ by
Assumption \ref{asn:flow_matching_coefs}.

Set the derivative order $m=1$ in the construction of \cite[Theorem 3.2]{yakovlev2026simultaneous}.
Let $g^\circ$ be the piecewise Taylor approximation of $g^*$ from Lemma C.1 of that work.
Then
\[
    \delta_g
    :=\|g^\circ-g^*\|_{L^\infty([0,1]^d)}
    \lesssim \sqrt{D}\eps^\beta,
    \quad
    \|g^\circ\|_{L^\infty([0,1]^d)}\leq2,
\]
provided that $\eps$ is sufficiently small.
Define
\[
    f^\circ(y)
    =
    \frac{
        \integral{[0,1]^d}
        g^\circ(u)
        \exp\left(
            -\frac{\|y-\beta_tg^\circ(u)\|^2}
            {2\tilde\alpha_t^2}
        \right)\,\dd u
    }{
        \integral{[0,1]^d}
        \exp\left(
            -\frac{\|y-\beta_tg^\circ(u)\|^2}
            {2\tilde\alpha_t^2}
        \right)\,\dd u
    }.
\]
We emphasize that the numerator contains $g^\circ$, not $\beta_tg^\circ$.
To obtain this function, apply the construction in
Lemmata 4.3 and 4.4 of \cite{yakovlev2026simultaneous}
with $\beta_tg^\circ$.
Choose
\begin{align}
\label{eq:R_def}
    R
    \asymp
    \exp\{\cO(\smin^{-2})\}D
    \left(\log(1/\eps)+\log(D\smin^{-2})\right),
\end{align}
and set
$\cK_\infty=[-R_\infty,R_\infty]^D$, where
$R_\infty=(R+2)\vee(5/2)$.
In particular,
$\cK_t\subseteq\cB(0,R+1)\subseteq\cK_\infty$.
The proof of Theorem 3.2 and Lemma 4.4 of
\cite{yakovlev2026simultaneous} yield a GELU network $\bar f \in \NN(L, W, S, B)$ with architecture
\begin{align*}
    L
    &\lesssim
    \log(DR\smin^{-2}\log(1/\eps))
    \lesssim
    \smin^{-2}\log(D\smin^{-2}\log(1/\eps)),
    \\
    \|W\|_\infty\vee S
    &\lesssim
    \eps^{-d}
    \left(
        DR\smin^{-2}\log(1/\eps)\log(DR\smin^{-2})
    \right)^{\cO(1)}
    \lesssim
    \eps^{-d}\exp\{\cO(\smin^{-2})\}
    (D\log(1/\eps))^{\cO(1)},
    \\
    \log B
    &\lesssim
    \left(
        R\smin^{-2}\log(1/\eps)\log(DR\smin^{-2})
    \right)^{\cO(1)}
    \lesssim
    \exp\{\cO(\smin^{-2})\}
    (D\log(1/\eps))^{\cO(1)}.
\end{align*}
Moreover, for every $1\leq l\leq D$,
\begin{align}
\label{eq:bar_f_f_circ_prox}
    \notag
    &\max_{j\in\{0,1\}}
    |\bar f_l-f_l^\circ|_{W^{j,\infty}(\cK_\infty)}
    \leq\eps^\beta,
    \\
    &\|\bar f_l\|_{L^\infty(\R^D)}\lesssim1,
    \\
    \notag
    &|\bar f_l|_{W^{1,\infty}(\R^D)}
    \lesssim
    \sdata^{-2}
    \left(
        D\log(1/\eps)\log(R\sdata^{-2})
    \right)^{\cO(1)}.
\end{align}
Here the first inequality controls both the function and
its first derivatives.
Differentiating the normalized integral defining $f^\circ$ implies that
\[
    \nabla f^\circ(y)
    =
    \frac{\beta_t}{\tilde\alpha_t^2}
    \operatorname{Cov}_{y,t}(g^\circ(U)),
\]
where the covariance is taken under the probability measure
on $[0,1]^d$ whose density is proportional to
$\exp(-\|y-\beta_tg^\circ(u)\|^2/(2\tilde\alpha_t^2))$.
Consequently, $\|\nabla f^\circ(y)\|\lesssim\sdata^{-2}$ for all
$y\in\R^D$.
Taking into account that $\cK_t\subseteq\cK_\infty$,
we obtain
\begin{align*}
    \sup_{y\in\cK_t}\|\nabla\bar f(y)\|
    &\leq
    D\max_{1\leq l\leq D}
    |\bar f_l-f_l^\circ|_{W^{1,\infty}(\cK_t)}
    +
    \sup_{y\in\cK_t}\|\nabla f^\circ(y)\|
    \\
    &\lesssim D\eps^\beta+\sdata^{-2}
    \lesssim\sdata^{-2}.
\end{align*}
The last inequality follows from
\eqref{eq:vel_appr_eps_constr}.
Furthermore, \eqref{eq:bar_f_f_circ_prox} yields
$\|\bar f\|_{L^\infty(\R^D)}\lesssim\sqrt D$.
Combining \eqref{eq:R_def} and
\eqref{eq:bar_f_f_circ_prox} also gives
\[
    \max_{1\leq l\leq D}
    |\bar f_l|_{W^{1,\infty}(\R^D)}
    \lesssim
    \left(\smin^{-2}D\log(1/\eps)\right)^{\cO(1)}.
\]
This establishes claim $(ii)$.

We now prove claim $(i)$.
For $\|y\|\leq R+1$, the logarithms of the Gaussian weights corresponding to $g^\circ$ and $g^*$ satisfy
\[
    \sup_{u\in[0,1]^d}
    \left|
        \frac{
            \|y-\beta_tg^\circ(u)\|^2
            -
            \|y-\beta_tg^*(u)\|^2
        }{2\tilde\alpha_t^2}
    \right|
    \lesssim \sdata^{-2}(R+3)\delta_g.
\]
Let $\eta = C\sdata^{-2}(R+3)\delta_g$.
If $\eta\leq1$, the ratio of the corresponding normalized
densities lies between $e^{-2\eta}$ and $e^{2\eta}$.
Subtracting the normalized integrals therefore gives
\begin{align*}
    \sup_{\|y\|\leq R+1}
    \|f^\circ(y)-f^*(y,t)\|
    \lesssim \delta_g + \eta
    \lesssim
    \sqrt D\,\sdata^{-2}(R+3)\eps^\beta.
\end{align*}
Together with \eqref{eq:bar_f_f_circ_prox}, this implies
\[
    \sup_{\|y\|\leq R+1}
    \|\bar f(y)-f^*(y,t)\|
    \lesssim
    \sqrt D\,\sdata^{-2}(R+3)\eps^\beta.
\]

For every $s\in[0,1-t_0]$, the distribution of $X_s$
coincides with that of
$\beta_sg^*(U)+\tilde\alpha_s Z$, with
$\|\beta_sg^*(U)\|\leq1$ and $\tilde\alpha_s\leq1$.
Hence, by Lemma \ref{lem:gaus_conc_norm} and a sufficiently
large constant in \eqref{eq:R_def},
\[
    \sup_{0\leq s\leq1-t_0}
    \p(\|X_s\|>R+1)
    \leq \p(\|Z\|>R)
    \lesssim\eps^{2\beta}.
\]
Since $\bar f$ and $f^*(\cdot,t)$ are globally bounded
by $C\sqrt D$, splitting the integral into
$\cB(0,R+1)$ and its complement yields
\begin{align*}
    \sup_{0\leq s\leq1-t_0}
    \|\bar f-f^*(\cdot,t)\|^2_{L^2(\sfp_s^*)}
    &\lesssim
    D\sdata^{-4}(R+3)^2\eps^{2\beta}
    +D\eps^{2\beta}
    \\&
    \lesssim
    \left(
        D\exp(\smin^{-2})\log(1/\eps)
    \right)^{\cO(1)}\eps^{2\beta}.
\end{align*}
Finally,
\[
    \bar b(y)-b^*(y,t)
    =
    c_2(t,\sdata)\bigl(\bar f(y)-f^*(y,t)\bigr),
\]
and $|c_2(t,\sdata)|\lesssim\sdata^{-2}$.
This proves $(i)$.

It remains to check the smallness requirements.
The assumptions of \cite[Lemma 4.3]{yakovlev2026simultaneous} are ensured by
\[
    \sqrt D(R+3)(\eps^\beta\vee\eps)
    \lesssim \sdata^2.
\]
Substituting \eqref{eq:R_def}, this follows from \eqref{eq:vel_appr_eps_constr}.
The proof is complete.

\myendproof

\subsection{Proof of Lemma \ref{lem:f_star_deriv}}
\label{sec:lem_f_star_deriv_proof}

Fix arbitrary $y\in\R^D$ and $t\in[0,1-t_0]$.
It follows immediately from \eqref{eq:f_star_def} that $\|f^*(y,t)\|\leq1$.
For brevity, write
\[
    a(t)=\frac{\beta_t}{\tilde\alpha_t^2},
    \quad
    q(t)=\frac{\beta_t^2}{\tilde\alpha_t^2},
\]
and introduce
\[
    w(u)
    =
    \exp\left(
        a(t)y^\top g^*(u)
        -\frac{q(t)}2\|g^*(u)\|^2
    \right),
    \quad
    J=\integral{[0,1]^d}w(u)\,\dd u.
\]
The dependence of $w$ and $J$ on $(y,t)$ is suppressed.
Then
\begin{align}
\label{eq:f_star_another_def}
    f^*(y,t)
    =
    \frac{1}{J}
    \integral{[0,1]^d}g^*(u)w(u)\,\dd u .
\end{align}
For products of integration variables, let us write $g_j=g^*(u_j)$ and $w_j=w(u_j)$.
Taking the gradient with respect to $y$ yields
\[
    \nabla f^*(y,t)
    =
    \frac{a(t)}{J^2}
    \integral{([0,1]^d)^2}
    P_y(u_1,u_2)w_1w_2\,\dd u_1\dd u_2,
\]
where $P_y(u_1,u_2) = g_1\otimes(g_1-g_2)$.
To obtain the constant in (i), observe that
\[
    \frac{1}{J^2}
    \integral{([0,1]^d)^2} P_y(u_1,u_2)w_1w_2\,\dd u_1\dd u_2
    =
    \frac{1}{J}
    \integral{[0,1]^d}
    g^*(u)\otimes g^*(u)\,w(u)\,\dd u
    -
    f^*(y,t)\otimes f^*(y,t).
\]
This is a covariance matrix under the probability density
$w/J$.
Its operator norm is at most one, since $\|g^*(u)\|\leq1$.
Consequently, $\|\nabla f^*(y,t)\| \leq a(t) \leq 2\sdata^{-2}$,
where the last inequality follows from
$\beta_t\leq1$ and $\tilde\alpha_t^2\geq\sdata^2/2$.
This proves (i).

Taking another derivative with respect to $y$, we obtain
\[
    \nabla^2 f^*(y,t)
    =
    \frac{a(t)^2}{J^3}
    \integral{([0,1]^d)^3}
    P_{yy}(u_1,u_2,u_3)
    w_1w_2w_3\,\dd u_1\dd u_2\dd u_3,
\]
where $P_{yy}(u_1,u_2,u_3) = P_y(u_1,u_2)\otimes(g_1+g_2-2g_3)$.
Using $\|P_y\|\leq2$ and $\|g_1+g_2-2g_3\|\leq4$, we find that
\[
    \|\nabla^2f^*(y,t)\|
    \leq 8a(t)^2
    \leq32\sdata^{-4},
\]
which proves (iii).

By differentiating \eqref{eq:f_star_another_def}
with respect to $t$, we obtain
\begin{align}
\label{eq:d_t_f_star}
    \partial_t f^*(y,t)
    =
    \frac{
        \integral{([0,1]^d)^2}
        P_t(u_1,u_2)w_1w_2\,\dd u_1\dd u_2
    }{J^2},
\end{align}
where
\begin{align}
\label{eq:deriv_pt_def}
    P_t(u_1,u_2)
    =
    a'(t)g_1y^\top(g_1-g_2)
    -\frac{q'(t)}2g_1
    \bigl(\|g_1\|^2-\|g_2\|^2\bigr).
\end{align}
Set
\[
    Q(y,t)
    =
    \|y\|\,|a'(t)|+\frac12|q'(t)|.
\]
Since $\|g_j\|\leq1$, we have
$\|P_t\|\leq2Q(y,t)$ and therefore
\[
    \|\partial_t f^*(y,t)\|
    \leq2Q(y,t)
    =
    2\|y\|\,|a'(t)|+|q'(t)|.
\]
Hence, claim $(ii)$ readily follows.

We next differentiate \eqref{eq:d_t_f_star} with respect to $y$:
\[
    \nabla\partial_t f^*(y,t)
    =
    \frac{1}{J^3}
    \integral{([0,1]^d)^3}
    P_{y,t}(u_1,u_2,u_3)
    w_1w_2w_3\,\dd u_1\dd u_2\dd u_3 ,
\]
where
\[
    P_{y,t}(u_1,u_2,u_3)
    =
    a(t)P_t(u_1,u_2)\otimes(g_1+g_2-2g_3)
    +
    \nabla_yP_t(u_1,u_2).
\]
In particular,
\[
    \nabla_yP_t(u_1,u_2)
    =
    a'(t)g_1\otimes(g_1-g_2),
    \quad
    \|\nabla_yP_t(u_1,u_2)\|\leq2|a'(t)|.
\]
Hence,
\begin{align*}
    \|\nabla\partial_t f^*(y,t)\|
    &\leq
    4a(t)\|P_t\|_{L^\infty(([0,1]^d)^2)}
    +2|a'(t)|
    \\
    &\leq
    8a(t)Q(y,t)+2|a'(t)|
    \\
    &\leq
    16\sdata^{-2}
    \left(
        \|y\|\,|a'(t)|+\frac12|q'(t)|
    \right)
    +2|a'(t)|.
\end{align*}
Therefore, claim $(iv)$ readily follows.

Finally, define
\[
    \phi_j
    =
    a'(t)y^\top g_j-\frac{q'(t)}2\|g_j\|^2,
    \quad |\phi_j|\leq Q(y,t).
\]
Differentiating \eqref{eq:d_t_f_star} once more in $t$ yields
\[
    \partial_{tt}f^*(y,t)
    =
    \frac{1}{J^3}
    \integral{([0,1]^d)^3}
    P_{tt}(u_1,u_2,u_3)
    w_1w_2w_3\,\dd u_1\dd u_2\dd u_3,
\]
where
\[
    P_{tt}(u_1,u_2,u_3)
    =
    \partial_tP_t(u_1,u_2)
    +
    P_t(u_1,u_2)(\phi_1+\phi_2-2\phi_3).
\]
From \eqref{eq:deriv_pt_def} we find that
\[
    \|\partial_tP_t(u_1,u_2)\|
    \leq
    2\|y\|\,|a''(t)|+|q''(t)|.
\]
Furthermore, $|\phi_1+\phi_2-2\phi_3|\leq4Q(y,t)$.
Consequently,
\begin{align*}
    \|\partial_{tt}f^*(y,t)\|
    &\leq
    2\|y\|\,|a''(t)|+|q''(t)|
    +4Q(y,t)\|P_t\|_{L^\infty(([0,1]^d)^2)}
    \\
    &\leq
    2\|y\|\,|a''(t)|+|q''(t)|
    +8\left(
        \|y\|\,|a'(t)|+\frac12|q'(t)|
    \right)^2.
\end{align*}
Claim $(v)$ follows immediately.
The proof is complete.

\myendproof

\subsection{Proof of Lemma \ref{lem:heun_error_bound}}
\label{sec:lem_heun_error_bound_proof}

For clarity, we divide the proof into several steps.

\noindent
\textbf{Bounding $\|Z_{k}\|$.\quad}
Let us first establish the last claim of the lemma.
For all $0 \leq k \leq N - 1$, it holds that
\begin{align*}
    \|Z_{k + 1}\|
    &\leq \|Z_k\| + \frac{h}{2} \left( 2M_f + M \|Z_k\| + M\|Y_{k + 1}\| \right)
    \\&
    \leq \|Z_k\| + \frac{h}{2} \left( 2M_f + M \|Z_k\| + M (\|Z_k\| + h(M_f + M \|Z_k\|) ) \right)
    \\&
    = \left(1 + Mh + \frac{M^2 h^2}{2}\right) \|Z_k\| + M_f h + \frac{M M_f h^2}{2} .
\end{align*}
Using the assumption that $Mh \leq 1$ and applying the discrete Grönwall lemma, we arrive at
\begin{align*}
    \max_{0 \leq k \leq N}\|Z_k\|
    \leq \left(1 + \frac{3 M h}{2}\right)^N \|y\| + \frac{3 M_f}{2} \left(1 + \frac{3 M h}{2}\right)^N
    \leq \exp\{3 M / 2\} \left( \|y\| + \frac{3 M_f}{2} \right) .
\end{align*}
In addition, we have that, for all $0 \leq k \leq N - 1$,
\begin{align*}
    \|Y_{k + 1}\| \leq \|Z_k\| + h \|f(Z_k, \tau_k)\|
    &\leq (1 + Mh)\|Z_k\| + M_f h
    \\&
    \leq 2 \exp\{3M / 2\}\left(\|y\| + \frac{3M_f}{2}\right) + M_f h .
\end{align*}
We also note that, for all $0 \leq t \leq 1 - t_0$ and $y \in \R^D$,
\begin{align*}
    \|Z(t)\| \leq \|y\| + \integral{0}^t \|f^*(Z(s), s)\| \, \dd s
    \leq \|y\| + M \integral{0}^t  \dd s + M \integral{0}^t \|Z(s)\| \, \dd s.
\end{align*}
Hence, the Grönwall lemma yields that
\begin{align}
\label{eq:ode_heun_sol_growth}
    \|Z(t)\| \leq (1 + \|y\|) e^{Mt},
    \quad \text{for all } y \in \R^D,
    \; 0 \leq t \leq 1 - t_0 ,
\end{align}
and the last claim follows.

\noindent
\textbf{ODE solver error bound.}\quad
Consider $y \in \R^D$ satisfying the inequality
\begin{align*}
    (1 + \|y\|) e^M \vee \left( 2\exp\{3M /2 \} \left(\|y\| + \frac{3 M_f}{2}\right) + M_f h \right) \leq R .
\end{align*}
Therefore, for all such $y$ we have that $Z(t) \in \cB(0, R)$ for all $0 \leq t \leq 1 - t_0$, $Z_k \in \cB(0, R)$ for all $0 \leq k \leq N$, and $Y_k \in \cB(0, R)$ for every $1 \leq k \leq N$.
We now proceed with the approximation error bound for the ODE solver.
For every $0 \leq k \leq N - 1$ we have that
\begin{align}
\label{eq:heun_error_decomp}
    \notag
    \|Z_{k + 1} - Z(\tau_{k + 1})\|
    &\leq \|Z_k - Z(\tau_k)\| + \frac{h}{2}\left( M_y\|Z_k - Z(\tau_k)\| + \|f(Z(\tau_k), \tau_k) - f^*(Z(\tau_k), \tau_k)\| \right)
    \\&\quad
    \notag
    + \frac{h}{2}\left( M_y \|Y_{k + 1} - Z(\tau_{k + 1})\| + \|f(Z(\tau_{k + 1}), \tau_{k + 1}) - f^*(Z(\tau_{k + 1}), \tau_{k + 1})\| \right)
    \\&\quad
    + \left\| \integral{\tau_k}^{\tau_{k + 1}} f^*(Z(\tau), \tau) \dd \tau - \frac{h}{2}(f^*(Z(\tau_k), \tau_k) + f^*(Z(\tau_{k + 1}), \tau_{k + 1})) \right\| .
\end{align}
Next, notice that the function
\begin{align*}
    g : u \mapsto \integral{\tau_k}^{\tau_k + u} f^*(Z(\tau), \tau) \dd \tau - \frac{u}{2}(f^*(Z(\tau_k), \tau_k) + f^*(Z(\tau_k + u), \tau_k + u))
\end{align*}
satisfies $g(0) = 0$ and $g'(0) = 0$.
Furthermore, its second derivative is
\begin{align*}
    \frac{\dd^2 g(u)}{\dd u^2}
    &= \frac{\dd}{\dd u} \left[ \frac{1}{2}(f^*(Z(\tau_k + u), \tau_k + u) - f^*(Z(\tau_k), \tau_k)) - \frac{u}{2} \frac{\dd f^*(Z(\tau_k + u), \tau_k + u)}{\dd u} \right]
    \\&
    = - \frac{u}{2} \frac{\dd^2 f^*(Z(\tau_k + u), \tau_k + u)}{\dd u^2} .
\end{align*}
Therefore, by Taylor's theorem, we have that
\begin{align*}
    &\left\| \integral{\tau_k}^{\tau_{k + 1}} f^*(Z(\tau), \tau) \dd \tau - \frac{h}{2}(f^*(Z(\tau_k), \tau_k) + f^*(Z(\tau_{k + 1}), \tau_{k + 1})) \right\|
    \\&
    \leq \frac{h^3}{4}\max_{0 \leq u \leq h} \left\| \frac{\dd^2 f^*(Z(\tau_k + u), \tau_k + u)}{\dd u^2} \right\| .
\end{align*}
We now observe that the second derivative satisfies
\begin{align*}
    \left\| \frac{\dd^2 f^*(Z(\tau_k + u), \tau_k + u)}{\dd u^2} \right\|
    &= \left\| \nabla^2 f^*[f^*,f^*]+2\nabla\partial_t f^*\,f^*+\partial_{tt}f^*+\nabla f^*(\nabla f^*\,f^*+\partial_t f^*)\right\|
    \\&
    \leq (M_y(M_yM+M_t)+M_{yy}M^2+2M_{ty}M+M_{tt})(1+\|Z(\tau_k+u)\|)^2 ,
\end{align*}
where we abbreviate $f^*$ for $f^*(Z(\tau_k + u), \tau_k + u)$.
Hence, we obtain from \eqref{eq:ode_heun_sol_growth} that
\begin{align*}
  &\max_{0 \leq u \leq h} \left\| \frac{\dd^2 f^*(Z(\tau_k + u), \tau_k + u)}{\dd u^2} \right\|
  \\&
  \leq 4 (M_y(M_yM+M_t)+M_{yy}M^2+2M_{ty}M+M_{tt})(1+\|y\|)^2 e^{2M (1 - t_0)} .
\end{align*}
Thus, it follows that
\begin{align}
    \label{eq:heun_int_bound}
    &\left\| \integral{\tau_k}^{\tau_{k + 1}} f^*(Z(\tau), \tau) \dd \tau - \frac{h}{2}(f^*(Z(\tau_k), \tau_k) + f^*(Z(\tau_{k + 1}), \tau_{k + 1})) \right\|
    \\&
    \notag
    \leq h^3 (M_y(M_yM+M_t)+M_{yy}M^2+2M_{ty}M+M_{tt})(1+\|y\|)^2 e^{2M (1 - t_0)} .
\end{align}
Following this, we bound the error corresponding to the auxiliary sequence.
For all $0 \leq k \leq N - 1$, we obtain
\begin{align}
\label{eq:heun_Y_kp1_bound}
    \notag
    \|Y_{k + 1} - Z(\tau_{k + 1})\|
    &\leq \|Z_k - Z(\tau_k)\| + M_y h \|Z_k - Z(\tau_k)\| + h \|f^*(Z(\tau_k), \tau_k) - f(Z(\tau_k), \tau_k)\|
    \\&\quad
    + \left\|\integral{\tau_k}^{\tau_{k + 1}}[f^*(Z(\tau), \tau) - f^*(Z(\tau_k), \tau_k)] \, \dd \tau \right\| .
\end{align}
The last term is evaluated as
\begin{align*}
    &\left\|\integral{\tau_k}^{\tau_{k + 1}}[f^*(Z(\tau), \tau) - f^*(Z(\tau_k), \tau_k)] \, \dd \tau \right\|
    \\&
    \leq \integral{\tau_k}^{\tau_{k + 1}} \integral{\tau_k}^\tau \left\|\frac{\dd f^*(Z(u), u)}{\dd u} \right\| \, \dd u \, \dd \tau
    \\&
    \leq \frac{h^2}{2} \max_{ u \in [\tau_k, \tau_{k+1}] }\|\nabla f^*(Z(u), u) f^*(Z(u), u) + \partial_t f^*(Z(u), u)\|
    \\&
    \leq \frac{h^2}{2}(M_yM+M_t)
    \left(1+\max_{u\in[\tau_k,\tau_{k+1}]}\|Z(u)\|\right).
\end{align*}
Using the bound \eqref{eq:ode_heun_sol_growth}, we arrive at
\begin{align*}
    \left\|\integral{\tau_k}^{\tau_{k + 1}}[f^*(Z(\tau), \tau) - f^*(Z(\tau_k), \tau_k)] \, \dd \tau \right\|
    \leq h^2 (M_y M + M_t)e^{M(1 - t_0)}(1 + \|y\|) .
\end{align*}
Therefore, combining this bound with \eqref{eq:heun_error_decomp}, \eqref{eq:heun_int_bound}, and \eqref{eq:heun_Y_kp1_bound}, we conclude that
\begin{align*}
    &\|Z_{k + 1} - Z(\tau_{k + 1})\|
    \\&
    \leq \left(1 + \frac{3}{2} hM_y \right) \|Z_k - Z(\tau_k)\|
    + h^3 M_y (M_y M + M_t)(1 + \|y\|)e^{M(1 - t_0)}
    \\&\quad
    + h(\|f(Z(\tau_k), \tau_k) - f^*(Z(\tau_k), \tau_k)\| + \|f(Z(\tau_{k + 1}), \tau_{k + 1}) - f^*(Z(\tau_{k + 1}), \tau_{k + 1})\| )
    \\&\quad
    + h^3 (M_y(M_yM+M_t)+M_{yy}M^2+2M_{ty}M+M_{tt})(1+\|y\|)^2 e^{2M (1 - t_0)} ,
\end{align*}
where we used $h M_y \leq 1$.
Thus, applying the discrete Grönwall lemma yields
\begin{align*}
    &\max_{0 \leq k \leq N} \|Z_k - Z(\tau_k)\|
    \\&
    \leq
    4e^{3M_y (1 - t_0) / 2}h^2(M_y(M_y M + M_t) + M_{yy}M^2 + 2M_{ty}M + M_{tt} )(1 + \|y\|)^2 e^{2M(1 - t_0)}
    \\&\quad
    + 4 e^{3M_y(1 - t_0) / 2} h \sum_{k = 0}^N \|f(Z(\tau_k), \tau_k) - f^*(Z(\tau_k), \tau_k)\| .
\end{align*}
The desired claim now follows immediately.

\myendproof

\subsection{Proof of Lemma \ref{lem:b_splines_basic_prop}}
\label{sec:lem_b_splines_basic_prop_proof}

Note that properties 1, 2, 3, and 5 follow directly from standard results in \cite{phillips2003interpolation}:
\begin{itemize}
    \item Property 1 is established by Problem 6.2.3;
    \item Property 2 follows from Theorem 6.2.1;
    \item Property 3 follows from Theorem 6.2.5;
    \item Property 5 follows from Theorem 6.2.2;
\end{itemize}
We now prove property 4.
First, observe that $x \mapsto B_3(x + 2)$ is symmetric about zero.
Consequently, at $x = 0$,
\begin{align*}
    \sum_{k = -1}^{K + 1} kh \, B_3( - (k - 2))
    = \sum_{k \in \Z} kh \, B_3(-(k - 2))
    = 0 .
\end{align*}
To extend this to all $x \in [0, b]$, it suffices to show that, for all such $x$,
\begin{align*}
    \frac{\dd}{\dd x} \left[ \sum_{k = -1}^{K + 1}(kh - x) B_3(h^{-1}x - (k - 2)) \right] = 0 .
\end{align*}
By property 5, we have that
\begin{align*}
     &\frac{\dd}{\dd x} \left[ \sum_{k = -1}^{K + 1}(kh - x) B_3(h^{-1}x - (k - 2)) \right]
     = -\sum_{k = -1}^{K + 1}B_3(h^{-1}x - (k - 2))
     \\&
     + h^{-1} \sum_{k = -1}^{K + 1}(kh - x) \left(B_2(h^{-1}x - (k - 2)) - B_2(h^{-1}x - (k - 1)) \right) .
\end{align*}
Applying property 3 and rearranging the terms, we arrive at
\begin{align*}
    &\frac{\dd}{\dd x} \left[ \sum_{k = -1}^{K + 1}(kh - x) B_3(h^{-1}x - (k - 2)) \right]
    = h^{-1}(-h - x)B_2(h^{-1}x + 3)
    \\&
    - h^{-1} (Kh + h - x) B_2(h^{-1}x - K) - B_2(h^{-1}x + 3) .
\end{align*}
Using property 2, we deduce that the derivative is zero on $[0, b]$, proving property 4.
The proof is finished.

\myendproof

\subsection{Proof of Lemma \ref{lem:b_splines_approx}}
\label{sec:lem_b_splines_approx_proof}

We prove claims $(i)$ and $(ii)$ separately.

\noindent
\textbf{Proving claim $(i)$.}\quad
First, assume that $f \in C^2([-h, 1 - t_0 + h])$ and set $\alpha_{-1} = f(-h)$, $\alpha_{K + 1} = f(1 - t_0 + h)$.
Let $\bar{f}$ be the quasi-interpolant with these coefficients to avoid confusions.
We will reduce the original problem to this setting later in the proof.
Using the Taylor expansion around $x$, we obtain
\begin{align*}
    f(y) = f(x) + f'(x)(y - x) + \integral{x}^y (y - t) f''(t) \, \dd t .
\end{align*}
Hence, substituting this into the definition of $\bar f$, we deduce that, for all $x \in [0,1-t_0]$,
\begin{align*}
    \bar f(x)
    = \sum_{k = -1}^{K + 1}\left(f(x) + f'(x)(kh - x) + \integral{x}^{kh}(kh - t) f''(t) \, \dd t \right) B_3(h^{-1}x - (k - 2)) .
\end{align*}
Using the fact that $\{ B_{3}(h^{-1} x - (k - 2))(x) : -1 \leq k \leq K + 1\}$ sum to one and
\begin{align*}
    \sum_{k = -1}^{K + 1} (kh - x) B_3(h^{-1}x - (k - 2)) = 0
\end{align*}
for $0 \leq x \leq 1 - t_0$ (see Lemma \ref{lem:b_splines_basic_prop}), we conclude that
\begin{align*}
    |\bar f(x) - f(x)|
    &\leq \sum_{k = -1}^{K + 1} \left| \integral{x}^{kh}(kh - t) f''(t) \, \dd t \right| B_3(h^{-1} x - (k - 2))
    \\&
    \leq |f|_{W^{2, \infty}([-h, 1 - t_0 + h])}\sum_{k = -1}^{K + 1} \frac{(kh - x)^2}{2} B_3(h^{-1}x - (k - 2)) .
\end{align*}
Notice that $B_3(h^{-1}x - (k - 2))$ is non-zero if and only if $(k - 2)h \leq x \leq (k + 2)h$ and, hence,
\begin{align}
\label{eq:b_spline_bar_f_approx}
    \notag
    |\bar f(x) - f(x)|
    &\leq |f|_{W^{2, \infty}([-h, 1 - t_0 + h])} \sum_{k = -1}^{K + 1} \frac{(2h)^2}{2} B_3(h^{-1}x - (k - 2))
    \\&
    = 2 h^2 |f|_{W^{2, \infty}([-h, 1 - t_0 + h])}
\end{align}
holds for all $x \in [0,1-t_0]$.
We now consider the case where the function $f$ is defined on $[0, 1 - t_0]$.
Extend $f$ to the interval $[-h, 1 - t_0]$ by parabolic extrapolation:
\[
    f(x) = f(0) + f'(0)x + f''(0)x^2 / 2,
    \quad x \in [-h, 0],
\]
and, similarly,
\[
    f(x) = f(1 - t_0) + f'(1 - t_0)(x - 1 + t_0) + f''(1 - t_0)(x - 1 + t_0)^2 / 2,
    \quad x \in [1 - t_0, 1 - t_0 + h] .
\]
Hence, the resulting extension belongs to $C^2([-h, 1 - t_0 + h])$.
Taking into account that $0 \leq B_3(x) \leq 1$ for all $x \in \R$, we deduce that
\begin{align*}
    |\tilde f(x) - \bar f(x)|
    \leq |f(-h) - (2f(0) - f(\tau_1))| + |f(1 - t_0 + h) - (2f(\tau_K) - f(\tau_{K - 1}))| .
\end{align*}
We now elaborate on the first term on the right-hand side of the above bound.
From the Taylor theorem it follows that
\begin{align}
\label{eq:b_spline_lin_interp_err}
    \notag
    &|f(-h) - 2f(0) + f(\tau_1)|
    \\&
    \notag
    = \left| f(0) + f'(0)(-h) + \integral{0}^{-h} (-h - t)f''(t) \, \dd t - 2 f(0) + f(0) + f'(0) h + \integral{0}^h (h - t) f''(t) \, \dd t \right|
    \\&
    \leq |f|_{W^{2, \infty}([-h, 1 - t_0 + h])} h^2 .
\end{align}
The same bound also holds for the right boundary term.
Therefore, this bound together with \eqref{eq:b_spline_bar_f_approx} implies that
\begin{align*}
    |\tilde f(x) - f(x)|
    \leq |\tilde f(x) - \bar f(x)| + |\bar f(x) - f(x)|
    \leq 4 |f|_{W^{2, \infty}([-h, 1 - t_0 + h])} h^2,
    \quad 0 \leq x \leq 1 - t_0 .
\end{align*}
Finally, observe that the extension satisfies
\begin{align*}
    |f|_{W^{2, \infty}([-h, 1 - t_0 + h])}
    = |f|_{W^{2, \infty}([0, 1 - t_0])},
\end{align*}
and the claim $(i)$ follows.

\noindent
\textbf{Proving claim $(ii)$.}\quad
As in the proof of claim $(i)$, we first assume that the function $f$ is defined on $[-h, 1 - t_0 + h]$.
By Lemma \ref{lem:b_splines_basic_prop}, for all $x \in [0, 1 - t_0]$, we have that
\begin{align*}
    \bar f'(x)
    &= h^{-1}\sum_{k = -1}^{K + 1} f(kh) B_3'(h^{-1}x - (k - 2))
    \\&
    = h^{-1} \sum_{k = -1}^{K + 1} f(kh)(B_2(h^{-1}x - (k - 2)) - B_2(h^{-1}x - (k - 1)))
    \\&
    = h^{-1} \sum_{k = 0}^{K + 1} [f(kh) - f((k - 1)h)] B_2(h^{-1}x - (k - 2)) ,
\end{align*}
where the last equality follows from the observation that $B_2(h^{-1}x + 3) = 0$ for all $x \geq 0$ and, similarly, $B_2(h^{-1}x - K) = 0$ for all $x \leq Kh$.
The Taylor theorem implies that
\begin{align*}
    &\bar f'(x) - f'(x)
    \\&
    = \sum_{k = 0}^{K + 1} \left[f'((k - 1)h) - f'(x) + h^{-1}\integral{(k - 1)h}^{kh} f''(t)(kh - t) \dd t \right] B_2(h^{-1}x - (k - 2)) .
\end{align*}
This yields the error bound
\begin{align*}
    |\bar f'(x) - f'(x)| \leq \sum_{k = 0}^{K + 1} \left| \integral{x}^{(k - 1)h} f''(t) \, \dd t \right| B_2(h^{-1}x - (k - 2))
    + \frac{h}{2}|f|_{W^{2, \infty}([-h, 1 - t_0 + h])} .
\end{align*}
Since $B_2(h^{-1}x - (k - 2))$ is non-zero only if $(k - 2)h \leq x \leq (k + 1)h$, then
\begin{align}
\label{eq:b_splines_bar_f_f}
    \notag
    |\bar f'(x) - f'(x)|
    &\leq \sum_{k = 0}^{K + 1} (3h) |f|_{W^{2, \infty}([-h, 1 - t_0 + h])} B_2(h^{-1}x - (k - 2)) + \frac{h}{2}|f|_{W^{2, \infty}([-h, 1 - t_0 + h])}
    \\&
    = \frac{7 h}{2} |f|_{W^{2, \infty}([-h, 1 - t_0 + h])} ,
\end{align}
for all $x \in [0, 1 - t_0]$.
Therefore, from \eqref{eq:b_spline_lin_interp_err} and Lemma \ref{lem:b_splines_basic_prop} we deduce that
\begin{align*}
    &|\tilde f'(x) - \bar f'(x)|
    \\&
    \leq h^{-1} \max_{y \in \R}|B_3'(y)| \cdot (|f(-h) - (2f(0) - f(\tau_1))| + |f(1 - t_0 + h) - (2f(\tau_K) - f(\tau_{K - 1}))|)
    \\&
    \leq 4 h |f|_{W^{2, \infty}([-h, 1 - t_0 + h])} .
\end{align*}
Hence, claim $(ii)$ follows immediately from \eqref{eq:b_splines_bar_f_f}, the above bound, and the fact that the extension of $f$ has the same $W^{2, \infty}$-seminorm.

\myendproof

\section{Proof of Theorem \ref{th:estim_main}}
\label{sec:th_estim_main_proof}

We divide the proof into several steps for clarity.

\noindent
\textbf{Verifying the Bernstein condition.}\quad
We begin by establishing that the excess loss function class fulfills Bernstein's condition.
Fix an arbitrary $v \in \cV(L, W, S, B, N, R, L_f)$ (see Definition \ref{def:av_vel_class}).
From \eqref{eq:self_dist_loss_def} we deduce that
\begin{align}
\label{eq:var_sd_fm}
    \Var[\ell_{SD}(v, X_1) - \ell_{SD}(v^*, X_1)]
     = \Var[\ell_{FM}(v, X_1) - \ell_{FM}(v^*, X_1)] .
\end{align}
Therefore, since $R(v)$ is independent of $X_1$, it suffices to verify the Bernstein condition for the flow matching objective.
Note that this is similar to the case of the denoising score matching objective, which was previously done by \cite{yakovlev2025generalization}.
The key technical step is formalized in the following lemma.

\begin{Lem}
\label{lem:bernstein_helper}
    Assume Assumption \ref{asn:flow_matching_coefs} holds, and let $f, \Breve f : \R^D \times [0, 1 - t_0] \to \R^D$ be such that
    \begin{align*}
        \|f\|_{L^\infty(\R^D \times [0, 1 - t_0])} \vee \|\Breve f\|_{L^\infty(\R^D \times [0, 1 - t_0])} \leq M
    \end{align*}
    for some $M > 0$.
    For all $(y, t) \in \R^D \times [0, 1 - t_0]$, consider the corresponding velocity fields
    \begin{align*}
        b(y, t) = c_1(t, \sigma)y + c_2(t, \sigma) f(y, t),
        \quad \Breve b(y, t) = c_1(t, \Breve\sigma)y + c_2(t, \Breve \sigma) \Breve f(y, t),
    \end{align*}
    where $\sigma, \Breve \sigma \in [\smin, 1]$, and the functions $c_1$ and $c_2$ are defined as in Proposition \ref{prop:true_vel_repr}.
    Then, for every $y \in \R^D$, it follows that
    \begin{align*}
        &(\ell_{FM}(b, y) - \ell_{FM}(\Breve b, y))^2
        \\&
        \lesssim \left((\|y\|^2 + D + M^2) \smin^{-4} + D |\dot \alpha_t|^2_{W^{0, \infty}([0, 1 - t_0])} \right)  \integral{0}^{1 - t_0} \E_{X_t | X_1 = y} \|b(X_t, t) - \Breve b(X_t, t)\|^2 \, \dd t .
    \end{align*}
    Furthermore, we have that
    \begin{align*}
        \integral{0}^{1 - t_0} \E_{X_t | X_1 = y} \|b(X_t, t) - \Breve b(X_t, t)\|^2 \, \dd t
        \lesssim (\|y\|^2 + D + M^2) \smin^{-4} .
    \end{align*}
    and
    \begin{align*}
        \integral{0}^{1 - t_0} \E_{X_t | X_1 = y} \|\Breve b(X_t, t) - \dot X_t\|^2 \, \dd t
        \lesssim (\|y\|^2 + D + M^2) \smin^{-4} + D |\dot \alpha_t|^2_{W^{0, \infty}([0, 1 - t_0])} .
    \end{align*}
    Here, the hidden constants in the inequalities are absolute.

\end{Lem}

The proof of Lemma \ref{lem:bernstein_helper} can be found in Appendix \ref{sec:lem_bernstein_helper_proof}.
With the aid of Lemma \ref{lem:bernstein_helper}, we show that the Bernstein condition holds for the excess loss class
\begin{align}
\label{eq:excess_loss_cls}
    \cL = \{ \ell_{SD}(v, \cdot) - \ell_{SD}(v^*, \cdot) \, : \, v \in \cV(L, W, S, B, N, R, L_f) \}.
\end{align}
Notice that, according to Definition \ref{def:av_vel_class}, we can set $M \asymp \sqrt{D}$ in Lemma \ref{lem:bernstein_helper}.
Now fix some $\varkappa \geq 2$, which will be determined later in the proof.
Then, in view of \eqref{eq:var_sd_fm}, Lemma \ref{lem:bernstein_helper}, and the Hölder inequality, we obtain that
\begin{align*}
    &\Var[\ell_{SD}(v, X_1) - \ell_{SD}(v^*, X_1)]
    \lesssim
    \left[ \integral{0}^{1 - t_0} \|v_{t, t} - b^*(\cdot, t)\|^2_{L^2(\sfp_t^*)} \, \dd t \right]^{1 - 1 / \varkappa}
    \\&\quad
    \cdot
    \E^{1 / \varkappa}\left[ \left((\|X_1\|^2 + D)\smin^{-4} + D \lambda^2(t_0)\right)^\varkappa \integral{0}^{1 - t_0} \E_{X_t | X_1} \|v_{t, t}(X_t) - b^*(X_t, t)\|^2 \, \dd t \right] .
\end{align*}
Here we introduced
\begin{align}
\label{eq:lambda_t0_def}
    \lambda(t_0) = |\alpha_t|_{W^{1, \infty}([0, 1 - t_0])} \vee e .
\end{align}
By the second part of Lemma \ref{lem:bernstein_helper}, we have that
\begin{align*}
    &\Var[\ell_{SD}(v, X_1) - \ell_{SD}(v^*, X_1)]
    \\&
    \lesssim \left[ \integral{0}^{1 - t_0} \|v_{t, t} - b^*(\cdot, t)\|^2_{L^2(\sfp_t^*)} \, \dd t \right]^{1 - 1 / \varkappa}
    \left[ \E\left((\|X_1\|^2 + D)\smin^{-4} + D\lambda^2(t_0) \right)^{\varkappa + 1} \right]^{1 / \varkappa} .
\end{align*}
By Assumption \ref{asn:relax_man}, it follows that
\begin{align}
\label{eq:X1_square_psi_1_bound}
    \left\| \|X_1\|^2 \right\|_{\psi_1}
    \lesssim \left\| \|g^*(U)\|^2 \right\|_{\psi_1} + \sdata^2 \left\| \|Z\|^2 \right\|_{\psi_1}
    \lesssim 1 + \sdata^2 D .
\end{align}
Therefore, applying \citep[Proposition 2.7.1]{vershynin2018high}, we obtain that
\begin{align*}
    &\left[ \E\left((\|X_1\|^2 + D)\smin^{-4} + D\lambda^2(t_0) \right)^{\varkappa + 1} \right]^{1 / \varkappa}
    \\&
    \lesssim (\varkappa + 1)^{1 + 1 / \varkappa} \left\{ \left( \left\| \|X_1\|^2 \right\|_{\psi_1} + D \right)\smin^{-4} +D \lambda^2(t_0) \right\}^{1 + 1 / \varkappa}
    \\&
    \lesssim \varkappa \left( D(\smin^{-4} + \lambda^2(t_0)) \right)^{1 + 1 / \varkappa} .
\end{align*}
Thus, we deduce that, for all $v \in \cV(L, W, S, B, N, R, L_f)$,
\begin{align*}
    &\Var[\ell_{SD}(v, X_1) - \ell_{SD}(v^*, X_1)]
    \\&
    \lesssim \varkappa \left( D(\smin^{-4} + \lambda^2(t_0)) \right)^{1 + 1 / \varkappa}  \left[ \integral{0}^{1 - t_0} \|v_{t, t} - b^*(\cdot, t)\|^2_{L^2(\sfp_t^*)} \, \dd t \right]^{1 - 1 / \varkappa} .
\end{align*}
Applying Proposition \ref{prop:fm_loss_dif} and using the fact that $R(v)$ is non-negative, we conclude that
\begin{align}
\label{eq:bernst_cond_sd}
    \notag
    &\Var[\ell_{SD}(v, X_1) - \ell_{SD}(v^*, X_1)]
    \\&
    \lesssim \varkappa \left( D(\smin^{-4} + \lambda^2(t_0)) \right)^{1 + 1 / \varkappa} \left\{ \E[\ell_{SD}(v, X_1) - \ell_{SD}(v^*, X_1)] \right\}^{1 - 1 / \varkappa} .
\end{align}
Consequently, the Bernstein condition is satisfied.

\noindent
\textbf{Large deviation bound.}\quad
Fix an arbitrary $v \in \cV(L, W, S, B, N, R, L_f)$ and denote
\begin{align*}
    \widehat{\E}[\ell_{SD}(v, X_1) - \ell_{SD}(v^*, X_1)]
    = \frac{1}{n}\sum_{i = 1}^n \left(\ell_{SD}(v, Y_i) - \ell_{SD}(v^*, Y_i)\right).
\end{align*}
Recall that $Y_1, \dots, Y_n$ are independent copies of $X_1$.
We now provide a high-probability upper bound on
\begin{align*}
    \left| \widehat{\E}[\ell_{SD}(v, X_1) - \ell_{SD}(v^*, X_1)] - \E[\ell_{SD}(v, X_1) - \ell_{SD}(v^*, X_1)] \right| .
\end{align*}
To achieve this, we invoke the Bernstein inequality for unbounded random variables \cite[Proposition 5.2]{Lecue12oracle}.
This inequality states that, for any $\delta \in (0, 1)$, with probability at least $(1 - \delta / 2)$,
\begin{align*}
    &\left| \widehat{\E}[\ell_{SD}(v, X_1) - \ell_{SD}(v^*, X_1)] - \E[\ell_{SD}(v, X_1) - \ell_{SD}(v^*, X_1)] \right|
    \\&
    \lesssim \sqrt{\frac{\Var[\ell_{SD}(v, X_1) - \ell_{SD}(v^*, X_1)] \log(4 / \delta)}{n}}
    \\&
    + \frac{\| \ell_{SD}(v, X_1) - \ell_{SD}(v^*, X_1) \|_{\psi_1} \log n \log(4 / \delta)}{n} .
\end{align*}
We now derive an upper bound on the $\psi_1$-norm of the excess loss $\ell_{SD}(v, X_1) - \ell_{SD}(v^*, X_1)$.
This bound is provided by the following lemma.
\begin{Lem}
\label{lem:psi_1_diameter}
    Assume Assumptions \ref{asn:relax_man} and \ref{asn:flow_matching_coefs} hold.
    Assume also $h \lesssim \smin^4$.
    Then, it follows that
    \begin{align*}
        \sup_{v \in \cV(L, W, S, B, N, R, L_f)} \left\|\ell_{SD}(v, X_1) - \ell_{SD}(v^*, X_1) \right\|_{\psi_1}
        \lesssim D \left( h^{-2} \exp\{\cO(\smin^{-2})\} + \lambda^2(t_0) \right) .
    \end{align*}
\end{Lem}

The proof of Lemma \ref{lem:psi_1_diameter} is deferred to Appendix \ref{sec:lem_psi_1_diameter_proof}.
Therefore, setting
\begin{align}
\label{eq:estim_kappa_def}
    \varkappa \asymp \log n + \log \lambda(t_0)
\end{align}
and applying this lemma together with \eqref{eq:bernst_cond_sd} yields, for each fixed candidate $v\in\cV$,
\begin{align}
\label{eq:bernst_fixed}
    \notag
    &\left| \widehat{\E}[\ell_{SD}(v, X_1) - \ell_{SD}(v^*, X_1)] - \E[\ell_{SD}(v, X_1) - \ell_{SD}(v^*, X_1)] \right|
    \\&
    \lesssim \sqrt{ \frac{C_b \{ \E[\ell_{SD}(v, X_1) - \ell_{SD}(v^*, X_1)] \}^{1 - 1 / \varkappa} \log(4 / \delta) }{n} }
    + \frac{C_b \log(4 / \delta)}{n} ,
\end{align}
with probability at least $(1 - \delta / 2)$.
Here, we introduced
\begin{align}
\label{eq:c_bern_def}
    C_b = D^2 h^{-2} \exp\{\cO(\smin^{-2})\}\lambda^2(t_0) \log(\lambda(t_0)) \log n .
\end{align}

\noindent
\textbf{Covering number evaluation.}\quad
We now derive a uniform Bernstein-type bound for the excess loss class, as defined in \eqref{eq:excess_loss_cls}.
To achieve this, we need to evaluate the covering number of this class, which is established by the following lemma.

\begin{Lem}
\label{lem:cover_exc_loss}
    Under Assumptions \ref{asn:relax_man} and \ref{asn:flow_matching_coefs}, fix $\eps \in (0, 1)$.
    Let the parameters $R$ and $L_f$ in the class $\cV(L, W, S, B, N, R, L_f)$ (see Definition \ref{def:av_vel_class}) be given as in Theorem \ref{th:avg_vel_approx}, and let the parameter $h \lesssim \smin^4$.
    Then, for every $\delta \in (0, 1)$, there exists a set $\cV_\eps \subseteq \cV(L, W, S, B, N, R, L_f)$ such that
    \begin{align*}
        &\sup_{v \in \cV} \inf_{v' \in \cV_\eps} \left\{ |\E[\ell_{SD}(v, X_1) - \ell_{SD}(v', X_1)]| + |\widehat{\E}[\ell_{SD}(v, X_1) - \ell_{SD}(v', X_1)]| \right\}
        \\&
        \lesssim h^{-2} \exp\{\cO(\smin^{-4})\} (D \log(1 / \eps))^{\cO(1)} \log(n / \delta)\lambda(t_0) \eps^{2\beta}
    \end{align*}
    holds with probability at least $(1 - \delta)$.
    Moreover, it holds that
    \begin{align*}
        \log |\cV_\eps| \lesssim SL \smin^{-2} h^{-1} \log(\eps^{-1} L (B \vee 1)(\|W\|_\infty + 1)D) + SLh^{-1}\log(e \vee \log(n / \delta)) .
    \end{align*}

\end{Lem}

The proof of Lemma \ref{lem:cover_exc_loss} is given in Appendix \ref{sec:lem_cover_exc_loss_proof}.
Define
\begin{align}
\label{eq:estim_tau_def}
    \tau = h^{-2} \exp\{\cO(\smin^{-4})\} (D \log(1 / \eps))^{\cO(1)} \log(n / \delta)\lambda(t_0) \eps^{2\beta} .
\end{align}
Assume $h \lesssim \smin^4$ throughout the remainder of the proof.
Then, with the aid of this lemma and the union bound, we obtain from \eqref{eq:bernst_fixed} that, with probability at least $(1 - \delta)$, the following inequality holds for all $v \in \cV(L, W, S, B, N, R, L_f)$ simultaneously:
\begin{align*}
    &\left| \widehat{\E}[\ell_{SD}(v, X_1) - \ell_{SD}(v^*, X_1)] - \E[\ell_{SD}(v, X_1) - \ell_{SD}(v^*, X_1)] \right|
    \\&
    \lesssim \tau + \sqrt{ \frac{ C_b\{ \tau + \E[\ell_{SD}(v, X_1) - \ell_{SD}(v^*, X_1)] \}^{1 - 1 / \varkappa} \log(4 |\cV_\eps| / \delta)}{n} }
    + \frac{C_b \log(4 |\cV_\eps| / \delta)}{n}
    \\&
    \lesssim \tau + \sqrt{\frac{C_b \{ \E[\ell_{SD}(v, X_1) - \ell_{SD}(v^*, X_1)] \}^{1 - 1 / \varkappa} \log(4 |\cV_\eps| / \delta)}{n}}
    + \frac{C_b \log(4 |\cV_\eps| / \delta)}{n}
    \\&\quad
    + \left(\frac{C_b \log(4 |\cV_\eps| / \delta)}{n}\right)^{\varkappa / (\varkappa + 1)},
\end{align*}
where the last inequality uses the Young inequality.

\noindent
\textbf{Deriving the final bound.}\quad
In this step, we aim to establish a high-probability upper bound on the expected excess loss
\begin{align*}
    \E_{X_1}[\ell_{SD}(\widehat{v}, X_1) - \ell_{SD}(v^*, X_1)],
\end{align*}
where $X_1$ is independent from the sample $Y_1, \dots, Y_n$.
Since $\widehat{v}$ minimizes the empirical risk $\widehat{\E}_{X_1} \ell_{SD}(v, X_1)$ with respect to $v \in \cV(L, W, S, B, N, R, L_f)$, we deduce that
\begin{align*}
    &\E_{X_1}[\ell_{SD}(\widehat{v}, X_1) - \ell_{SD}(v^*, X_1)]
    \leq \E_{X_1}[\ell_{SD}(\tilde{v}, X_1) - \ell_{SD}(v^*, X_1)]
    \\&
    + \left| \widehat{\E}_{X_1}[\ell_{SD}(\widehat v, X_1) - \ell_{SD}(v^*, X_1)] - \E_{X_1}[\ell_{SD}(\widehat v, X_1) - \ell_{SD}(v^*, X_1)] \right|
    \\&
    + \left| \widehat{\E}_{X_1}[\ell_{SD}(\tilde v, X_1) - \ell_{SD}(v^*, X_1)] - \E_{X_1}[\ell_{SD}(\tilde v, X_1) - \ell_{SD}(v^*, X_1)] \right|,
\end{align*}
where $\tilde v \in \cV(L, W, S, B, N, R, L_f)$ is arbitrary.
By Young's inequality, with probability at least $(1 - \delta)$,
\begin{align*}
    &\E_{X_1}[\ell_{SD}(\widehat{v}, X_1) - \ell_{SD}(v^*, X_1)]
    \lesssim \E_{X_1}[\ell_{SD}(\tilde v, X_1) - \ell_{SD}(v^*, X_1)]
    + \tau
    + \frac{C_b \log(4 |\cV_\eps| / \delta)}{n}
    \\&
    + \sqrt{\frac{C_b \{ \E_{X_1}[\ell_{SD}(\widehat v, X_1) - \ell_{SD}(v^*, X_1)] \}^{1 - 1 / \varkappa} \log(4 |\cV_\eps| / \delta)}{n}}
    + \left(\frac{C_b \log(4 |\cV_\eps| / \delta)}{n}\right)^{\varkappa / (\varkappa + 1)} .
\end{align*}
Using the fact that $x \leq a \sqrt{x^{1 - 1 / \varkappa}} + b$ implies $x \leq 4 a^{2\varkappa / (1 + \varkappa)} + 2b$ for all non-negative $a$, $b$, and $x$, we deduce that
\begin{align*}
    &\E_{X_1}[\ell_{SD}(\widehat{v}, X_1) - \ell_{SD}(v^*, X_1)]
    \lesssim \E_{X_1}[\ell_{SD}(\tilde v, X_1) - \ell_{SD}(v^*, X_1)]
    + \tau
    \\&
    + \left(\frac{C_b \log(4 |\cV_\eps| / \delta)}{n}\right)^{\varkappa / (\varkappa + 1)}
    + \frac{C_b \log(4 |\cV_\eps| / \delta)}{n} ,
\end{align*}
with probability at least $(1 - \delta)$.
We now utilize the approximation result provided by Theorem \ref{th:avg_vel_approx}.
Fix an arbitrary $\eps \in (0, 1)$ satisfying the conditions of Theorem \ref{th:avg_vel_approx}, and let $\tilde v$ be the corresponding approximation.
Then, taking into account that $h^{-2} \asymp \eps^{-\beta} \smin^{-4}$, we obtain from \eqref{eq:estim_tau_def} that
\begin{align*}
    \tau
    &= h^{-2} \exp\{\cO(\smin^{-4})\} (D \log(1 / \eps))^{\cO(1)} \log(n / \delta)\lambda(t_0) \eps^{2\beta}
    \\
    &\lesssim \exp\{\cO(\smin^{-4})\} (D \log(1 / \eps))^{\cO(1)} \log(n / \delta)\lambda(t_0) \eps^\beta .
\end{align*}
Consequently, with probability at least $(1 - \delta)$, it holds that
\begin{align*}
    &\E_{X_1}[\ell_{SD}(\widehat{v}, X_1) - \ell_{SD}(v^*, X_1)]
    \lesssim
    (D \log(1 / \eps))^{\cO(1)} \exp\{\cO(\smin^{-4})\} \log(n / \delta) \lambda(t_0) \eps^\beta
    \\&
    + \left(\frac{C_b \log(4 |\cV_\eps| / \delta)}{n}\right)^{\varkappa / (\varkappa + 1)}
    + \frac{C_b \log(4 |\cV_\eps| / \delta)}{n} .
\end{align*}
Now it remains to bound $\log(|\cV_\eps|)$.
We must ensure that the approximation $\tilde v$ belongs to $\cV(L, W, S, B, N, R, L_f)$.
Therefore, it suffices to take the architecture and the parameters $L_f$ and $R$ as specified in Theorem \ref{th:avg_vel_approx}.
Thus, by Lemma \ref{lem:cover_exc_loss}, we have that
\begin{align*}
    \log |\cV_\eps|
    &\lesssim \eps^{-d}h^{-1}\exp\{\cO(\smin^{-2})\} (D \log(1 / \eps))^{\cO(1)} \log(n / \delta)
    \\&
    \lesssim \eps^{-d - \beta / 2} \exp\{\cO(\smin^{-2})\} (D \log(1 / \eps))^{\cO(1)} \log(n / \delta).
\end{align*}
Since $\varkappa \gtrsim \log n$ due to \eqref{eq:estim_kappa_def}, we deduce that, with probability at least $(1 - \delta)$,
\begin{align*}
    &\E_{X_1}[\ell_{SD}(\widehat{v}, X_1) - \ell_{SD}(v^*, X_1)]
    \\&
    \lesssim
    (D \log(1 / \eps))^{\cO(1)} \exp\{\cO(\smin^{-4})\}
    \log(n / \delta) \left\{ \lambda(t_0) \eps^\beta
    + \frac{(C_b \vee 1)\eps^{-d - \beta / 2}}{n} \right\} .
\end{align*}
From \eqref{eq:estim_kappa_def} and \eqref{eq:c_bern_def} we conclude that
\begin{align*}
    &\E_{X_1}[\ell_{SD}(\widehat{v}, X_1) - \ell_{SD}(v^*, X_1)]
    \\&
    \lesssim (D \log(1 / \eps))^{\cO(1)} \exp\{\cO(\smin^{-4})\} \lambda(t_0) \left\{ \eps^\beta + \frac{\lambda(t_0) \eps^{-d - 3\beta / 2}}{n} \right\} \log(\lambda(t_0)) \log(4 / \delta) \log^2 n
    ,
\end{align*}
with probability at least $(1 - \delta)$.
Thus, setting $\eps = (\lambda(t_0) / n)^{1 / (d + 5\beta / 2)}$ yields that, with probability at least $(1 - \delta)$,
\begin{align*}
    &\E_{X_1}[\ell_{SD}(\widehat{v}, X_1) - \ell_{SD}(v^*, X_1)]
    \\&
    \lesssim (D \log n)^{\cO(1)} \exp\{\cO(\smin^{-4})\}
    \lambda(t_0)
    \left(\frac{\lambda(t_0)}{n} \right)^{\frac{2\beta}{5\beta + 2d}} \log(\lambda(t_0)) \log(4 / \delta).
\end{align*}
It now remains to ensure that this choice of $\eps$ satisfies the constraints of Theorem \ref{th:avg_vel_approx}.
This requires the sample size $n$ to be sufficiently large, specifically
\begin{align*}
    D^{3/2} \exp\{\cO(\smin^{-2})\}\log(nD\smin^{-2})
    \left(\frac{\lambda(t_0)}{n}\right)^{\frac{2(\beta\wedge1)}{5\beta+2d}}\lesssim1,
    \quad
    \left(\frac{\lambda(t_0)}{n}\right)^{\frac{2\beta}{5\beta+2d}} \lesssim \smin^4,
    \quad n \gtrsim \lambda(t_0).
\end{align*}
Furthermore, Theorem \ref{th:avg_vel_approx} suggests that the architecture of the class $\cV(L, W, S, B, N, R, L_f)$ is the following:
\begin{align*}
    &L \lesssim \smin^{-2} \log(D \smin^{-2} \log n),
    \quad \|W\|_\infty \vee S \lesssim \left( \frac{n}{\lambda(t_0)} \right)^{\frac{2d}{5\beta + 2d}} \exp\{\cO(\smin^{-2})\} (D \log n)^{\cO(1)},
    \\
    &\log B \lesssim \exp\{\cO(\smin^{-2})\} (D \log n)^{\cO(1)} .
\end{align*}
Finally, the theorem also implies that
\begin{align*}
    L_f \asymp (\smin^{-2}D \log n)^{\cO(1)},
    \quad
    N \asymp \left( \frac{n}{\lambda(t_0)} \right)^{\frac{\beta}{5\beta + 2d}} \smin^{-2},
    \quad
    R \asymp \exp\{\cO(\smin^{-2})\}D \log(nD\smin^{-2})
\end{align*}
in Definition \ref{def:av_vel_class}.
The explicit smallness condition on $\eps^\beta$ and the sufficiently large constant in $N$ ensure $h\leq c\smin^4$, as required.
The proof is now complete.

\myendproof

\subsection{Proof of Lemma \ref{lem:bernstein_helper}}
\label{sec:lem_bernstein_helper_proof}

We begin by decomposing the loss difference
\begin{align*}
    \ell_{FM}(b, y) - \ell_{FM}(\Breve b, y)
    &= \integral{0}^{1 - t_0}\E_{X_t | X_1 = y}\|b(X_t, t) - \Breve b(X_t, t)\|^2 \, \dd t
    \\&\quad
    + 2 \integral{0}^{1 - t_0} \E_{X_t\mid X_1=y}\langle b(X_t, t) - \Breve b(X_t, t), \Breve b(X_t, t) - \dot{X}_t \rangle \, \dd t .
\end{align*}
Applying the Cauchy-Schwarz inequality yields
\begin{align}
\label{eq:fm_loss_diff_sq_aux}
    &(\ell_{FM}(b, y) - \ell_{FM}(\Breve b, y))^2
    \leq 2 \left(  \integral{0}^{1 - t_0}\E_{X_t | X_1 = y}\|b(X_t, t) - \Breve b(X_t, t)\|^2 \, \dd t \right)^2
    \\&\quad
    \notag
    + 8 \integral{0}^{1 - t_0} \E_{X_t | X_1 = y} \|b(X_t, t) - \Breve b(X_t, t)\|^2 \, \dd t \integral{0}^{1 - t_0} \E_{X_t | X_1 = y} \|\Breve b(X_t, t) - \dot{X}_t\|^2 \, \dd t .
\end{align}
Following this, by Proposition \ref{prop:true_vel_repr} and Assumption \ref{asn:flow_matching_coefs}, we have that
\begin{align*}
    \max_{\smin \leq \sigma \leq 1} \left( |c_1(\cdot, \sigma)|_{W^{0, \infty}([0, 1 - t_0])} \vee |c_2(\cdot, \sigma)|_{W^{0, \infty}([0, 1 - t_0])}\right) \lesssim \smin^{-2} .
\end{align*}
Hence, it follows that
\begin{align*}
    &\|b(X_t, t) - \Breve b(X_t, t)\|^2
    \\&
    \leq 3 \|X_t\|^2(c_1(t, \sigma) - c_1(t, \Breve \sigma))^2 + 3c_2^2(t, \sigma)\|f(X_t, t)\|^2 + 3c_2^2(t, \Breve \sigma) \|\Breve f(X_t, t)\|^2
    \\&
    \lesssim (\|X_t\|^2 + M^2)\smin^{-4} .
\end{align*}
Using the moment bound
\begin{align*}
    \E_{X_t | X_1 = y}\|X_t\|^2
    \leq 2\|\beta_t y\|^2 + 2 \E_{X_t | X_1 = y}\|X_t - \beta_t y\|^2
    = 2\beta_t^2 \|y\|^2 + 2D \alpha_t^2 ,
\end{align*}
we conclude that
\begin{align}
\label{eq:vel_field_diff_cond_aux}
    \integral{0}^{1 - t_0} \E_{X_t | X_1 = y} \|b(X_t, t) - \Breve b(X_t, t)\|^2 \, \dd t
    \lesssim (\|y\|^2 + D + M^2) \smin^{-4} ,
\end{align}
thereby establishing the second claim of the lemma.
We next observe that
\begin{align*}
    \E_{X_t | X_1 = y} \| \dot X_t \|^2
    \leq 2 \dot\beta_t^2 \|y\|^2 + 2 \dot\alpha_t^2 D .
\end{align*}
Therefore, it follows that
\begin{align*}
    &\integral{0}^{1 - t_0} \E_{X_t | X_1 = y} \|\Breve b(X_t, t) - \dot X_t\|^2 \, \dd t
    \\&
    \leq 3 \integral{0}^{1 - t_0} \left( c_1^2(t, \Breve \sigma)\E_{X_t | X_1 = y} \|X_t\|^2 + c_2^2(t, \Breve \sigma) M^2 + \E_{X_t | X_1 = y}\|\dot X_t\|^2 \right)\, \dd t
    \\&
    \lesssim (\|y\|^2 + D + M^2) \smin^{-4} + D |\dot \alpha_t|^2_{W^{0, \infty}([0, 1 - t_0])} .
\end{align*}
Consequently, the third claim is established.
Combining this with \eqref{eq:fm_loss_diff_sq_aux} and \eqref{eq:vel_field_diff_cond_aux}, we conclude that
\begin{align*}
    &(\ell_{FM}(b, y) - \ell_{FM}(\Breve b, y))^2
    \\&
    \lesssim \left((\|y\|^2 + D + M^2) \smin^{-4} + D |\dot \alpha_t|^2_{W^{0, \infty}([0, 1 - t_0])} \right)  \integral{0}^{1 - t_0} \E_{X_t | X_1 = y} \|b(X_t, t) - \Breve b(X_t, t)\|^2 \, \dd t .
\end{align*}
The proof is finished.

\myendproof

\subsection{Proof of Lemma \ref{lem:psi_1_diameter}}
\label{sec:lem_psi_1_diameter_proof}

For brevity, we suppress the full notation $\cV(L, W, S, B, N, R, L_f)$ and simply write $\cV$.
It holds that
\begin{align*}
    \sup_{v \in \cV}\left\|\ell_{SD}(v, X_1) - \ell_{SD}(v^*, X_1) \right\|_{\psi_1}
    &\lesssim \sup_{v \in \cV} \left\| \ell_{FM}(v, X_1) - \ell_{FM}(v^*, X_1) \right\|_{\psi_1}
    + \sup_{v \in \cV} R(v) .
\end{align*}
We proceed by analyzing each term individually.
Applying the third claim of Lemma \ref{lem:bernstein_helper} with $M \asymp \sqrt{D}$, we obtain that
\begin{align}
\label{eq:psi_1_diam_triang_aux}
    \sup_{v \in \cV} \left\|\ell_{FM}(v, X_1) - \ell_{FM}(v^*, X_1) \right\|_{\psi_1}
    \notag
    &\leq \sup_{v \in \cV} \left\| \ell_{FM}(v, X_1) \right\|_{\psi_1}
    + \sup_{v \in \cV}\left\| \ell_{FM}(v^*, X_1) \right\|_{\psi_1}
    \\&
    \lesssim \left(\left\| \|X_1\|^2 \right\|_{\psi_1} + D \right) \smin^{-4} + D \lambda^2(t_0) .
\end{align}
Using \eqref{eq:X1_square_psi_1_bound}, we deduce that
\begin{align}
\label{eq:l_fm_psi_1_diam}
    \sup_{v \in \cV} \left\| \ell_{FM}(v, X_1) - \ell_{FM}(v^*, X_1) \right\|_{\psi_1}
    \lesssim D(\smin^{-4} + \lambda^2(t_0)) .
\end{align}
We now analyze the consistency penalty term.
Let $\Psi_{0, t}(y) = y + t \, v_{0, t}(y)$ be as in Definition \ref{def:av_vel_class}.
Then it follows that
\begin{align*}
    R(v)
    &\lesssim \integral{0}^{1 - t_0} \E \|\partial_t \Psi_{0, t}(X_0)\|^2 \, \dd t + \integral{0}^{1 - t_0} \E \|v_{t, t}(\Psi_{0, t}(X_0))\|^2 \, \dd t
    \\&
    \lesssim \integral{0}^{1 - t_0} h^{-2} \E \max_{0 \leq k \leq N} \|Z_k\|^2 \, \dd t
    + \integral{0}^{1 - t_0} (c_1(t, \sigma)^2 \E \|\Psi_{0, t}(X_0)\|^2 + c_2(t, \sigma)^2 D) \, \dd t .
\end{align*}
In view of Proposition \ref{prop:true_vel_repr}, we have that
\begin{align*}
    |c_j(\cdot, \sigma)|_{W^{0, \infty}([0, 1 - t_0])} \lesssim \smin^{-2},
    \quad j \in \{1, 2\} .
\end{align*}
This implies that
\begin{align*}
    R(v) \lesssim \integral{0}^{1 - t_0} \left( (h^{-2} + \smin^{-4}) \E \max_{0 \leq k \leq N} \|Z_k\|^2 + \smin^{-4}D \right) \, \dd t .
\end{align*}
Applying the second claim of Lemma \ref{lem:heun_error_bound} with $M \asymp \smin^{-2}$ and $M_f = \smin^{-2}\sqrt{D}$, we conclude that
\begin{align*}
    \sup_{v \in \cV} R(v) \lesssim D h^{-2} \exp\{\cO(\smin^{-2})\} .
\end{align*}
Combining this with \eqref{eq:psi_1_diam_triang_aux} and \eqref{eq:l_fm_psi_1_diam}, we obtain that
\begin{align*}
    \sup_{v \in \cV} \left\| \ell_{SD}(v, X_1) - \ell_{SD}(v^*, X_1) \right\|_{\psi_1}
    \lesssim D \left( h^{-2} \exp\{\cO(\smin^{-2})\} + \lambda^2(t_0) \right) .
\end{align*}
This completes the proof.

\myendproof

\subsection{Proof of Lemma \ref{lem:cover_exc_loss}}
\label{sec:lem_cover_exc_loss_proof}

Let $\cB(0, \tilde R) \subseteq \R^D$, where $\tilde R > 0$ will be chosen later in the proof, and consider average velocity field candidates $v^{(1)}, v^{(2)} \in \cV(L, W, S, B, N, R, L_f)$ with corresponding
\begin{align*}
    f^{(j)}(y, t), \quad \{f_k^{(j)}(y, t) \, : \, 0 \leq k \leq N\},
    \quad \smin \leq \sigma^{(j)} \leq 1,
    \quad j \in \{1, 2\}
\end{align*}
as in Definition \ref{def:av_vel_class}.
Then, we have that
\begin{align}
\label{eq:cov_triang_first}
    |\ell_{SD}(v^{(1)}, y) - \ell_{SD}(v^{(2)}, y)| \leq |\ell_{FM}(v^{(1)}, y) - \ell_{FM}(v^{(2)}, y)| + |R(v^{(1)}) - R(v^{(2)})| .
\end{align}
We proceed by examining each term individually.

\noindent
\textbf{The distance between flow matching losses.}\quad
Applying Lemma \ref{lem:bernstein_helper} with $M \asymp \sqrt{D}$, we obtain that for all $y \in \R^D$,
\begin{align}
\label{eq:l_fm_prox_start}
    \notag
    &|\ell_{FM}(v^{(1)}, y) - \ell_{FM}(v^{(2)}, y)|
    \\&
    \lesssim \left( \smin^{-2}(\|y\| + \sqrt{D}) + \sqrt{D} \lambda(t_0) \right)
    \left( \integral{0}^{1 - t_0} \E_{X_t | X_1 = y} \|v_{t, t}^{(1)}(X_t) - v_{t, t}^{(2)}(X_t)\|^2 \, \dd t \right)^{1/2} .
\end{align}
Furthermore, Lemma \ref{lem:bernstein_helper} implies that, for all $y \in \R^D$,
\begin{align}
\label{eq:cov_global_lin_grow}
    |\ell_{FM}(v^{(1)}, y) - \ell_{FM}(v^{(2)}, y)|
    \lesssim \left( \smin^{-2}(\|y\| + \sqrt{D}) + \sqrt{D} \lambda(t_0) \right) \smin^{-2}(\|y\| + \sqrt{D}) .
\end{align}
We now note that
\begin{align*}
    &\integral{0}^{1 - t_0} \E_{X_t | X_1 = y} \|v_{t, t}^{(1)}(X_t) - v_{t, t}^{(2)}(X_t)\|^2 \, \dd t
    \\&
    \lesssim \integral{0}^{1 - t_0} (c_2(t, \sigma^{(1)}) - c_2(t, \sigma^{(2)}))^2 \|f^{(1)}\|^2_{L^\infty(\R^D \times [0, 1 - t_0])} \, \dd t
    \\&
    + \integral{0}^{1 - t_0} (c_1(t, \sigma^{(1)}) - c_1(t, \sigma^{(2)}))^2 \E_{X_t | X_1 = y} \|X_t\|^2 \, \dd t
    \\&
    + \integral{0}^{1 - t_0} c_2^2(t, \sigma^{(2)}) \E_{X_t | X_1 = y} \|f^{(1)}(X_t, t) - f^{(2)}(X_t, t)\|^2 \, \dd t .
\end{align*}
Next, from Proposition \ref{prop:true_vel_repr} we find that, for all $0 \leq t \leq 1 - t_0$,
\begin{align}
\label{eq:cj_proximity}
    |c_j(t, \sigma^{(1)}) - c_j(t, \sigma^{(2)})| \lesssim \smin^{-4} |\sigma^{(1)} - \sigma^{(2)}|,
    \quad |c_j(t, \sigma^{(2)})| \lesssim \smin^{-2},
    \quad j \in \{1, 2\} .
\end{align}
Leveraging the bound $\E_{X_t | X_1 = y}\|X_t\|^2 \lesssim \alpha_t^2 D + \beta_t^2 \|y\|^2$, we obtain
\begin{align*}
    &\integral{0}^{1 - t_0} \E_{X_t | X_1 = y} \|v_{t, t}^{(1)}(X_t) - v_{t, t}^{(2)}(X_t)\|^2 \, \dd t
    \lesssim (D + \|y\|^2)\smin^{-8}(\sigma^{(1)} - \sigma^{(2)})^2
    \\&\quad
    + \smin^{-4} \|f^{(1)} - f^{(2)}\|_{L^\infty(\cB(0, \tilde R) \times [0, 1 - t_0])}^2
    + \smin^{-4} D \max_{0 \leq t \leq 1 - t_0} \p(\|X_t\| \geq \tilde R \, |\, X_1 = y ) .
\end{align*}
The Gaussian concentration bound in Lemma \ref{lem:gaus_conc_norm} implies that choosing the radius $\tilde R$ satisfying
\begin{align}
\label{eq:cover_tilde_R_def}
    \tilde R \asymp \sqrt{\log(1 / \eps)} + \sqrt{D}
\end{align}
ensures that
\begin{align*}
    &\integral{0}^{1 - t_0} \E_{X_t | X_1 = y} \|v_{t, t}^{(1)}(X_t) - v_{t, t}^{(2)}(X_t)\|^2 \, \dd t
    \lesssim (D + \|y\|^2)\smin^{-8}(\sigma^{(1)} - \sigma^{(2)})^2
    \\&\quad
    + \smin^{-4} \|f^{(1)} - f^{(2)}\|_{L^\infty(\cB(0, \tilde R + \|y\|) \times [0, 1 - t_0])}^2
    + \smin^{-4} D \eps^{4\beta} .
\end{align*}
Substituting this into \eqref{eq:l_fm_prox_start}, we conclude that, for all $y \in \R^D$,
\begin{align*}
    &|\ell_{FM}(v^{(1)}, y) - \ell_{FM}(v^{(2)}, y)|
    \lesssim \left( \smin^{-2}(\|y\| + \sqrt{D}) + \sqrt{D} \lambda(t_0) \right)
    \\&
    \cdot \left( \smin^{-4}(\sqrt{D} + \|y\|)|\sigma^{(1)} - \sigma^{(2)}|
    + \smin^{-2} \|f^{(1)} - f^{(2)}\|_{L^\infty(\cB(0, \tilde R + \|y\|) \times [0, 1 - t_0])} + \smin^{-2} \sqrt{D} \eps^{2\beta} \right) .
\end{align*}
By applying Gaussian concentration again, we deduce that setting
\begin{align}
\label{eq:cover_R_data_def}
    R_\text{data} \asymp 1+\sdata \sqrt{\log(n / \delta) + \log(1 / \eps)} + \sdata\sqrt{D}
\end{align}
ensures that, with probability at least $(1 - \delta)$, the following inequality holds:
\begin{align}
\label{eq:cov_fm_ball_aux}
    &\|\ell_{FM}(v^{(1)}, \cdot) - \ell_{FM}(v^{(2)}, \cdot)\|_{L^\infty(\cB(0, R_\text{data}))}
    \lesssim \smin^{-6} \left( \sqrt{\log(n\eps^{-1} / \delta)} + \sqrt{D}\lambda(t_0) \right)
    \\&
    \notag
    \cdot \left( (\sqrt{\log(n \eps^{-1} / \delta)} + \sqrt{D})|\sigma^{(1)} - \sigma^{(2)}| + \|f^{(1)} - f^{(2)}\|_{L^\infty(\cB(0, \tilde R + R_\text{data}) \times [0, 1 - t_0])} + \sqrt{D} \eps^{2\beta} \right) .
\end{align}
By Assumption \ref{asn:relax_man}, we have the bound $\E^{1/2}(\|X_1\| + \sqrt{D})^4 \lesssim D$.
Hence, using \eqref{eq:cov_global_lin_grow}, we obtain that
\begin{align}
\label{eq:cov_fm_exp_prox_aux}
    \notag
    &|\E\ell_{FM}(v^{(1)}, X_1) - \E \ell_{FM}(v^{(2)}, X_1)|
    \\&
    \notag
    \lesssim \|\ell_{FM}(v^{(1)}, \cdot) - \ell_{FM}(v^{(2)}, \cdot)\|_{L^\infty(\cB(0, R_\text{data}))} + \eps^{2\beta} \sqrt{ \E(\ell_{FM}(v^{(1)}, X_1) - \ell_{FM}(v^{(2)}, X_1))^2 }
    \\&
    \lesssim \|\ell_{FM}(v^{(1)}, \cdot) - \ell_{FM}(v^{(2)}, \cdot)\|_{L^\infty(\cB(0, R_\text{data}))}
    + \eps^{2\beta} \smin^{-2} D (\smin^{-2} + \sqrt{D}\lambda(t_0)) .
\end{align}

\noindent
\textbf{The distance between consistency penalties.}\quad
We now analyze the proximity between the penalty functions $R(v^{(1)})$ and $R(v^{(2)})$.
Introduce $\Psi_{0, t}^{(j)}(y) = y + t \, v_{0, t}^{(j)}(y)$ as in Definition \ref{def:av_vel_class}.
By Cauchy-Schwarz, we have that
\begin{align}
\label{eq:r_prox_cs}
    \notag
    &|R(v^{(1)}) - R(v^{(2)})|
    \\&
    \leq  \left( \integral{0}^{1 - t_0} \E \| \partial_t[\Psi_{0,t}^{(1)}(X_0) - \Psi^{(2)}_{0, t}(X_0)]  - v_{t, t}^{(1)}(\Psi_{0, t}^{(1)}(X_0)) + v_{t, t}^{(2)}(\Psi_{0, t}^{(2)}(X_0))\|^2  \, \dd t \right)^{1/2}
    \\&\quad
    \notag
    \cdot
    \left( \integral{0}^{1 - t_0} \E \| \partial_t[\Psi_{0,t}^{(1)}(X_0) + \Psi^{(2)}_{0, t}(X_0)]  - v_{t, t}^{(1)}(\Psi_{0, t}^{(1)}(X_0)) - v_{t, t}^{(2)}(\Psi_{0, t}^{(2)}(X_0))\|^2 \, \dd t \right)^{1 / 2} .
\end{align}
By Lemma \ref{lem:b_splines_basic_prop}, we have that
\begin{align}
\label{eq:flow_der_prox}
    \integral{0}^{1 - t_0} \E \|\partial_t[\Psi_{0, t}^{(1)}(X_0) - \Psi_{0, t}^{(2)}(X_0)]\|^2 \, \dd t
    \lesssim h^{-2} \integral{0}^{1 - t_0} \E \max_{0 \leq k \leq N} \| Z_k^{(1)} - Z_k^{(2)} \|^2 \, \dd t ,
\end{align}
where $Z_k^{(j)}$ is defined in Definition \ref{def:av_vel_class} with $\{f_k^{(j)} \, : \, 0 \leq k \leq N \}$ and $Z_0^{(j)} = X_0$ for $j \in \{1, 2\}$.
In addition, it holds that
\begin{align}
\label{eq:r_second_term_prox}
    \notag
    &\integral{0}^{1 - t_0} \E \| v_{t, t}^{(1)}(\Psi_{0, t}^{(1)}(X_0)) - v_{t, t}^{(2)}(\Psi_{0, t}^{(2)}(X_0)) \|^2 \, \dd t
    \lesssim \integral{0}^{1 - t_0} \E \|v_{t, t}^{(1)}(\Psi_{0, t}^{(1)}(X_0)) - v_{t, t}^{(2)}(\Psi_{0, t}^{(1)}(X_0))\|^2 \, \dd t
    \\&\quad
    + \smin^{-4}(1 + L_f^2) \integral{0}^{1 - t_0} \E \|\Psi_{0, t}^{(1)}(X_0) - \Psi_{0, t}^{(2)}(X_0)\|^2 \, \dd t .
\end{align}
Applying Lemma \ref{lem:heun_error_bound} with $M \asymp \smin^{-2}$ and $M_f \asymp \smin^{-2}\sqrt{D}$ together with Lemma \ref{lem:b_splines_basic_prop}, we deduce that
\begin{align}
    \label{eq:psi_j_global}
    \|\Psi_{0, t}^{(j)}(y)\| \lesssim \max_{0 \leq k \leq N} \|Z_k^{(j)}\|
    \lesssim \sqrt{D} \exp\{\cO(\smin^{-2})\} (1 + \|y\|),
    \quad \text{for all } y \in \R^D,
    \; j \in \{1, 2\} .
\end{align}
Furthermore, by Lemma \ref{lem:heun_error_bound}, we have that, for all $X_0$ satisfying
\begin{align}
\label{eq:covering_X_0_cond}
    \|X_0\| \lesssim \sqrt{D} \exp\{\cO(\smin^{-2})\} \left( \log(1 / \eps) + \log(D\smin^{-2}) \right),
\end{align}
it follows that $\|\Psi^{(j)}_{0,t}(X_0)\|\leq R$ for every $0 \leq t \leq 1 - t_0$ and $j \in \{1, 2\}$.
Then, using \eqref{eq:cj_proximity}, we obtain that
\begin{align*}
    &\integral{0}^{1 - t_0} \E \|v_{t, t}^{(1)}(\Psi_{0, t}^{(1)}(X_0)) - v_{t, t}^{(2)}(\Psi_{0, t}^{(1)}(X_0))\|^2 \, \dd t
    \lesssim \smin^{-8}(\sigma^{(1)} - \sigma^{(2)})^2 \integral{0}^{1 - t_0} \E\|\Psi_{0, t}^{(1)}(X_0)\|^2 \, \dd t
    \\&\quad
    + D\smin^{-8} (\sigma^{(1)} - \sigma^{(2)})^2
    + \smin^{-4} \|f^{(1)} - f^{(2)}\|^2_{L^\infty(\cB(0, R) \times [0, 1 - t_0])}
    \\&\quad
    + \smin^{-4} D \cdot \p(\|X_0\| \gtrsim \sqrt{D}\exp\{\cO(\smin^{-2})\}\log(1 / \eps)) .
\end{align*}
Using the fact that
\begin{align*}
    \integral{0}^{1 - t_0} \E \| \Psi_{0, t}^{(1)}(X_0) - \Psi_{0, t}^{(2)}(X_0) \|^2 \, \dd t
    \lesssim \integral{0}^{1 - t_0} \E \max_{0 \leq k \leq N} \|Z_k^{(1)} - Z_k^{(2)}\|^2 \, \dd t ,
\end{align*}
we conclude from \eqref{eq:flow_der_prox}, \eqref{eq:r_second_term_prox}, and \eqref{eq:psi_j_global} that
\begin{align}
\label{eq:reg_prox_diff_sq}
    \notag
    &\integral{0}^{1 - t_0} \E \| \partial_t[\Psi_{0,t}^{(1)}(X_0) - \Psi^{(2)}_{0, t}(X_0)]  - v_{t, t}^{(1)}(\Psi_{0, t}^{(1)}(X_0)) + v_{t, t}^{(2)}(\Psi_{0, t}^{(2)}(X_0))\|^2  \, \dd t
    \\&
    \lesssim (h^{-2} + \smin^{-4}(1 + L_f^2)) \integral{0}^{1 - t_0} \E \max_{0 \leq k \leq N} \|Z_k^{(1)} - Z_k^{(2)}\|^2 \, \dd t
    \\&
    + \smin^{-8} (\sigma^{(1)} - \sigma^{(2)})^2 D^2 \exp\{\cO(\smin^{-2})\}
    \notag
    + \smin^{-4} \|f^{(1)} - f^{(2)}\|^2_{L^\infty(\cB(0, R) \times [0, 1 - t_0])}
    + \eps^{4\beta} .
\end{align}
Similarly, from \eqref{eq:cj_proximity} we have that
\begin{align*}
    &\integral{0}^{1 - t_0} \E \| \partial_t[\Psi_{0,t}^{(1)}(X_0) + \Psi^{(2)}_{0, t}(X_0)]  - v_{t, t}^{(1)}(\Psi_{0, t}^{(1)}(X_0)) - v_{t, t}^{(2)}(\Psi_{0, t}^{(2)}(X_0))\|^2 \, \dd t
    \\&
    \lesssim \sum_{j = 1}^2\integral{0}^{1 - t_0} \E \left[\|\partial_t \Psi_{0, t}^{(j)}(X_0)\|^2 + \smin^{-4} \|\Psi_{0, t}^{(j)}(X_0)\|^2 + \smin^{-4} D \right] \, \dd t .
\end{align*}
Proceeding similarly to \eqref{eq:psi_j_global}, we deduce that, for all $j \in \{1, 2\}$,
\begin{align*}
    \|\partial_t \Psi_{0, t}^{(j)}(y)\|
    \lesssim h^{-1} \max_{0 \leq k \leq N}\|Z_k^{(j)}\|
    \lesssim h^{-1} \sqrt{D} \exp\{\cO(\smin^{-2})\} (1 + \|y\|),
    \quad \text{for all } y \in \R^D .
\end{align*}
This bound together with \eqref{eq:psi_j_global} implies that
\begin{align*}
    &\integral{0}^{1 - t_0} \E \| \partial_t[\Psi_{0,t}^{(1)}(X_0) + \Psi^{(2)}_{0, t}(X_0)]  - v_{t, t}^{(1)}(\Psi_{0, t}^{(1)}(X_0)) - v_{t, t}^{(2)}(\Psi_{0, t}^{(2)}(X_0))\|^2 \, \dd t
    \\&
    \lesssim h^{-2} D^2\exp\{\cO(\smin^{-2})\} .
\end{align*}
Using the above bound together with \eqref{eq:r_prox_cs} and \eqref{eq:reg_prox_diff_sq}, we conclude that
\begin{align}
\label{eq:r_v1_r_v2_almost_final}
    &|R(v^{(1)}) - R(v^{(2)})|
    \notag
    \lesssim h^{-1}D \exp\{\cO(\smin^{-2})\} \Bigg\{ (h^{-1} + L_f) \left[\integral{0}^{1 - t_0}\E \max_{0 \leq k \leq N} \|Z_k^{(1)} - Z_k^{(2)} \|^2 \, \dd t\right]^{1/2}
    \\&
    + D \smin^{-4}\exp\{\cO(\smin^{-2})\} |\sigma^{(1)} - \sigma^{(2)}|
    + \smin^{-2} \|f^{(1)} - f^{(2)}\|_{L^\infty(\cB(0, R) \times [0, 1 - t_0])} + \eps^{2\beta}
    \Bigg\} .
\end{align}

\noindent
\textbf{The distance between $Z^{(1)}_k$ and $Z^{(2)}_k$.}\quad
For $j\in\{1,2\}$, write
\[
 b_k^{(j)}(z)=c_1(\tau_k,\sigma^{(j)})z+c_2(\tau_k,\sigma^{(j)})f_k^{(j)}(z),\qquad
 e_k=\|Z_k^{(1)}-Z_k^{(2)}\|.
\]
Let $\Delta_\sigma=|\sigma^{(1)}-\sigma^{(2)}|$ and
$\Delta_f=\max_k\|f_k^{(1)}-f_k^{(2)}\|_{L^\infty(\cB(0,R))}$.
On $\cB(0,R)$ the velocity fields $b_k^{(j)}$ have Lipschitz constant $\ell \asymp \smin^{-4}$, and \eqref{eq:cj_proximity} yields
\[
 \|b_k^{(1)}(z)-b_k^{(2)}(z)\|
 \lesssim \smin^{-4}\left( (\|z\|+\sqrt D)\Delta_\sigma+\Delta_f \right).
\]
For convenience, let $\eta(z)$ denote the right-hand side of this estimate.
For initial points satisfying \eqref{eq:covering_X_0_cond}, choose the radius constants so that all sequences lie in $\cB(0,R)$.
We have that
\[
 \|Y_{k+1}^{(1)}-Y_{k+1}^{(2)}\|\leq(1 + h\ell)e_k+h\eta(Z_k^{(1)}).
\]
The corrected Heun update then implies
\begin{align*}
 e_{k+1}&\leq(1+h\ell+h^2\ell^2/2)e_k
 +\frac h2(1+h\ell)\eta(Z_k^{(1)})+\frac h2\eta(Y_{k+1}^{(1)})\\
 &\leq(1+3h\ell/2)e_k+h\eta(Z_k^{(1)})+\frac h2\eta(Y_{k+1}^{(1)}),
\end{align*}
where $h\ell \leq 1$.
The growth part of Lemma \ref{lem:heun_error_bound} bounds both terms by
$C\sqrt D\exp\{C\smin^{-2}\}(1+\|X_0\|)$.
Consequently, discrete Gr\"onwall lemma, applied with coefficient $1+3h\ell/2$ and $e_0=0$, yields
\begin{align*}
 \max_{0\leq k\leq N}e_k
 \lesssim\sqrt D\exp\{\cO(\smin^{-4})\}\Delta_\sigma(1+\|X_0\|)
 +\exp\{\cO(\smin^{-4})\}\Delta_f.
\end{align*}
Using Lemma \ref{lem:heun_error_bound}, we obtain that
\begin{align*}
    &\integral{0}^{1 - t_0}\E\max_{0 \leq k \leq N} \|Z_k^{(1)} - Z_k^{(2)}\|^2 \, \dd t
    \lesssim D \exp\{\cO(\smin^{-4})\} (\sigma^{(1)} - \sigma^{(2)})^2 R^2
    \\&\quad
    + \exp\{\cO(\smin^{-4})\} \max_{0 \leq l \leq N}\|f_l^{(1)} - f_l^{(2)}\|_{L^\infty(\cB(0, R))}^2
    \\&\quad
    + D\exp\{\cO(\smin^{-2})\} \integral{0}^{1 - t_0} \E[(1 + \|X_0\|^2)\1(\|X_0\| \gtrsim \sqrt{D} \exp\{\cO(\smin^{-2})\} \log(1 / \eps) )] \, \dd t .
\end{align*}
Gaussian concentration bound in Lemma \ref{lem:gaus_conc_norm} in conjunction with the Cauchy-Schwarz inequality implies that
\begin{align*}
    &\integral{0}^{1 - t_0}\E\max_{0 \leq k \leq N} \|Z_k^{(1)} - Z_k^{(2)}\|^2 \, \dd t
    \lesssim D \exp\{\cO(\smin^{-4})\} (\sigma^{(1)} - \sigma^{(2)})^2 R^2
    \\&\quad
    + \exp\{\cO(\smin^{-4})\} \max_{0 \leq l \leq N}\|f_l^{(1)} - f_l^{(2)}\|_{L^\infty(\cB(0, R))}^2
    + \eps^{4\beta} ,
\end{align*}
Substituting the bound for $R$ from the premise of the lemma, we obtain that
\begin{align*}
    \left(\integral{0}^{1 - t_0}\E\max_{0 \leq k \leq N} \|Z_k^{(1)} - Z_k^{(2)}\|^2 \, \dd t \right)^{1/2}
    &\lesssim D^{3/2} \log(1 / \eps) \exp\{\cO(\smin^{-4})\} |\sigma^{(1)} - \sigma^{(2)}| + \eps^{2\beta}
    \\&\quad
    + \exp\{\cO(\smin^{-4})\} \max_{0 \leq l \leq N} \| f^{(1)}_l - f^{(2)}_l \|_{L^\infty(\cB(0, R))} .
\end{align*}
Therefore, using \eqref{eq:r_v1_r_v2_almost_final} and the expression for $L_f$ given by the initial conditions of the lemma, we conclude that
\begin{align}
\label{eq:cov_reg_prox_aux}
    &|R(v^{(1)}) - R(v^{(2)})|
    \lesssim h^{-2}\exp{\cO(\smin^{-4})} (D \log(1 / \eps))^{\cO(1)}
    \\&
    \notag
    \cdot
    \left( |\sigma^{(1)} - \sigma^{(2)}| + \|f^{(1)} - f^{(2)}\|_{L^\infty(\cB(0, R) \times [0, 1 - t_0])} + \max_{0 \leq l \leq N} \|f^{(1)}_l - f^{(2)}_l\|_{L^\infty(\cB(0, R))} + \eps^{2\beta} \right) .
\end{align}

\noindent
\textbf{Covering number evaluation.}\quad
The following result evaluates the covering number of the class of feedforward neural networks with GELU activation function.

\begin{Lem}[\cite{yakovlev2025implicit}, Lemma B.1, adapted]
\label{lem:cov_nn_compact}
For every $R > 0$ and every
\begin{align*}
    0 < \tau \leq \sup \{\|f - g\|_{L^\infty([-R, R]^D)} \, : \, f, g \in \NN(L, W, S, B) \},
\end{align*}
the covering number satisfies the bound
\begin{align*}
    \log \cN(\tau, \NN(L, W, S, B), \|\cdot\|_{L^\infty([-R, R]^D)})
    \lesssim SL \log(\tau^{-1}L(B \vee 1)(\|W\|_\infty + 1)(R \vee 1)) .
\end{align*}

\end{Lem}

Let $\cF$ be the subclass of $\NN(L,W,S,B)$ satisfying the two constraints and the global Lipschitz bound in Definition \ref{def:av_vel_class}.
Use the compact domain $\cB(0,R+\tilde R+R_{\rm data})$.
Apply Lemma \ref{lem:cov_nn_compact} at half the desired radius.
For each covering ball that meets $\cF$, choose a representative from its intersection with $\cF$.
This yields a net $\cF_{\eps,k}\subseteq\cF$ of radius $\sqrt D\eps^{2\beta}$, up to a fixed constant, and
\begin{align}
\label{eq:cover_f_eps}
 \log|\cF_{\eps,k}|
 &\lesssim SL\log\bigl(\eps^{-1}L(B\vee1)(\|W\|_\infty+1)D\bigr)
 +SL\smin^{-2}+SL\log(e\vee\log(n/\delta)).
\end{align}
Here, the parameter $R$ is taken from the initial conditions of the lemma.
Now let $\cG_\eps$ be an $\eps^{2\beta}$-net of the interval $[\smin, 1]$.
Evidently, $|\cG_\eps| \lesssim \eps^{-2\beta}$.
Choose the scalar net inside $[\smin,1]$ and form $\cV_\eps$ from its product with the constrained subnetwork nets.
Every resulting candidate lies in $\cV$.
From \eqref{eq:cover_f_eps} it follows that
\begin{align*}
    \log |\cV_\eps|
    &\leq \log |\cG_\eps| + \sum_{k = 0}^N\log|\cF_{k, \eps}|
    \\&
    \lesssim \smin^{-2}(N + 1)SL \log(\eps^{-1} L (B \vee 1)(\|W\|_\infty + 1)D) + (N+1)SL\log(e \vee \log(n / \delta)) .
\end{align*}
By Lemma \ref{lem:b_splines_basic_prop}, we have that
\begin{align}
\label{eq:cover_f_fk_red}
    \|f^{(1)} - f^{(2)}\|_{L^\infty(\cB(0, R + \tilde R + R_\text{data}) \times [0, 1 - t_0])}
    \lesssim \max_{0 \leq k \leq N} \|f_k^{(1)} - f_k^{(2)}\|_{L^\infty(\cB(0, R + \tilde R + R_\text{data}))} .
\end{align}
Combining this with \eqref{eq:cov_fm_ball_aux}, we obtain that, with probability at least $(1 - \delta)$,
\begin{align*}
    &\sup_{v \in \cV}\inf_{v' \in \cV_\eps} \|\ell_{FM}(v, \cdot) - \ell_{FM}(v', \cdot)\|_{L^\infty(\cB(0, R_\text{data}))}
    \lesssim (D\smin^{-2})^{\cO(1)} \left(\log(n / \delta) + \log(1 / \eps)\right)\lambda(t_0) \eps^{2\beta} .
\end{align*}
Similarly, \eqref{eq:cov_fm_exp_prox_aux} yields that
\begin{align*}
    &\sup_{v \in \cV} \inf_{v' \in \cV_\eps} |\E\ell_{FM}(v, X_1) - \E\ell_{FM}(v', X_1)|
    \lesssim (D\smin^{-2})^{\cO(1)} \left(\log(n / \delta) + \log(1 / \eps)\right)\lambda(t_0) \eps^{2\beta} .
\end{align*}
Finally, from \eqref{eq:cov_reg_prox_aux} and \eqref{eq:cover_f_fk_red} we deduce that
\begin{align*}
    \sup_{v \in \cV}\inf_{v' \in \cV_\eps}|R(v) - R(v')|
    \lesssim h^{-2} \exp\{\cO(\smin^{-4})\} (D \log(1 / \eps))^{\cO(1)} \eps^{2\beta} .
\end{align*}
Since $N+1\lesssim h^{-1}$, the entropy bound follows, and \eqref{eq:cov_triang_first} proves the lemma.

\myendproof

\section{Proof of Theorem \ref{th:distr_estim}}
\label{sec:th_distr_estim_proof}

Let $\Psi_{0, t}(y) = y + t \cdot v_{0, t}(y)$, where $v \in \cV(L, W, S, B, N, R, L_f)$ is arbitrary.
It holds that
\begin{align}
\label{eq:distr_w2_triang}
    W_2(\mathrm{Law}(\Psi_{0, 1 - t_0}(X_0)), \sfp^*)
    \leq W_2(\mathrm{Law}(\Psi_{0, 1 - t_0}(X_0)), \sfp_{1 - t_0}^*)
    + W_2(\sfp_{1 - t_0}^*, \sfp^*).
\end{align}
We now bound each term individually.
For every $t \in [0, 1 - t_0]$ define $E_t = \E \| \Psi^*_{0, t}(X_0) - \Psi_{0, t}(X_0) \|^2$.
Then, it follows that
\begin{align*}
    \frac{\dd E_t}{\dd t}
    &= 2 \E \langle \Psi^*_{0, t}(X_0) - \Psi_{0, t}(X_0), b^*(\Psi^*_{0, t}(X_0), t) - \partial_t \Psi_{0, t}(X_0) \rangle
    \\&
    = 2\E \langle \Psi^*_{0, t}(X_0) - \Psi_{0, t}(X_0), b^*(\Psi^*_{0, t}(X_0), t) - v_{t, t}(\Psi^*_{0, t}(X_0)) \rangle
    \\&\quad
    + 2 \E \langle \Psi^*_{0, t}(X_0) - \Psi_{0, t}(X_0), v_{t, t}(\Psi^*_{0, t}(X_0)) - v_{t, t}(\Psi_{0, t}(X_0)) \rangle
    \\&\quad
    + 2 \E \langle \Psi^*_{0, t}(X_0) - \Psi_{0, t}(X_0), v_{t, t}(\Psi_{0, t}(X_0)) - \partial_t \Psi_{0, t}(X_0) \rangle .
\end{align*}
Applying the Cauchy-Schwarz inequality and the Young inequality, we obtain
\begin{align}
\label{eq:distr_et_der}
    \notag
    &\frac{\dd E_t}{\dd t}
    \leq 2 E_t + \E \|b^*(\Psi^*_{0, t}(X_0), t) - v_{t, t}(\Psi^*_{0, t}(X_0))\|^2
    + \E \|v_{t, t}(\Psi_{0, t}(X_0)) - \partial_t \Psi_{0, t}(X_0)\|^2
    \\&
    + 2 \E \left[ \|\Psi^*_{0, t}(X_0) - \Psi_{0, t}(X_0)\| \cdot \|v_{t, t}(\Psi^*_{0, t}(X_0)) - v_{t, t}(\Psi_{0, t}(X_0))\| \right] .
\end{align}
Using Lemma \ref{lem:heun_error_bound} with $M \asymp \smin^{-2}$ and $M_f \asymp \smin^{-2} \sqrt{D}$ and Lemma \ref{lem:b_splines_basic_prop}, we have that, for all $X_0$ satisfying
\begin{align*}
    \|X_0\| \lesssim \sqrt{D} \exp\{\cO(\smin^{-2})\} \log(n D\smin^{-2}),
\end{align*}
it follows that $\|\Psi^*_{0, t}(X_0)\| \vee \|\Psi_{0, t}(X_0)\| \leq R$ for every $0 \leq t \leq 1 - t_0$.
Furthermore, we also have that
\begin{align}
\label{eq:distr_psi_bounds}
    \|\Psi_{0, t}(y)\| \vee \|\Psi_{0, t}^*(y)\|
    \lesssim \exp\{\cO(\smin^{-2})\} \sqrt{D} (1 + \|y\|),
    \quad \text{for all } y \in \R^D .
\end{align}
For $z,z'\in\cB(0,R)$, the linear term and the spline coefficients yield
\begin{align*}
 \|v_{t,t}(z)-v_{t,t}(z')\|
 &\lesssim\smin^{-2}\|z-z'\|+\smin^{-2}\max_k\|f_k(z)-f_k(z')\|\\
 &\lesssim\smin^{-4}\|z-z'\|.
\end{align*}
Here the fixed constants account for the extrapolated boundary coefficients, as in Lemma \ref{lem:b_splines_basic_prop}.
Apply this bound with $z=\Psi_{0,t}(X_0)$ and $z'=\Psi^*_{0,t}(X_0)$ on the event above.
Hence, using the fact that $v_{t, t}$ is globally Lipschitz, we obtain that
\begin{align*}
    &\E \left[ \|\Psi^*_{0, t}(X_0) - \Psi_{0, t}(X_0)\| \cdot \|v_{t, t}(\Psi^*_{0, t}(X_0)) - v_{t, t}(\Psi_{0, t}(X_0))\| \right]
    \\&
    \lesssim \smin^{-4} \E \|\Psi_{0, t}(X_0) - \Psi_{0, t}^*(X_0)\|^2
    \\&\quad
    + \smin^{-2}(1 + L_f) \E \left[ \|\Psi^*_{0, t}(X_0) - \Psi_{0, t}(X_0)\|^2 \cdot \1[\|X_0\| \gtrsim \sqrt{D} \exp\{\cO(\smin^{-2})\} \log(n D\smin^{-2}) ] \right] .
\end{align*}
By the Cauchy-Schwarz inequality, the bound  from \eqref{eq:distr_psi_bounds}, and the fact that $X_0 \sim \cN(0, I_D)$ we obtain
\begin{align*}
    &\E \left[ \|\Psi^*_{0, t}(X_0) - \Psi_{0, t}(X_0)\| \cdot \|v_{t, t}(\Psi^*_{0, t}(X_0)) - v_{t, t}(\Psi_{0, t}(X_0))\| \right]
    \\&
    \lesssim \smin^{-4} \E \|\Psi_{0, t}(X_0) - \Psi_{0, t}^*(X_0)\|^2
    \\&\quad
    + \smin^{-2}(1 + L_f) \exp\{\cO(\smin^{-2})\} D^{\cO(1)} \p^{1/2}(\|X_0\| \gtrsim \sqrt{D} \exp\{\cO(\smin^{-2})\} \log(n D\smin^{-2}) ) .
\end{align*}
Substituting the expression for $L_f$ and using the Gaussian concentration bound from Lemma \ref{lem:gaus_conc_norm}, we conclude that
\begin{align*}
    \E \left[ \|\Psi^*_{0, t}(X_0) - \Psi_{0, t}(X_0)\| \cdot \|v_{t, t}(\Psi^*_{0, t}(X_0)) - v_{t, t}(\Psi_{0, t}(X_0))\| \right]
    \lesssim \smin^{-4} E_t + n^{-1}.
\end{align*}
Substituting this bound into \eqref{eq:distr_et_der} and taking into account that $\mathrm{Law}(\Psi_{0, t}^*(X_0)) = \sfp_t^*$, we obtain
\begin{align*}
    \frac{\dd E_t}{\dd t}
    \lesssim \smin^{-4} E_t + n^{-1}
    + \E \|b^*(X_t, t) - v_{t, t}(X_t)\|^2
    + \E \|v_{t, t}(\Psi_{0, t}(X_0)) - \partial_t \Psi_{0, t}(X_0)\|^2 .
\end{align*}
By the Grönwall lemma, we have that
\begin{align*}
    E_{1 - t_0}
    &\lesssim \exp\{\cO(\smin^{-4})\}\left( \integral{0}^{1 - t_0} \E\|v_{t, t}(X_t) - b^*(X_t, t)\|^2 \, \dd t + R(v) + n^{-1} \right)
    \\&
    \lesssim \exp\{\cO(\smin^{-4})\} \left( \E[\ell_{SD}(v, X_1) - \ell_{SD}(v^*, X_1)] + n^{-1} \right),
\end{align*}
where the last inequality uses Proposition \ref{prop:fm_loss_dif}.
Therefore, observing that $W_2^2(\mathrm{Law}(\Psi_{0, 1 - t_0}(X_0)), \sfp_{1 - t_0}^*) \leq E_{1 - t_0}$ and using the generalization bound established in Theorem \ref{th:estim_main}, we conclude that
\begin{align}
\label{eq:w2_sq_stopped}
    \notag
    &W_2^2(\mathrm{Law}(\widehat{\Psi}_{0, 1 - t_0}(X_0)), \sfp^*_{1 - t_0})
    \\&
    \lesssim (D \log n)^{\cO(1)} e^{\cO(\smin^{-4})}
    \lambda(t_0)
    \left(\frac{\lambda(t_0)}{n} \right)^{\frac{2\beta}{5\beta + 2d}} \log(\lambda(t_0)) \log(4 / \delta),
\end{align}
on an event of probability at least $(1 - \delta)$.
As for the second term in \eqref{eq:distr_w2_triang}, by Assumption \ref{asn:relax_man}, we have that
\begin{align*}
    W_2(\sfp^*_{1 - t_0}, \sfp^*)
    \leq \E^{1/2} \|\alpha_{1 - t_0}X_0 + \beta_{1 - t_0}X_1 - X_1\|^2
    \lesssim \sqrt{D}(\alpha_{1 - t_0} \vee (1 - \beta_{1 - t_0})) .
\end{align*}
Combining this bound with \eqref{eq:distr_w2_triang} and \eqref{eq:w2_sq_stopped} yields the desired result.

\myendproof

\section{Tools from probability theory}

\begin{Lem}[Gaussian concentration of the norm \cite{vershynin2018high}, Theorem 3.1.1, adapted]
\label{lem:gaus_conc_norm}
Let $Z \sim \cN(0, I_D)$ be a standard Gaussian random vector in $\R^D$.
There exist absolute constants $C, c > 0$ such that,
\begin{align*}
    \p\left(\big|\|Z\| - \sqrt{D}\big| \geq t \right) \leq C \exp(-c t^2),
    \quad \text{for all $t \geq 0$}.
\end{align*}
\end{Lem}

\end{document}